\documentclass[11pt]{article}
\usepackage{graphicx}
\usepackage{color}
\usepackage{epsfig}
\usepackage{amsmath}
\usepackage{amsthm}
\usepackage{amsfonts}
\usepackage{amssymb}
\usepackage[english]{babel}

\newtheorem{theorem}{Theorem}
\newtheorem{corollary}[theorem]{Corollary}
\newtheorem{lemma}[theorem]{Lemma}
\newtheorem{proposition}[theorem]{Proposition}
\theoremstyle{definition}
\newtheorem{definition}[theorem]{Definition}
\newtheorem{remark}[theorem]{Remark}
\newtheorem{example}[theorem]{Example}

\newcommand{\R}{\mathbb{R}}

\newcommand{\argmin}{\mbox{\rm argmin\,}}

\DeclareMathOperator{\Dist}{Dist}

\begin{document}

%%%%%%%%%%%%%%%%%%%%%%%%%%%%%%%%%%%%%%%%%%%%%%%%%%%%%
%%%%%                                           %%%%%
%%%%%    Characterizations of Lojasiewicz       %%%%%
%%%%%       inequalities and applications       %%%%%
%%%%%                                           %%%%%
%%%%% J. Bolte, A. Daniilidis, O. Ley, L. Mazet %%%%%
%%%%%                                           %%%%%
%%%%%            February 6, 2008               %%%%%
%%%%%          (version soumise  - tours)       %%%%%
%%%%%%%%%%%%%%%%%%%%%%%%%%%%%%%%%%%%%%%%%%%%%%%%%%%%%

\begin{center}
{\LARGE Characterizations of \L ojasiewicz inequalities and applications
\\ \medskip
} \vspace{0.2cm}

\vspace{0.8cm}

\textbf{J\'{e}r\^{o}me BOLTE},\textbf{\ Aris DANIILIDIS},\textbf{\ Olivier LEY
\& Laurent MAZET}

\bigskip
\end{center}

\noindent\textbf{Abstract }The classical \L ojasiewicz inequality
and its extensions for partial differential equation problems
(Simon) and to o-minimal structures (Kurdyka) have a considerable
impact on the analysis of gradient-like methods and related
problems: minimization methods, complexity theory, asymptotic
analysis of dissipative partial differential equations, tame
geometry. This paper provides alternative characterizations of
this type of inequalities for nonsmooth lower semicontinuous
functions defined on a metric or a real Hilbert space. In a metric
context, we show that a generalized form of the \L ojasiewicz inequa\-li\-ty
(hereby called the Kurdyka-\L ojasiewicz inequa\-li\-ty) relates to metric
regularity and to the Lipschitz continuity of the sublevel
mapping, yielding applications to discrete methods (strong
convergence of the proximal algorithm). In a Hilbert setting we
further establish that asymptotic properties of the semiflow
generated by $-\partial f$ are strongly linked to this inequality.
This is done by introducing the notion of a piecewise subgradient
curve: such curves have uniformly bounded lengths if and only if
the Kurdyka-\L ojasiewicz inequality is satisfied. Further
characterizations in terms of \textit{talweg} lines ---a concept
linked to the location of the less steepest points at the level
sets of $f$--- and integrability conditions are given. In the
convex case these results are significantly reinforced, allowing
in particular to establish the asymptotic equivalence of discrete
gradient methods and continuous gradient curves. On the other
hand, a counterexample of a convex $C^{2}$ function in
$\mathbb{R}^{2}$ is constructed to illustrate the fact that,
contrary to our intuition, and unless a specific growth condition is
satisfied, convex functions may fail to fulfill the Kurdyka-\L
ojasiewicz inequality.

\bigskip


\bigskip

\textit{Secondary }03C64, 37N40, 49J52, 65K10.

\bigskip

\noindent\textbf{Acknowledgement }The first two authors
acknowledge support of the ANR grant ANR-05-BLAN-0248-01 (France).
The second author acknowledge support of the MEC grant
MTM2005-08572-C03-03 (Spain). During the preparation of this work,
several research visits of the co-authors have been
realized, respectively to the CRM (Mathematical Research Center in
Barcelona), the University Autonomous of Barcelona, the University
of Paris~6 and the University of Tours. In each case the concerned
author wishes to acknowledge their hosts for hospitality.

\newpage
\tableofcontents
\newpage

\section{Introduction}

The \L ojasiewicz inequality is a powerful tool to analyze
convergence of gradient-like methods and related problems. Roughly
speaking, this inequality is satisfied by a $C^{1}$ function $f$,
if for some $\theta\in\lbrack\frac {1}{2},1)$ the quantity
\[
|f-f(\bar{x})|^{\theta}\,\Vert\nabla f\Vert^{-1}
\]
remains bounded away from zero around any (possibly critical)
point $\bar {x}$. This result is named after S.~\L ojasiewicz
\cite{Loja63}, who was the first to establish its validity for the
classes of real--analytic and $C^{1}$ subanalytic functions. At
the same time, it has been known that the \L ojasiewicz inequality
would fail for $C^{\infty}$ functions in general (see the
classical example of the function $x\longmapsto \exp(-1/x^{2})$,
if $x\neq0$ and $0$, if $x=0$ around the point $\bar{x}=0$).

\smallskip

A generalized form of this inequality has been introduced by
K.~Kurdyka in \cite{Kurdyka98}. In the framework of a $C^{1}$
function $f$ defined on a real Hilbert space
$[H,\langle\cdot,\cdot\rangle]$, and assuming for simplicity that
$\bar{f}=0$ is a critical value, this generalized inequality (that
we hereby call the Kurdyka--\L ojasiewicz inequality, or in short,
the K\L--inequality) states that
\begin{equation}
||\nabla(\varphi\circ f)(x)||\geq1, \label{Loja0}
\end{equation}
for some continuous function $\varphi:[0,r)\rightarrow\mathbb{R}$,
$C^{1}$ on $(0,r)$ with $\varphi^{\prime}>0$ and all $x$ in
$[0<f<r]:=\{y\in H:0<f(y)<r\}$. The class of such functions
$\varphi$ will be further denoted by $\mathcal{K}(0,\bar{r})$, see
(\ref{Kr}). Note that the \L ojasiewicz inequality corresponds to
the case $\varphi(t)=t^{1-\theta}$.

\smallskip

In finite-dimensional spaces it has been shown in \cite{Kurdyka98}
that (\ref{Loja0}) is satisfied by a much larger class of
functions, namely, by those that are definable in an o-minimal
structure \cite{Coste99}, or even more generally by functions
belonging to analytic-geometric categories \cite{Dries-Miller96}.
In the meantime the original \L ojasiewicz result was used to
derive new results in the asymptotic analysis of nonlinear heat
equations \cite{Sim83} and damped wave equations \cite{Har01}.
Many results related to partial differential equations followed,
see the monograph of Huang \cite{Huang} for an insight. Other
fields of application of (\ref{Loja0}) are nonconvex optimization
and nonsmooth analysis. This was one of the motivations for the
nonsmooth K\L--inequalities developed in \cite{BDLS2006,BDLS2007}.
Due to its considerable impact on several field of applied
mathematics: minimization and algorithms
\cite{Absil,AttBol,BDLS2006, Lageman07}, asymptotic theory of
differential inclusions \cite{NisQui}, neural networks \cite{Quinc}, complexity theory
\cite{NestPol} (see \cite[Definition~3]{NestPol} where functions
satisfying a K\L--type inequality are called gradient dominated
functions), partial differential equations
\cite{Sim83,Har01,Huang}, we hereby tackle the problem of
characterizing such inequalities in an nonsmooth
infinite-dimensional setting and provide further clarification in
several application aspects. Our framework is rather broad
(infinite dimensions, nonsmooth functions), nevertheless, to the
best of our knowledge, most of the present results are also new in
a smooth finite-dimensional framework: readers who feel unfamiliar
with notions of nonsmooth and variational analysis may, at a first
stage, consider that all functions involved 
are differentiable and replace subdifferentials by usual derivatives
and subgradient systems by smooth ones.

\smallskip

A first part of this work (Section~2) is devoted to the analysis
of metric versions of the K\L--inequality. The underlying space
$H$ is only assumed to be a complete metric space (without any
linear structure), the function
$f:H\rightarrow\mathbb{R}\cup\{+\infty\}$ is lower semicontinuous
and possibly real-extended valued and the notion of a gradient is
replaced by the variational notion of a strong-slope
\cite{DMT1980,AzeCor2004}. Indeed, introducing the multivalued
mapping $F(x)=[f(x),+\infty)$ (whose graph is the epigraph of
$f$), the K\L--inequality (\ref{Loja0}) appears to be equivalent
to the metric regularity of $F:H\rightrightarrows\mathbb{R}$ on an
adequate set, where $\mathbb{R}$ is endowed with the metric
$d_{\varphi}(r,s)=|\varphi (r)-\varphi(s)|$. This fact is strongly
connected to famous classical results in this area (see
\cite{DLR02,mordukhovich93,Ioffe,penot89} for example) and in
particular to the notion of $\rho$-metric regularity introduced in
\cite{Ioffe} by A.~Ioffe. The particularity of our result is due
to the fact that $F$ takes its values in a totally ordered set
which is not the case in the general theory. Using results on
global error-bounds of Az\'e-Corvellec \cite{AzeCor2004} and Zorn's lemma, we
establish indeed that some global forms of the K\L-inequality and
metric regularity are both equivalent to the ``Lipschitz
continuity'' of the sublevel mapping
\[
\left\{
\begin{array}
[c]{lll} \mathbb{R} & \rightrightarrows &  H\\ r & \mapsto &
[f\leq r]:=\{x\in H:f(x)\leq r\},
\end{array}
\right.
\]
where $(0,r)\subset(0,+\infty)$ is endowed with $d_{\varphi}$ and the
collection of subsets of $H$ with the ``Hausdorff distance''. As it is shown
in a section devoted to applications (Section 3.4), this reformulation is
particularly adapted for the analysis of proximal methods involving nonconvex
criteria: these results are in the line of \cite{combettes,AttBol}.

\smallskip

In the second part of this work (Section~3), $H$ is a proper real
Hilbert space and $f$ is assumed to be a semiconvex function,
\emph{i.e.} $f$ is the difference of a proper lower semicontinuous
convex function and a function proportional to the canonical
quadratic form. Although this assumption is not particularly
restrictive, it does not aim at full generality. Semiconvexity is used here to provide a convenient framework in which the
 formulation and the study of subdifferential evolution equations are simple and elegant
(\cite{ac99,DMT1985}). Using the Fr\'{e}chet subdifferential (see
Definition~\ref{Definition_Frechet}), the corresponding
subgradient dynamical system indeed reads
\begin{equation}
\left\{
\begin{array}
[c]{l} \dot{x}(t)+\partial f(x(t))\ni0,\mbox{ a.e. on
}(0,+\infty),\\ x(0)\in\mbox{\rm dom\,}f
\end{array}
\right.  \label{SD}
\end{equation}
where $x(\cdot)$ is an \emph{absolutely continuous curve} called
\emph{subgradient curve}. Relying on several works
\cite{DMT1985,Thibault,Brezis}, if $f$ is semiconvex, such curves
exist and are unique. The asymptotic properties of the semiflow
associated to this evolution equation are strongly connected to
the K\L-inequality. This can be made precise by introducing 
 the following notion: for $T\in (0,+\infty],$ a piecewise
absolutely continuous curve $\gamma:[0,T)\rightarrow H$ (with
countable pieces) is called a \emph{piecewise subgradient curve}
if $\gamma$ is a solution to (\ref{SD}) where in addition
$t\mapsto(f\circ \gamma)(t)$ nonincreasing (see
Definition~\ref{piecewise} for details). Consider all piecewise
subgradient curves lying in a ``K\L--neighborhood'', \textit{e.g.}
a slice of level sets. Under a compactness assumption and a
condition of Sard type (automatically satisfied in finite
dimensions if $f$ belongs to an o-minimal class), their lengths
are uniformly bounded if and only if $f$ satisfies the
K\L--inequality in its nonsmooth form (see \cite{BDLS2007}), that
is, for all $x\in [0<f<r]$,
\[
||\partial(\varphi\circ f)(x)||_{-}:=\inf\{||p||: p\in \partial(\varphi\circ f)\}\geq1,
\]
where $\varphi:(0,r)\rightarrow \R$  is $C^1$ function bounded from below such that $\varphi'>0$ (see (\ref{Kr})). A byproduct of this result
(through not an equivalent statement, as we show in Section~4.3
--see Remark~\ref{Remark38}~(c)) is the fact that bounded
subgradient curves have finite lengths and hence converge to a
generalized critical point.

\smallskip

Further characterizations are given involving several aspects among which, an
integrability condition in terms of the inverse function of the minimal
subgradient norm associated to each level set $[f=r]$ of $f,$ as well as
connections to the following \emph{talweg} selection problem:
Find a piecewise absolutely continuous curve $\theta:(0,r)\rightarrow H$ with
finite length such that
\[
\theta(r)\in\left\{  x\in\lbrack f=r]:
||\partial(\varphi\circ f)(x)||_{-} \leq
R\inf_{y\in\lbrack f=r]}||\partial(\varphi\circ f)(y)||_{-}
\right\}  ,\;\mbox
{ with }R>1.
\]
The curve $\theta$ is called a \emph{talweg}. Early connections between the
K\L-inequality and this old concept can be found in \cite{Kurdyka98}, and even
more clearly in \cite{Dacunto}. Indeed, under mild assumptions the existence
of such a selection curve $\theta$ characterizes the K\L-inequality. The proof
relies strongly on the property of the semiflow associated to $-\partial f$.
Recent developments of the metric theory of ``gradient'' curve (\cite{Amb})
open the way to a more general approach of these characterizations, and
hopefully to new applications in the line of \cite{Amb,DMT1980}.

\smallskip

The analysis of the convex case (that is, $f$ is a convex
function) in Section 4, reveals interesting phenomena. In this
case, the K\L-inequality, whenever true on a slice of level sets,
will be true on the whole space $H$ (globalization) and, in
addition, the involved function $\varphi$ can be taken to be
concave (Theorem~\ref{convex}). This is always the case if a
specific growth assumption near the set of minimizers of $f$ is
assumed. On the other hand, arbitrary convex functions do not
satisfy the K\L--inequality: this is a straightforward consequence
of a classical counterexample, due to J.-B.~Baillon
\cite{baillon78}, of the existence of a convex function $f$ in a
Hilbert space, having a subgradient curve which is not 
strongly converging to $0\in\arg\min f$.
However, surprisingly, even smooth finite-dimensional coercive
convex functions may fail to satisfy the K\L-inequality, and this
even in the case that the lengths of their gradient curves are
uniformly bounded. Indeed, using the above mentioned
characterizations and results from \cite{Torralba}, we construct a
counterexample of a $C^{2}$ convex function whose set of
minimizers is compact and has a nonempty interior (Section~4.3).

\smallskip

As another application we consider abstract \emph{explicit}
gradient schemes for convex functions with a Lipschitz continuous
gradient. A common belief is that the analysis of gradient curves
and their explicit discretization used in numerical optimization
are somehow disconnected problems. We hereby show that this is not always
the case, by establishing that the piecewise gradient iterations
are uniformly bounded if and only if the piecewise subgradient
curves are so. This aspect sheds further light on the
(theoretical) stability of convex gradient-like methods and the
interest of relating the K\L--inequality to the asymptotic study
of subgradient-type methods.

\bigskip

\noindent\textbf{Notation.} (Multivalued mappings) Let $X,Y$ be
two metric spaces and $F:X\rightrightarrows Y$ be a multivalued
mapping from $X$ to $Y.$ We denote by
\begin{equation}
\mathrm{Graph\,}F:=\{(x,y)\in X\times Y:y\in F(x)\} \label{aris3}
\end{equation}
the \emph{graph} of the multivalued mapping $F$ (subset of
$X\times Y$) and by
\begin{equation}
\mathrm{dom\,}F:=\{x\in X:\mathrm{\,}\exists y\in
Y,\mathrm{\,}(x,y)\in \mathrm{Graph\,}F\} \label{aris2}
\end{equation}
its \emph{domain} (subset of $X$).

\smallskip

\noindent(Single--valued functions) Given a function $f:X\longrightarrow
\mathbb{R}\cup\{+\infty\}$ we define its \emph{epigraph} by
\begin{equation}
\mathrm{epi\,}f:=\{(x,\beta)\in X\times\mathbb{R}:f(x)\leq\beta\}.
\label{aris5}
\end{equation}
We say that the function $f$ is \emph{proper} (respectively, \emph{lower
semicontinuous}) if the above set is nonempty (respectively, closed). Let us
recall that the domain of the function $f$ is defined by
\[
\mathrm{dom\,}f:=\{x\in X: f(x)<+\infty\}.
\]

\smallskip

\noindent(Level sets) Given $r_{1}\leq r_{2}$ in $[-\infty,+\infty]$ we set
\[
\lbrack r_{1}\leq f\leq r_{2}]:=\{x\in X:r_{1}\leq f(x)\leq r_{2}\}.
\]
When $r_{1}=r_{2}$ (respectively $r_{1}=-\infty$), the above set
will be simply denoted by $[f=r_{1}]$ (respectively $[f\leq
r_{2}]$).

\smallskip

\noindent(Strong slope) Let us recall from \cite{DMT1980} (see also
\cite{Ioffe}, \cite{AzeCor2004}) the notion of \emph{strong slope }defined for
every $x\in\mathrm{dom\,}f$ as follows:
\begin{equation}
|\nabla f|(x)=\underset{y\rightarrow
x}{\mathrm{\,}\lim\sup\mathrm{\,}} \frac{\left(  f(x)-f(y)\right)
^{+}}{d(x,y)}, \label{slope}
\end{equation}
where for every $a\in\mathbb{R}$ we set $a^{+}=\max\,\{a,0\}$.

If $[X,||\cdot||]$ is a Banach space with (topological) dual space
$[X^{*},||\cdot||_{*}]$ and $f$ is a $C^{1}$ finite-valued function then
\[
|\nabla f|(x)=||\nabla f(x)||_{*},
\]
for all $x$ in $X$, where $\nabla f(\cdot)$ is the differential map of $f$.

\smallskip

\noindent(Hausdorff distance) We define the \emph{distance} of a
point $x\in X$ to a subset $S$ of $X$ by
\[
\mathrm{dist\,}(x,S):=\mathrm{\,}\underset{y\in S}{\inf}\mathrm{\,}d(x,y),
\]
where $d$ denotes the distance on $X.$ The \emph{Hausdorff
distance} $\mathrm{Dist}(S_{1},S_{2})$ of two subsets $S_{1}$ and
$S_{2}$ of $X$ is given by
\begin{equation}
\mathrm{Dist}(S_{1},S_{2})\,:=\,\max\,\left\{  \underset{x\in
S_{1}}{\sup }\,\mathrm{dist\,}(x,S_{2}),\,\underset{x\in
S_{2}}{\sup}\,\mathrm{dist\,} (x,S_{1})\right\}
\,.\label{Hausdorff}
\end{equation} Let us denote by
$\mathcal{P}(X)$ the collection of all subsets of $X$. In general
$\mathrm{Dist}(\cdot,\cdot)$ can take infinite values and does not
define a distance on $\mathcal{P}(X)$. However if $K(X)$ denotes
the collection of nonempty compact subsets of $X$, then
$\mathrm{Dist}(\cdot ,\cdot)$ defines a proper notion of distance
on $K(X)$. In the sequel we deal with multivalued mappings
$F:X\rightrightarrows Y$ enjoying the following property
\[
\mbox{\rm Dist}\,(F(x),F(y))\leq k\;d(x,y)
\]
where $k$ is a positive constant. For simplicity such functions are called
Lipschitz continuous, although $[\mathcal{P}(Y),\,\mbox{\rm Dist}\,]$ is not a
metric space in general.

\smallskip

\noindent(Desingularization functions) Given
$\bar{r}\in(0,+\infty]$, we set
\begin{equation}
\mathcal{K}(0,\bar{r}):=\left\{  \phi\in C([0,\bar{r}))\cap
C^{1}(0,\bar {r}):\;\phi(0)=0,\;\mbox{ and
}\phi^{\prime}(r)>0,\forall r\in(0,\bar {r})\right\}  ,\label{Kr}
\end{equation}
where $C([0,\bar{r}])$ (respectively, $C^{1}(0,\bar{r})$) denotes the set of
continuous functions on $[0,\bar{r}]$ (respectively, $C^{1}$ functions on
$(0,\bar{r})$).

\smallskip

Finally throughout this work, $B(x,r)$ will stand for the usual
open ball of center $x$ and radius $r>0$ and $\bar{B}(x,r)$ will
denote its closure. If $H$ is a Hilbert space, its inner product
will be denoted by $\langle\cdot ,\cdot\rangle$ and the
corresponding norm by $||\cdot||.$

\section{K\L--inequality is a metric regularity condition}

\label{S:reg}

Let $X,Y$ be two \emph{complete} metric spaces, $F:X\rightrightarrows Y$ a
multivalued mapping and $(\bar{x},\bar{y})\in\mathrm{Graph\,}F.$ Let us recall
from \cite[Definition 1 (loc)]{Ioffe} the following definition.

\begin{definition}
[metric regularity of multifunctions]\label{DefinitionMetricRegF}
Let $k\in\lbrack0,+\infty)$. (i) The multivalued mapping $F$ is
called $k$-metrically regular at
$(\bar{x},\bar{y})\in\mbox{Graph}\;F$, if there exist
$\varepsilon,\delta>0$ such that for all $(x,y)\in B(\bar{x}
,\varepsilon)\times B(\bar{y},\delta)$ we have
\begin{equation}
\mathrm{dist\,}(x,F^{-1}(y))\mathrm{\,}\leq\mathrm{\,}k\mathrm{\,dist\,}
(y,F(x)). \label{aris1}
\end{equation}
(ii) Let $V$ be a nonempty subset of $X\times Y$. The multivalued mapping $F $
is called $k$-metrically regular on $V$, if $F$ is metrically regular at
$(\bar{x},\bar{y})$ for every $(\bar{x},\bar{y})\in\mathrm{Graph\,}F\cap V.$
\end{definition}

\subsection{Metric regularity and global error bounds}

The following theorem is an essential result: it will show that
Kurdyka-\L ojasiewicz inequality and metric regularity are equivalent concepts
(see Corollary~\ref{Corollary_KurMetric} and Remark~\ref{top}). The
equivalence [(ii)$\Leftrightarrow$(iii)] is due to Az\'e-Corvellec (see \cite[Theorem
2.1]{AzeCor2004}).

\begin{theorem}
\label{PropositionIC} Let $X$ be a complete metric space, $f:X\longrightarrow
\mathbb{R}\cup\{+\infty\}$ a proper lower semicontinuous function and
$r_{0}>0$. The following assertions are equivalent:\smallskip

\noindent(i) The multivalued mapping
\[
F:\left\{
\begin{array}
[c]{cll} X & \rightrightarrows & \mathbb{R}\\ x & \longmapsto &
[f(x),+\infty)
\end{array}
\right.
\]
is $k$-metrically regular on $[0<f<r_{0}]\times(0,r_{0})$ $;$\smallskip

\noindent(ii) For all $r\in(0,r_{0})$ and $x\in\lbrack0<f<r_{0}]$
\begin{equation}
\mathrm{dist\,}(x,[f\leq r])\mathrm{\,}\leq\mathrm{\,}k\mathrm{\,}
(f(x)-r)^{+}\text{ }; \label{aris6}
\end{equation}

\noindent(iii) For all $x\in\lbrack0<f<r_{0}]$
\[
|\nabla f|(x)\geq\frac{1}{k}.
\]
\end{theorem}

\noindent\textbf{Proof.} The equivalence of (ii) and (iii) follows from
\cite[Theorem 2.1]{AzeCor2004} and is based on  Ekeland variational
principle. Definition~\ref{DefinitionMetricRegF} (metric regularity of
multifunctions) yields the following restatement for (i):

\smallskip

\noindent$(i)_{1}$ For every
$(\bar{x},\bar{r})\in\mathrm{Graph\,}F$ with
$\bar{x}\in\lbrack0<f<r_{0}]$ and $\bar{r}\in(0,r_{0}),$ there
exist $\varepsilon>0$ and $\delta>0$ such that
\begin{equation}
(x,r)\in\left(  B(\bar{x},\varepsilon)\cap\lbrack0<f<r_{0}]\right)
\times\lbrack(\bar{r}-\delta,\bar{r}+\delta)\cap(0,r_{0})]\Longrightarrow
\mathrm{dist\,}(x,[f\leq r])\mathrm{\,}\leq\mathrm{\,}k\mathrm{\,}
(f(x)-r)^{+}.\label{aris7}
\end{equation}
Clearly $(i)\Rightarrow(i)_{1}$. Now, in order to prove $(i)_{1}
\Rightarrow(i)$, consider
$(\bar{x},\bar{r})\in\mathrm{Graph\,}F\cap
\lbrack0<f<r_{0}]\times(0,r_{0})$. Take $\varepsilon$ and $\delta$
positive given by $(i)_{1}$ such that
$0<\bar{r}-\delta<\bar{r}+2\delta<r_{0}$, $\varepsilon\leq
k(r_{0}-\bar{r}-2\delta)$ and $f$ is positive in $B(\bar
{x},\varepsilon)$ ($f$ is lower semicontinuous so $[f>0]$ is
open). For any $(x,r)\in
B(\bar{x},\varepsilon)\times(\bar{r}-\delta,\bar{r}+\delta)$, we
have $r\in(0,r_{0})$ and $f(x)>0.$ Thus if $f(x)<r_{0}$ by
$(i)_{1}$ we have
\[
\mathrm{dist\,}(x,[f\leq r])\leq k(f(x)-r)^{+}=k\,\mathrm{dist\,}(r,F(x)).
\]
If $f(x)\geq r_{0},$ then
\begin{align*}
\mathrm{dist\,}(x,[f\leq
r])\leq\mathrm{dist\,}(x,\bar{x})+\mathrm{dist\,} (\bar{x},[f\leq
r]) &  \leq\varepsilon+k\,(f(\bar{x})-r)^{+}\\ &
\leq\varepsilon+k\delta\\ &  \leq k(r_{0}-\bar{r}-\delta)\\ &
\leq k(r_{0}-r)\\ &  \leq k(f(x)-r)^{+}=k\,\mathrm{dist\,}(r,F(x)).
\end{align*}
Thus $(i)_{1}\Rightarrow(i)$.

\smallskip

It is now straightforward to see that $(ii)\Longrightarrow(i),$ thus it
remains to prove that $(i)_{1}\Longrightarrow(ii).$ To this end, fix any
$k^{\prime}>k$, $r_{1}\in(0,r_{0})$ and $x_{1}\in[f=r_{1}]$. We shall prove
that
\[
\mbox{\rm dist}\,(x_{1},[f\leq s])\leq k^{\prime}(r_{1}-s),
\]
for all $s\in(0,r_{1}].$

\smallskip

\noindent\textbf{Claim 1} Let $r\in(0,r_{0})$ and $x\in\lbrack f=r]$. Then
there exist $r^{-}<r$ and $x^{-}\in\lbrack f=r^{-}]$ such that
\begin{equation}
d(x,x^{-})\leq k^{\prime}(r-r^{-})\label{xmoins}
\end{equation}
with
\[
\mbox{\rm dist}\,(x,[f\leq s])\leq k^{\prime}(r-s)\text{,\quad for
all } s\in\lbrack r^{-},r].
\]
[\textit{Proof of Claim 1.} Apply $(i)_{1}$ at $(x,r)\in\mbox{Graph}\;F$ to
obtain the existence of $\rho\in(0,r)$ such that $\mbox{\rm dist}\,(x,[f\leq
s])\leq k(r-s)$ for all $s\in\lbrack\rho,r]$. Since $k^{\prime}>k$ there
exists $x^{-}\in\lbrack f\leq\rho]$ satisfying
\[
d(x,x^{-})<\frac{k^{\prime}}{k}\mathrm{dist\,}(x,[f\leq\rho]),
\]
which in view of (\ref{aris7}) yields
\[
d(x,x^{-})<k^{\prime}\mathrm{\,}(r-\rho).
\]
To conclude, set $r^{-}=f(x^{-})\leq\rho$ and observe that for any
$s\in\lbrack r^{-},\rho]$ we have
\[
\mbox{\rm dist}\,(x,[f\leq s])\leq d(x,x^{-})\leq k^{\prime}(r-\rho)\leq
k^{\prime}(r-s)=k^{\prime}(f(x)-s).
\]
This completes the proof of the claim.$\hfill\diamondsuit$]

\medskip

Let $\mathcal{A}$ be the set of all families $\{(x_{i},r_{i})\}_{i\in
I}\subset\lbrack f\leq r_{1}]\times\mathbb{R}$ containing $(x_{1},r_{1})$ such that\medskip

-- (P$_{1}$)\quad $f(x_{i})=r_{i}$ for all $i\in I$ {and} $r_{i}\neq r_{j},$
for $i\neq j$ ;\medskip

-- (P$_{2}$)\quad If $i,j\in I$ and $r_{i}<r_{j}$ then $d(x_{j},x_{i} )\le
k^{\prime}\mathrm{\,}(r_{j}-r_{i}).$ ;\medskip

-- (P$_{3}$)\quad For $r^{\ast}=\inf\{r_{i}:i\in I\}$ and for $s\in(r^{\ast
},r_{1}]$ we have:
\[
\mbox{\rm dist}\,(x_{1},[f\leq s])\leq k^{\prime}(r_{1}-s).
\]
\smallskip The set $\mathcal{A}$ is nonempty (it contains the one--element
family $\{(x_{1},r_{1})\}$) and can be ordered by the inclusion relation (that
is, $\mathcal{J}_{1}\preceq\mathcal{J}_{2}$ if, and only if, $\mathcal{J}
_{1}\subset\mathcal{J}_{2}$). Under this relation $\mathcal{A}$ becomes a
totally ordered set: every totally ordered chain in $\mathcal{A}$ has an upper
bound in $\mathcal{A}$ (its union). Thus, by Zorn lemma, there exists a
maximal element $\mathcal{M}=\{(x_{i},r_{i})\}_{i\in I}$ in $\mathcal{A}$.

\smallskip

\noindent\textbf{Claim 2}. Any maximal element $\mathcal{M}=\{(x_{i}
,r_{i})\}_{i \in I}$ of $\mathcal{A}$ satisfies
\begin{equation}
r^{*}=\inf_{i\in I} r_{i} \leq0. \label{aris11}
\end{equation}

\noindent\lbrack\textit{Proof of the Claim 2.} Let us assume,
towards a contradiction, that (\ref{aris11}) is not true,
\textit{i.e. }$r^{\ast}>0$. Let us first assume that there exists
$j\in I$ such that $r^{\ast}=r_{j}$. Define
$r^{-}:=r_{j}^{-}<r_{j}$ and $x_{j}^{-}=x^{-}\in\lbrack f=r^{-}]$
as specified in Claim 1 and consider the family
$\mathcal{M}_{1}=\mathcal{M} \cup\{(x^{-},r^{-})\}$. Then
$\mathcal{M}_{1}$ clearly complies with (P$_{1} $). To see that
$\mathcal{M}_{1}$ satisfies (P$_{2}$), simply observe that for
each $i\in I$,
\[
d(x^{-},x_{i})\leq d(x^{-},x_{j})+d(x_{j},x_{i})\leq k^{\prime}(r_{i}-r^{-}).
\]
Let $s\in\lbrack r^{-},r_{j}]$. By using the properties of the
couple $(x^{-},r^{-})$, one obtains
\[
\mbox{\rm dist}\,(x_{1},[f\leq s])\leq\mbox{\rm dist}\,(x_{1},x_{j})+\mbox{\rm
dist}\,(x_{j},[f\leq s])\leq k^{\prime}(r_{1}-r_{j})+k^{\prime}(r_{j}-s)\leq
k^{\prime}(r_{1}-s).
\]
This means that $\mathcal{M}_{1}\in\mathcal{A}$ which is contradicts the
maximality of $\mathcal{M}$.

\smallskip

\noindent Thus it remains to treat the case when the infimum $r^{\ast}$ is not
attained. Let us take any decreasing sequence $\{r_{i_{n}}\}_{n\geq1}
,\;i_{n}\in I$ satisfying $r_{i_{1}}=r_{1}$ and $r_{i_{n}}\searrow r^{\ast}.$
For simplicity the sequences $\{r_{i_{n}}\}_{n}$ and $\{x_{i_{n}}\}_{n}$ will
be denoted, respectively, by $\{r_{n}\}_{n}$ and $\{x_{n}\}_{n}$. Applying
(P$_{2}$) we obtain
\begin{equation}
d(x_{n},x_{n+m})\leq k^{\prime} \mathrm{\,}(r_{n}-r_{n+m}).
\label{aris10}
\end{equation}
It follows that $\{x_{n}\}_{n\geq1}$ is a Cauchy sequence, thus it converges
to some $x^{\ast}$. Taking the limit as $m\rightarrow+\infty$ we deduce from
(\ref{aris10}) that $d(x_{n},x^{\ast})\leq k^{\prime}\mathrm{\,}(r_{n}
-r^{\ast}),$ for all $n\in\mathbb{N}^{\ast}$. For any $i\in I$, there exists
$n$ such that $r_{n}<r_{i}$ and therefore
\begin{equation}
\mathrm{dist\,}(x^{\ast},x_{i})\leq
d(x^{\ast},x_{n})+d(x_{n},x_{i})\leq
k^{\prime}(r_{i}-r^{\ast})\leq k^{\prime}(r_{i}-f(x^{\ast})),
\label{treize}
\end{equation}
where the last inequality follows from the lower semicontinuity of $f$. Set
$f(x^{\ast})=\rho^{\ast}\leq r^{\ast}$ and $\mathcal{M}_{1}=\mathcal{M}
\cup\{(x^{\ast},\rho^{\ast})\}$. Since the infimum is not attained in
$\inf\{r_{i}:i\in I\}$ the family $\mathcal{M}_{1}$ satisfies (P$_{1}$).
Further by using (\ref{treize}), we see that $\mathcal{M}_{1}$ complies also
with (P$_{2}$). Take $s\in\lbrack\rho^{\ast},r^{\ast}]$. Since $x^{\ast}
\in\lbrack f\leq s]$, we have
\[
\mbox{\rm dist}\,(x_{1},[f\leq s])\leq\mbox{\rm dist}\,(x_{1},x^{\ast})\leq
k^{\prime}(r_{1}-r^{\ast})\leq k^{\prime}(r_{1}-s).
\]
Hence $\mathcal{M}_{1}$ belongs to $\mathcal{A}$ which contradicts the
maximality of $\mathcal{M}$.$\hfill\diamondsuit$]

\medskip

The desired implication follows easily by taking the limit as $k^{\prime}$
goes to $k$. This completes the proof.$\hfill\Box$

\begin{remark}
[Sublevel mapping and Lipschitz continuity]\label{RemarkLipLevel}
It is straightforward to see that statement (ii) above is
equivalent to the ``Lipschitz continuity'' (see (\ref{Hausdorff}))
of the sublevel set application

\[
\left\{
\begin{array}
[c]{cll} (0,r_{0}) & \rightrightarrows &  X\\ r & \longmapsto &
\lbrack f\leq r]
\end{array}
\right.
\]
for the Hausdorff ``metric'' given in (\ref{Hausdorff}).
Note that $F^{-1}$ is exactly the sublevel mapping given above, and thus in
this context the Lipschitz continuity of $F^{-1}$ is equivalent to the Aubin
property of $F^{-1}$, see \cite{dqz06, Ioffe}.
\end{remark}

\subsection{Metric regularity and K{\L} inequality}

As an immediate consequence of Theorem~\ref{PropositionIC} and
Remark~\ref{RemarkLipLevel}, we have the following result.

\begin{corollary}
[K\L-inequality and sublevel set mapping]\label{Corollary_KurMetric}Let
$f:X\longrightarrow\mathbb{R}\cup\{+\infty\}$ be a lower semicontinuous
function defined on a complete metric space $X$ and let $\varphi\in
\mathcal{K}(0,{r}_{0})$ (see (\ref{Kr})). The following assertions are equivalent:\smallskip

(i) the multivalued mapping
\[
\left\{
\begin{array}
[c]{cll} X & \rightrightarrows & \mathbb{R}\\ x & \mapsto &
[(\varphi\circ f)(x),+\infty)
\end{array}
\right.
\]

is $k$-metrically regular on $[0<f<r_{0}]\times(0,\varphi(r_{0}))$ ;\smallskip

(ii) for all $r_{1},r_{2}\in(0,r_{0})$
\[
\mathrm{Dist\,}([f\leq r_{1}],[f\leq
r_{2}])\leq\mathrm{\,}k\mathrm{\,}
|\varphi(r_{1})-\varphi(r_{2})|\text{ ;}
\]

(iii) for all $x\in\lbrack0<f<r_{0}]$
\[
|\nabla(\varphi\circ f)|(x)\geq\frac{1}{k}.
\]
\end{corollary}

It might be useful to observe the following:

\begin{remark}
[Change of metric]\label{top} Let $\varphi\in\mathcal{K}(0,{r}_{0})$ and
assume that it can be extended continuously to an increasing function still
denoted $\varphi:\mathbb{R}_{+}\rightarrow\mathbb{R}_{+}$. Set $d_{\varphi
}(r,s)=|\varphi(r)-\varphi(s)|$ for any $r,s\in\mathbb{R}_{+}$ and assume that
$\mathbb{R}_{+}$ is endowed with the metric $d_{\varphi}$. Endowing $\R_+$ with this new
metric, assertions (i), (ii) and (iii) can be reformulated very simply:\smallskip

(i ') The multivalued mapping
\[
\left\{
\begin{array}
[c]{cll} X & \rightrightarrows & \mathbb{R_{+}}\\ x & \mapsto &
[f(x),+\infty)
\end{array}
\right.
\]
is $k$-metrically regular on $[0<f<r_{0}]\times(0,r_{0})$.\smallskip

(ii') The sublevel mapping
\[
\mathbb{R}_{+}\ni r\mapsto\lbrack f\leq r],
\]
is $k$ Lipschitz continuous on $(0,r_{0})$.\smallskip

(iii') For all $x\in\lbrack0<f<r_{0}]$
\[
|\nabla_{\varphi}f|(x)\geq\frac{1}{k},
\]
where $|\nabla_{\varphi}f|$ denotes the strong slope of the restricted
function $\bar{f}:[0<f]\rightarrow\lbrack\mathbb{R}_{+},d_{\varphi}]$.
\end{remark}

Given a lower semicontinuous function
$f:X\longrightarrow\mathbb{R} \cup\{+\infty\}$ we say that $f$ is
\emph{strongly slope-regular}, if for each point $x$ in its domain
$\mathrm{dom\,}f$ one has
\begin{equation}
|\nabla f|(x)=|\nabla(-f)|(x).\label{ssreg}
\end{equation}
Note that all $C^{1}$ functions are strongly slope-regular according to the
above definition.

\begin{proposition}
[Level mapping and Lipschitz
continuity]\label{PropositionLipLev}Assume
$f:X\rightarrow\mathbb{R}$ is continuous and strongly
slope-regular. Then any of the assertions (i)--(iii) of
Theorem~\ref{PropositionIC} is equivalent to the fact that the
\emph{level set application}
\[
\left\{
\begin{array}
[c]{cll} \mathbb{R} & \rightrightarrows &  X\\ r & \mapsto & [f=r]
\end{array}
\right.
\]
is Lipschitz continuous on $(0,r_{0})$ with respect to the Hausdorff metric.
\end{proposition}

\noindent\textbf{Proof}. The result follows by applying
Theorem~\ref{PropositionIC} twice. (Details are left to the reader.)$\hfill
\Box$

\bigskip

Let us finally state the following important corollary.

\begin{corollary}
[K\L-inequality and level set mapping]\label{CorollaryKurLojLipLev} Let
$f:X\longrightarrow\mathbb{R}$ be a continuous function which is strongly
slope-regular on $[0<f<r_{0}]$ and let $\varphi\in\mathcal{K}(0,r_{0})$
(recall (\ref{Kr})). Then the following assertions are equivalent:\smallskip

\noindent(i) $\varphi\circ f$ is $k$-metrically regular on
$[0<f<r_{0} ]\times(0,\varphi(r_{0}))$;\newline \noindent(ii) for
all $r_{1},r_{2} \in(0,r_{0})$
\[
\mathrm{Dist\,}([f=r_{1}],[f=r_{2}])\leq\mathrm{\,}k\mathrm{\,} |\varphi
(r_{1})-\varphi(r_{2})|;
\]
(iii) for all $x\in\lbrack0<f<r_{0}]$
\[
|\nabla(\varphi\circ f)|(x)\geq\frac{1}{k}.
\]
\end{corollary}

\noindent\textbf{Proof.} It follows easily by combining
Theorem~\ref{PropositionIC} with Proposition~\ref{PropositionLipLev}.
$\hfill\Box$

\section{K{\L}--inequality in Hilbert spaces}

From now on, we shall work on a real Hilbert space
$[H,\langle\cdot ,\cdot\rangle]$. Given a vector $x$ in $H$, the
norm of $x$ is defined by $||x||=\sqrt{\langle x,x\rangle}$ while
for any subset $C$ of $H$, we set
\begin{equation}\label{c-form}
||C||_{-}=\mbox{\rm dist}\;(0,C)=\inf\{||x||:x\in C\}\in\mathbb{R}
\cup\{+\infty\}.
\end{equation}
Note that $C=\emptyset$ implies $||C||_{-}=+\infty$.

\subsection{Elements of nonsmooth analysis}

Let us first recall the notion of Fr\'{e}chet subdifferential (see
\cite{CLSW98, Morduk}).

\begin{definition}
[Fr\'{e}chet subdifferential]\label{Definition_Frechet} Let $f:H\rightarrow
\mathbb{R}\cup\{+\infty\}$ be a real-extended-valued function. We say that
$p\in H$ is a (Fr{\'{e}}chet) subgradient of $f$ at $x\in\mbox{dom\,}f$ if
\[
\liminf_{y\rightarrow x,\;y\neq x}\,\frac{f(y)-f(x)-\langle p,y-x\rangle
}{||y-x||}\geq0.
\]
\end{definition}

\noindent We denote by $\partial f(x)$ the set of Fr\'{e}chet subgradients of
$f$ at $x$ and set $\partial f(x)=\emptyset$ for $x\notin\mbox{dom\,}f$. Let
us now define the notion of critical point in variational analysis.

\begin{definition}
[critical point/values]\label{Definition_critical} (i) A point $x_{0}\in H$ is called
\emph{critical} for the function $f,$ if $0\in\partial f(x_{0}).$\\
(ii) The value
$r\in f(H)$ is called a critical value, if $[f=r]$ contains at least one
critical point.
\end{definition}

In this section we shall mainly deal with the class of
\emph{semiconvex} functions. Let us give the corresponding
definition. (The reader should be aware that the terminology is not
yet completely fixed in this area, so that the notion of
semiconvex function may vary slightly from one author to another.)

\begin{definition}
[semiconvexity]\label{Definition_semiconvex}A proper lower semicontinuous
function $f$ is called \emph{semiconvex} (or convex up to a square) if for
some $\alpha>0$ the function
\[
x\longmapsto f(x)+\frac{\alpha}{2}||x||^{2}
\]
is convex.
\end{definition}

\begin{remark}
\label{Remark_semiconvex} \ \newline  (i) For each $x\in H$, $\partial f(x)$ is a
(possibly empty) closed convex subset of $H$ and $\partial f(x)$ is nonempty
for $x\in\mathrm{int\,dom\,}f.$\smallskip

\noindent(ii) It is straightforward from the above definition that the
multivalued operator $x\longmapsto\partial f(x)+\alpha x$ is (maximal)
monotone (see \cite[Definition~12.5]{Rock98} for the definition). \smallskip

\noindent(iii) For general properties of semiconvex functions, see
\cite{ac99}. Let us mention that Definition \ref{Definition_semiconvex} is
equivalent to the fact that
\begin{equation}
f(y)-f(x)\geq\langle p, y-x\rangle-\alpha||x-y||^{2},
\label{croiss-grad-semiconv}
\end{equation}
for all $x,y\in H$ and all $p\in\partial f(x)$ (where $\alpha>0$).

\noindent(iii) According to
Definition~\ref{Definition_semiconvex}, semiconvex functions are
contained in several important classes of (nonsmooth) functions,
as for instance $\phi$-convex functions (\cite{DMT1985}), weakly
convex functions (\cite{ADT2005}) and primal--lower--nice
functions (\cite{Thibault}). Although an important part of the
forthcoming results is extendable to these more general classes,
we shall hereby sacrifice extreme generality in sake of simplicity
of presentation.
\end{remark}

Given a real-extended-valued function $f$ on $H,$ we define the
\emph{remoteness} (\textit{i.e.}, distance to zero) of its subdifferential
$\partial f$ at $x\in H$ as follows:
\begin{equation}
||\partial f(x)||_{-}=\,\underset{p\in\partial f(x)}{\inf}
\,||p||\,=\,\mathrm{dist} \,(0,\partial f(x)). \tag{remoteness}
\end{equation}

\begin{remark}\label{rmq-min-norm} (minimal norm)
 \newline (i) If $\partial f(x)\not = \emptyset,$ the infimum in the above
definition is achieved since $\partial f(x)$ is a nonempty closed convex set. If we
define $\partial^{0} f(x)$ as the projection of $0$ on the closed convex set
$\partial f(x)$ we of course have
\begin{equation}
||\partial f(x)||_{-}=||\partial^{0} f(x)||.
\end{equation}
Some properties of $H\ni x\mapsto ||\partial f(x)||_{-}$ are given in
Section~\ref{Annex} (Annex).\smallskip

\noindent(ii) If $f$ is a semiconvex function, then $||\partial f(x)||_{-}$
coincides with the notion of strong slope $|\nabla f|(x)$ introduced in
(\ref{slope}), see Lemma~\ref{Lemma_slope} (Annex).
\end{remark}

\subsection{Subgradient curves: basic properties}

\label{sec:subgradient-curves}

Let $f:H\rightarrow\mathbb{R}\cup\{+\infty\}$ be a proper lower semicontinuous
semiconvex function. The purpose of this subsection is to recall the main properties of the
trajectories (subgradient curves) of the corresponding differential
inclusion:
\[
\left\{
\begin{array}
[c]{l} \dot{\chi}_{x}(t)\in-\partial f(\chi_{x}(t))\quad\text{a.e.
on }(0,+\infty),\\
\\
\chi_{x}(0)=x\in\mbox{dom\,}f.
\end{array}
\right.
\]

The following statement aggregates useful results concerning
existence and uniqueness of solutions. These results are
essentially known even for a more general class of functions (see
\cite[Theorem~2.1, Proposition~2.14, Theorem~3.3]{Thibault} for
instance for the class of primal--lower--nice functions). It should also be noticed
 that the integration of measurable curves of the form $\R \ni t\rightarrow\gamma(t)\in H$ relies on Bochner integration/measurability  theory 
(basic properties can be found in \cite{Brezis}). 

\begin{theorem}
[subgradient curves]\label{Theorem_thibault} For every $x\in\mbox{dom\,}f$
there exists a unique absolutely continuous curve (called trajectory or
subgradient curve) $\chi_{x}:[0,+\infty)\rightarrow H$ that satisfies
\begin{equation}
\left\{
\begin{array}
[c]{l} \dot{\chi}_{x}(t)\in-\partial f(\chi_{x}(t))\quad\text{a.e.
on }(0,+\infty),\\
\\
\chi_{x}(0)=x\in\mbox{dom\,}f.
\end{array}
\right.  \label{subgr}
\end{equation}
Moreover the trajectory satisfies:\smallskip

\begin{itemize}
\item [(i)]$\chi_{x}(t)\in\mbox{dom\,}\partial f$ for all $t\in(0,+\infty)$.\smallskip

\item[(ii)] For all $t>0$, the right derivative $\dot{\chi}_{x}(t^{+})$ of
$\chi_{x}$ is well defined and equal to
\[
\dot{\chi}_{x}(t^{+})=-\partial^{0}f(\chi_{x}(t)).
\]
In
particular $\dot{\chi}_{x}(t)=-\partial^{0}f(\chi_{x}(t))$, for almost all $t$.\smallskip

\item[(iii)] The mapping $t\mapsto||\partial f(\chi_{x}(t))||_{-}$ is
right-continuous at each $t\in(0,+\infty)$.\smallskip

\item[(iv)] The function $t\longmapsto f(\chi_{x}(t))$ is nonincreasing and
continuous on $[0,+\infty)$. Moreover, for all $t,\tau\in\lbrack0,+\infty)$
with $t\le\tau$, we have
\[
f(\chi_{x}(t))-f(\chi_{x}(\tau))\,\geq\,\int_{t}^{\tau}||\dot{\chi}
_{x}(u)||^{2}\,du\,,
\]
and equality holds if $t>0$.

\item[(v)] The function $t\longmapsto f(\chi_{x}(t))$ is Lipschitz continuous
on $[\eta,+\infty)$ for any $\eta>0.$ Moreover
\[
\frac{d}{dt}f(\chi_{x}(t))=-||\dot{\chi}_{x}(t)||^{2}\mbox{ a.e on
} (\eta,+\infty).
\]
\end{itemize}
\end{theorem}

\noindent\textbf{Proof}. The only assertion that does not appear explicitly in
\cite{Thibault} is the continuity of the function $f\circ\chi_{x}$ at $t=0$
when $x\in\mbox{dom\,}f\diagdown\mbox{dom\,}\partial f$, but this is an easy
consequence of the fact that $f$ is lower semicontinuous, $\chi_{x}$ is
(absolutely) continuous and $f\circ\chi_{x}$ is decreasing. For the rest of
the assertions we refer to \cite{Thibault}.$\hfill\Box$

\bigskip

The following result asserts that the semiflow mapping associated with the
differential inclusion (\ref{subgr}) is continuous. This type of result can be
established by standard techniques and therefore is essentially known (see
\cite{Brezis,Thibault} for example). We give here an outline of proof (in case
that $f$ is semiconvex) for the reader's convenience.

\begin{theorem}
[continuity of the semiflow]\label{continuous} For any semiconvex function $f$
the semiflow mapping
\[
\left\{
\begin{array}
[c]{lll} \mathbb{R}_{+}\times\mbox{{\rm dom\,}}f & \rightarrow &
H\\ (t,x) & \mapsto & \chi_{x}(t)
\end{array}
\right.
\]
is (norm) continuous on each subset of the form $[0,T]\times(B(0,R)\cap\lbrack
f\leq r])$ where $T,R>0$ and $r\in\mathbb{R}$.
\end{theorem}

\noindent\textbf{Proof}. Let us fix $x,y\in\mbox{dom\,}f$ and $T>0$. Then for
almost all $t\in\lbrack0,T],$ there exist $p(\chi_{x}(t))\in\partial
f(\chi_{x}(t))$ and $q(\chi_{y}(t))\in\partial f(\chi_{y}(t))$ such that
\[
\frac{d}{dt}||\chi_{x}(t)-\chi_{y}(t)||^{2}=2\langle\chi_{x}(t)-\chi
_{y}(t),\dot{\chi}_{x}(t)-\dot{\chi}_{y}(t)\rangle=-2\langle\chi_{x}
(t)-\chi_{y}(t),p(\chi_{x}(t))-q(\chi_{y}(t))\rangle.
\]
It follows by \eqref{croiss-grad-semiconv} that
\[
\frac{d}{dt}||\chi_{x}(t)-\chi_{y}(t)||^{2}\leq2\alpha||\chi_{x}(t)-\chi
_{y}(t)||^{2},
\]
which implies (using Gr\"{o}nwall's lemma) that for all $0\leq t\leq T$ we
have
\begin{equation}
||\chi_{x}(t)-\chi_{y}(t)||^{2}\leq\exp(2\alpha T)||x-y||^{2}.
\label{1}
\end{equation}
For any $0\leq t\leq s \leq T,$ using Cauchy--Schwartz inequality and
Theorem~\ref{Theorem_thibault} we deduce that
\begin{equation}
||\chi_{x}(s)-\chi_{x}(t)||\leq\int_{t}^{s}||\dot{\chi}_{x}(\tau)||d\tau
\leq\sqrt{s-t}\sqrt{\int_{s}^{t}||\dot{\chi}_{x}(\tau)||^{2}d\tau}\leq
\sqrt{s-t}\sqrt{f(x)}. \label{2}
\end{equation}
The result follows by combining (\ref{1}) and (\ref{2}).$\hfill\Box$

\bigskip

%We shall also use the following result (see \cite{Thibault}) concerning the
%asymptotic behavior of the norm precompact trajectories.
%
%\begin{proposition}
%[asymptotic behavior of trajectories]\label{asymptotics} Let us assume that
%$f$ is semiconvex and that for some $x\in\mbox{dom\,}\partial f$ the
%trajectory $\{\chi_{x}(t)\}_{t\in\lbrack0,+\infty)}$ is relatively compact for
%the strong topology. Then the strong cluster points of $t\longmapsto\chi
%_{x}(t)$ as $t$ goes to infinity are critical points of $f$ and
%\[
%l:=\lim_{t\rightarrow+\infty}f(\chi_{x}(t))
%\]
%is a critical value of $f$ (\textit{i.e.} there exists $x_{\infty}\in H$ such
%that $\partial f(x_{\infty})\ni0$ and $f(x_{\infty})=l$).
%\end{proposition}
%
%\bigskip

Let us introduce the notions of a \emph{piecewise absolutely continuous curve}
and of a \emph{piecewise subgradient curve}. This latter notion, due to its
robustness, will play a central role in our study.

\begin{definition}
\label{piecewise}Let $a,b\in\lbrack-\infty,+\infty]$ with $a<b$.\newline
\emph{(Piecewise absolutely continuous curve)} A curve $\gamma
:(a,b)\rightarrow H$ is said to be \emph{piecewise absolutely continuous} if
there exists a countable partition of $(a,b)$ into intervals $I_{k}$ such that
the restriction of $\gamma$ to each $I_{k}$ is absolutely
continuous.\smallskip\newline \emph{(Length of a curve)} Let $\gamma
:(a,b)\rightarrow H$ be a piecewise absolutely continuous curve. The length of
$\gamma$ is defined by
\[
\mbox{\rm length}\;[\gamma]:=\int_{a}^{b}||\dot{\gamma}(t)||\;dt.
\]
\smallskip\emph{(Piecewise subgradient curve)} Let $T\in(0,+\infty]$. A curve
$\gamma:[0,T)\rightarrow H$ is called a piecewise subgradient curve for
(\ref{subgr}) if there exists a countable partition of $[0,T]$ into
(nontrivial) intervals $I_{k}$ such that:

\indent-- the restriction $\gamma|_{I_{k}}$ of $\gamma$ to each interval
$I_{k}$ is a subgradient curve ;

\indent-- for each disjoint pair of intervals $I_{k},I_{l}$, the intervals
$f(\gamma(I_{k}))$ and $f(\gamma(I_{l}))$ have at most one point in common.
\end{definition}

Note that piecewise subgradient curves are piecewise absolutely continuous.
Observe also that subgradient curves satisfy the above definition in a trivial way.

\subsection{Characterizations of the K\L -inequality}

In this section we state and prove one of the main results of this work. Let
$f:H\rightarrow\mathbb{R}\cup\{+\infty\}$ and $\bar{x}\in\lbrack f=0]$ be a
critical point. Throughout this section the following assumptions will be used:\smallskip

-- There exist $\bar{r},\bar{\epsilon}>0$ such that
\begin{equation}
x\in\bar{B}(\bar{x},\bar{\epsilon})\cap\lbrack0<f\leq\bar{r}]\;\Longrightarrow
\;0\notin\partial f(x)\quad\quad\quad\quad\mathrm{(0\ is\ a\
locally\ upper} \text{ }\mathrm{isolated}\text{ }\mathrm{critical\
value)}.\label{crit}
\end{equation}

-- There exist $\bar{r},\bar{\epsilon}>0$ such that
\begin{equation}
\bar{B}(\bar{x},\bar{\epsilon})\cap\lbrack f\leq\bar{r} ]\
\mathrm{is\ (norm)\
compact}\quad\quad\quad\quad\mathrm{(local}\text{
}\mathrm{sublevel}\text{
}\mathrm{compactness)}.\label{compacite-ss-niv}
\end{equation}

\begin{remark}\ \\
(i) The first condition can be seen as a Sard-type condition.\smallskip

\noindent (ii) Assumption \eqref{compacite-ss-niv} is always satisfied in
finite-dimensional spaces, but is also satisfied in several
interesting cases involving infinite-dimensional spaces. Here are
two elementary examples.\smallskip\newline 
\noindent(ii)$_1$ The (convex)
function $f:\ell^{2}(\mathbb{N})\rightarrow\mathbb{R}$ defined by
\[
f(x)=\sum_{n\geq1}n^{2}x_{i}^{2}
\]
has compact lower level sets.\smallskip\newline 
(ii)$_2$ Let $g:\mathbb{R}
\rightarrow\mathbb{R}\cup\{+\infty\}$ be a proper lower
semicontinuous semiconvex function and let
$\Phi:L^{2}(\Omega)\rightarrow\mathbb{R} \cup\{+\infty\}$ be as
follows (\cite{brezis71})
\[
\Phi(x)=\left\{
\begin{array}
[c]{l} \frac{1}{2}\int_{\Omega}||\nabla
x||^{2}+\int_{\Omega}g(x)\mbox { if }x\in H^{1}(\Omega)\\
+\infty\mbox{ otherwise.}
\end{array}
\right.
\]
The above function is a lower semicontinuous semiconvex function
and the sets of the form $[\Phi\leq r]\cap B(\bar{x},R)$ are
relatively compact in $L^{2}(\Omega)$ (use the compact embedding
theorem of $H^{1}(\Omega )\hookrightarrow
L^{2}(\Omega)$).
\end{remark}

As shown in Theorem~\ref{Theorem_main}, Kurdyka-\L ojasiewicz inequality can
be characterized in terms of boundedness of the length of ``worst (piecewise
absolutely continuous) curves'', that is those defined by the points of less
steepest descent.

\begin{definition}
[Talweg/Valley]\label{Definition_valley}Let $\bar{x}\in[f=0]$ be a critical
point of $f$ and assume that \eqref{crit} holds for some $\bar{r},
\bar{\epsilon}>0.$ Let $D$ be any closed bounded set that contains $B(\bar
{x},\bar{\epsilon})\cap\lbrack0<f\leq\bar{r}].$ For any $R>1$ the $R$-valley
$\mathcal{V}_{R}(\cdot)$ of $f$ around $\bar{x}$ is defined as follows:
\begin{equation}
\mathcal{V}_{R}(r)=\left\{  x\in[f=r]\cap D:\,||\partial
f(x)||_{-} \,\leq\,R\,\underset{y\in[f=r]\cap D}{\inf}\,||\partial
f(y)||_{-}\right\} ,\; \text{ for all }r\in(0,\bar{r}].
\label{valley}
\end{equation}
A selection $\theta: (0,\bar{r}]\rightarrow H$ of $\mathcal{V}_{R}$,
\textit{i.e.} a curve such that $\theta(r)\in\mathcal{V}_{R}(r), \forall r
\in(0,\bar{r}]$, is called an $R$-talweg or simply a talweg.
\end{definition}

We are ready to state the main result of this work.

\begin{theorem}
[Subgradient inequality -- local
characterization]\label{Theorem_main} Let
$f:H\rightarrow\mathbb{R}\cup\{+\infty\}$ be a lower
semicontinuous semiconvex function and $\bar{x}\in[f=0]$ be a
critical point. Assume that there exist $\bar{\epsilon},\bar{r}>0$
such that \eqref{crit} and \eqref{compacite-ss-niv} hold.

Then, the following statements are equivalent:\smallskip

\noindent(i) \textbf{[Kurdyka-\L ojasiewicz inequality]} There
exist $r_{0} \in(0,\bar{r}),\;\epsilon\in(0,\bar{\epsilon})$ and
$\varphi\in\mathcal{K} (0,r_{0})$ such that
\begin{equation}
||\partial(\varphi\circ f)(x)||_{-}\geq1,\qquad\text{for all }
x\in\bar {B}(\bar{x},\epsilon)\cap\lbrack0<f\leq r_{0}].
\label{Kurdyk}
\end{equation}

\noindent(ii) \textbf{[Length boundedness of subgradient curves]} There exist
$r_{0}\in(0,\bar{r}),\;\epsilon\in(0,\bar{\epsilon})$ and a strictly
increasing continuous function $\sigma:[0,r_{0}]\rightarrow\lbrack0,+\infty)$
with $\sigma(0)=0$ such that for all subgradient curves $\chi_{x}$ of
(\ref{subgr}) satisfying $\chi_{x}([0,T))\subset\bar{B}(\bar{x},\epsilon
)\cap\lbrack0<f\leq r_{0}]$ ($T\in(0,+\infty]$) we have
\[
\int_{0}^{T}||\dot{\chi}_{x}(t)||dt\leq\sigma(f(x))-\sigma(f(\chi_{x}(T))).
\]

\noindent(iii) \textbf{[Piecewise subgradient curves have finite length]}
There exist $r_{0}\in(0,\bar{r}),\;\epsilon\in(0,\bar{\epsilon})$ and $M>0$
such that for all piecewise subgradient curves $\gamma:[0,T)\rightarrow H$ of
(\ref{subgr}) satisfying $\gamma([0,T))\subset\bar{B}(\bar{x},\epsilon
)\cap\lbrack0<f\leq r_{0}]$ ($T\in(0,+\infty]$) we have
\[
\mbox{\rm length}[\gamma]:=\int_{0}^{T}||\dot{\gamma}(\tau)||d\tau<M.
\]

\noindent(iv) \textbf{[Talwegs of finite length] }For every $R>1$, there exist
$r_{0}\in(0,\bar{r}),\;\epsilon\in(0,\bar{\epsilon}),$ a closed bounded subset
$D$ containing $B(\bar{x},\epsilon)\cap\lbrack0<f\leq r_{0}]$ and a piecewise
absolutely continuous curve $\theta:(0,r_{0}]\rightarrow H$ of \emph{finite
length} which is a selection of the valley $\mathcal{V}_{R}(r),$ that is,
\[
\theta(r)\in\mathcal{V}_{R}(r),\text{ for all }r\in(0,r_{0}].
\]

\noindent(v) \textbf{[Integrability condition]} There exist $r_{0}\in
(0,\bar{r})$ and $\epsilon\in(0,\bar{\epsilon})$ such that the function
\[
u(r)=\frac{1}{\underset{x\in\bar{B}(\bar{x} ,\epsilon)\cap[f=r]}{\inf
}\,{||\partial f(x)||_{-}}},\;\;r\in(0,r_{0}]
\]
is finite-valued and belongs to $L^{1}(0,r_{0})$.
\end{theorem}

\begin{remark} \hspace*{0mm}  \newline
\label{rmk-thm111} 
\noindent(i) As it appears clearly in the proof,
statement $(iv)$ can be replaced by $(iv^{\prime})$ ``\emph{There
exist $R>1$,} $r_{0}
\in(0,\bar{r}),\;\epsilon\in(0,\bar{\epsilon}),$ a closed bounded
subset $D$ containing $B(\bar{x},\epsilon)\cap\lbrack0<f\leq
r_{0}]$ and a piecewise absolutely continuous curve
$\theta:(0,r_{0}]\rightarrow H$ of \emph{finite length} which is a
selection of the valley $\mathcal{V}_{R}(r),$ that is,
\[
\theta(r)\in\mathcal{V}_{R}(r),\text{ for all }r\in(0,r_{0}]^{\prime\prime}.
\]
(ii) The compactness assumption \eqref{compacite-ss-niv} is only
used in the proofs of $(iii)\Rightarrow(ii)$ and
$(ii)\Rightarrow(iv)$. Hence if this assumption is removed, we
still have:
\[
(iv)\Longrightarrow(iv^{\prime})\Longrightarrow(v)\Longleftrightarrow
(i)\Longrightarrow(ii)\Longrightarrow(iii).
\]
(iii) Note that (i) implies condition (\ref{crit}). This follows immediately
from the chain rule (see Annex, Lemma~\ref{chain}).
\end{remark}

\noindent\textbf{Proof of Theorem \ref{Theorem_main}.} \textbf{[(i)}$\Rightarrow$\textbf{(ii)]} Let
$\epsilon,r_{0},\varphi$ be as in (i) such that (\ref{Kurdyk}) holds. Let
further $\chi_{x}$ be a subgradient curve of (\ref{subgr}) for $x\in
\lbrack0<f\leq r_{0}]$ and assume that $\chi_{x}([0,T))\subset\bar{B} (\bar
{x},\epsilon)\cap\lbrack0<f\leq r_{0}]$ for some $T>0$.\smallskip

\noindent Let us first assume that $x\in\mbox{dom\,}\partial f$. Since
$\varphi$ is $C^{1}$ on $(0,r_{0})$, by Theorem~\ref{Theorem_thibault}(v) and
Lemma~\ref{chain} (Annex) we deduce that the curve $t\mapsto\varphi(f(\chi
_{x}(t))$ is absolutely continuous with derivative
\[
\frac{d}{dt}(\varphi\circ f\circ\chi_{x})(t)=-\varphi^{\prime}(f(\chi
_{x}(t))||\dot{\chi}_{x}(t)||^{2}\mbox{ a.e.
on }(0,T).
\]
Integrating both terms on the interval $(0,T)$ and recalling (\ref{Kurdyk}),
$\chi_{x}(0)=x$ we get
\begin{align*}
\varphi(f(x))-\varphi(f(\chi_{x}(T))) & =-\int_{0}^{T}\frac{d}{dt}
(\varphi\circ f\circ\chi_{x})(t)dt\\ &
=\int_{0}^{T}\varphi^{\prime}(f(\chi_{x}(t))||\dot{\chi}_{x}(t)||^{2}
dt\geq\int_{0}^{T}||\dot{\chi}_{x}(t)||dt.
\end{align*}
Thus (ii) holds true for $\sigma:=\varphi$ and for all subgradient curves
starting from points in $\mbox{dom\,}\partial f.$ Let now $x\in\mbox
{dom\,}f\diagdown\mbox{dom\,}\partial f$ and fix any $\delta\in(0,T).$ Since
$\chi_{x}([\delta,T])\subset\mbox{dom\,}\partial f$ we deduce from the above
that
\[
\int_{\delta}^{T}||\dot{\chi}_{x}(t)||dt\leq\sigma(f(\chi_{x}(\delta
))-\sigma(f(\chi_{x}(T))).
\]
Thus the result follows by taking $\delta\searrow0^{+}$ and using the
continuity of the mapping $t\longmapsto f(\chi_{x}(t))$ at $0$
(Theorem~\ref{Theorem_thibault}(ii)).

\medskip

\noindent\textbf{[(ii)}$\Rightarrow$\textbf{(iii)]} Let $\gamma$ be a
piecewise subgradient curve as in (iii) and let $I_{k}$ be the associated
partition of $[0,T]$ (\textit{cf}. Definition~\ref{piecewise}). Let
$\{a_{k}\}$ and $\{b_{k}\}$ be two sequences of real numbers such that
$\mbox{int}\;I_{k} =(a_{k},b_{k})$. Since the restriction $\gamma|_{I_{k}}$ of
$\gamma$ onto $I_{k}$ is a subgradient curve, applying (ii) on $(a_{k}
,b_{k})$ we get
\[
\mbox{length}\;[\gamma|_{I_{k}}]\leq\sigma(f(\gamma(a_{k})))-\sigma
(f(\gamma(b_{k}))).
\]
Let $m$ be an integer and $I_{k_{1}},\dots,I_{k_{m}}$  a finite subfamily
of the partition. We may assume that these intervals are ordered as follows $0\le
a_{k_{1}}\le b_{k_{1}}\le\cdots\le a_{k_{m}}\le b_{k_{m}}$. Hence
\[
\sum_{1}^{m}\left[  \sigma(f(\gamma(a_{k_{i}})))-\sigma(f(\gamma(b_{k_{i}})))
\right]  \leq\sigma(f(\gamma(a_{k_{1}})))\le\sigma(r_{0}).
\]
Thus the family $\{\sigma(f(\gamma(a_{k})))-\sigma(f(\gamma(b_{k})))\}$ is
summable, hence using the definition of Bochner integral (see \cite{Brezis})
\[
\mathrm{length}\;[\gamma]\;=\;\sum_{k\in\mathbb{N}}\mathrm{length}
\;[\gamma|_{I_{k}}]\;\leq\;\sigma(r_{0}).
\]

\smallskip

\noindent\textbf{[(iii)}$\Rightarrow$\textbf{(ii)]} Let $\epsilon,r_{0}$ be as
in (iii), pick any $0\leq r^{\prime}<r\leq r_{0}$ and denote by $\Gamma
_{r^{\prime},r}$ the (nonempty) set of piecewise subgradient curves
$\gamma:[0,T)\rightarrow H\;$ (where $T\in(0,+\infty]$) such that
\[
\gamma([0,T))\subset\bar{B}(\bar{x},\epsilon)\cap\lbrack r^{\prime}<f\leq r].
\]
Note that, by Theorem~\ref{Theorem_thibault}(iv) and Proposition
\ref{slopes}(iii), $T=+\infty$ is possible only when $r^{\prime}=0$. Set
further
\[
\psi(r^{\prime},r):=\sup_{\gamma\in\Gamma_{r^{\prime},r}} \,\mbox{\rm
length}[\gamma]\qquad\text{and}\qquad\sigma(r):=\psi(0,r).
\]
Note that (iii) guarantees that $\psi$ and $\sigma$ have finite values. We can
easily deduce from Definition~\ref{piecewise} that
\begin{equation}
\psi(0,r^{\prime})+\psi(r^{\prime},r)=\psi(0,r). \label{sigma}
\end{equation}
Thus for each $x\in\bar{B}(\bar{x},\epsilon)\cap\lbrack0<f\leq r_{0}]$ and
$T>0$ such that $\chi_{x} ([0,T])\subset B(\bar{x},\epsilon)\cap\lbrack0<f\leq
r_{0}],$ we have
\begin{equation}
\int_{0}^{T}||\dot{\chi}_{x}(\tau)||d\tau+\sigma(f(\chi_{x}(T))\leq
\sigma(f(x)). \label{leng}
\end{equation}
Since the function $\sigma$ is nonnegative and increasing it can be extended
continuously at $0$ by setting $\sigma(0)=\lim_{t\downarrow0}\sigma(t)\geq0$.
Since the property (\ref{leng}) remains valid if we replace $\sigma(\cdot)$ by
$\sigma(\cdot)-\sigma(0)$, there is no loss of generality to assume
$\sigma(0)=0$.

\smallskip

To conclude it suffices to establish the continuity of $\sigma$ on
$(0,r_{0} ]$. Fix $\tilde{r}$ in $(0,r_{0})$ and take a
subgradient curve $\chi :[0,T)\rightarrow H$ satisfying
$\chi([0,T))\subset\bar{B}(\bar{x} ,\epsilon)\cap[f\leq r_{0}]$,
where $T\in(0,+\infty]$. Set $f(\chi(0))=r$ and
$\lim_{t\rightarrow T} f(\chi(t))=r^{\prime}$ and assume that
$\tilde{r}\leq r^{\prime}\leq r\leq r_{0}.$

\smallskip From Theorem~\ref{Theorem_thibault}(iv) and Proposition
\ref{slopes}(iii) (Annex), we deduce that $T<+\infty$ so that
$\chi ([0,T])\subset\bar{B}(\bar{x},\epsilon)\cap\lbrack
r^{\prime}\leq f\leq r]$. Using assumption (\ref{crit}) together
with Theorem~\ref{Theorem_thibault} (i),(v), we deduce that the
absolutely continuous function $f\circ
\chi:[0,T]\rightarrow\lbrack r^{\prime},r]$ is invertible and
\begin{equation}
\frac{d}{d\rho}[f\circ\chi]^{-1}(\rho)=\frac{-1}{||\dot{\chi}([f\circ
\chi]^{-1}(\rho)||^{2}}\geq\frac{-1}{\displaystyle\mathop{\rm
inf}_{x\in \bar{B}(\bar{x},\epsilon)\cap\lbrack\tilde{r}\leq f\leq
r_{0}]}||\partial f(x)||_{-}^{2}}:=-K,\label{inverse}
\end{equation}
for almost all $\rho\in(r,r^{\prime})$. By Proposition~\ref{slopes}(iii)
(Annex) we get that $K<+\infty$ and therefore the function $\rho
\longmapsto\lbrack f\circ\chi]^{-1}(\rho)$ is Lipschitz continuous with
constant $K$ on $[r^{\prime},r]$. Using the Cauchy-Schwarz inequality and
Theorem~\ref{Theorem_thibault}(iv) we obtain
\begin{align*}
\mathrm{length}\;[\chi]\, &
=\,\int_{0}^{T}||\dot{\chi}||\;\leq\;\sqrt
{T}\sqrt{\int_{0}^{T}||\dot{\chi}||^{2}}\;=\,\sqrt{[f\circ\chi]^{-1}
(r)-[f\circ\chi]^{-1}(r^{\prime})}\,\sqrt{\int_{0}^{T}||\dot{\chi}||^{2}}\\
& \\ &
\leq\sqrt{K(r-r^{\prime})}\sqrt{r-r^{\prime}}=\sqrt{K}(r-r^{\prime}).
\end{align*}
This last inequality implies that each piecewise subgradient curve
$\gamma:[0,T)\rightarrow H$ such that $\gamma([0,T))\subset\bar{B}(\bar
{x},\epsilon)\cap\lbrack r^{\prime}\leq f\leq r]$ satisfies
\[
\mathrm{length}\;[\gamma]\;\leq\;\sqrt{K}(r-r^{\prime}),
\]
thus using (\ref{sigma}) we obtain $\sigma(r)-\sigma(r^{\prime})\leq\sqrt
{K}(r-r^{\prime})$, which yields the continuity of $\sigma$.

\medskip

\noindent\textbf{[(ii)}$\Rightarrow$\textbf{(iv)]} Let us assume
that (ii) holds true for $\epsilon$ and $r_{0}$. In a first step
we establish the existence of a closed bounded subset $D$ of
$[0<f\leq r_{0}]$ satisfying
\begin{equation}
x\in D,\;t\geq0,\;f(\chi_{x}(t))>0\;\Rightarrow\;\chi_{x}(t)\in D.
\label{stable}
\end{equation}
Let $r_{0}\geq r_{1}>0$ be such that $\sigma(r_{1})<\epsilon/3$ and let us
set
\[
D:=\{y\in\bar{B}(\bar{x},\epsilon)\cap\lbrack0<f\leq r_{1}]:\exists x\in
\bar{B}(\bar{x},\epsilon/3)\cap\lbrack0<f\leq r_{1}],\exists t\ge0\mbox{
such that }\chi_{x}(t)=y\}.
\]

Let us first show that $D$ enjoys property (\ref{stable}). It suffices to
establish that
\[
x\in\bar{B}(\bar{x},\epsilon/3)\cap[0<f\leq r_{1}],\; t \geq0,\;f(\chi
_{x}(t))>0\Rightarrow\chi_{x}(t)\in D.
\]
To this end, fix $x\in\bar{B}(\bar{x},\epsilon/3)\cap[0<f\leq r_{1}].$ By
continuity of the flow, we observe that $\chi_{x}(t)\in\bar{B}(\bar
{x},\epsilon)$ for small $t>0$ and for all $t\geq0$ such that $\chi
_{x}([0,t])\subset\bar{B}(\bar{x},\epsilon)$ with $f(\chi_{x}(t))>0$,
assumption (ii) yields
\begin{equation}
||\chi_{x}(t)-\bar{x}||\;\leq\;||\chi_{x}(t)-x||\,+\,||x-\bar{x}||\;\leq
\;\int_{0}^{t}||\dot{\chi}_{x}(\tau)||d\tau+\epsilon/3\;\leq\;\sigma
(r_{1})+\epsilon/3\leq2\epsilon/3. \label{well}
\end{equation}
Thus $D$ satisfies \eqref{stable} and $\bar{B}(\bar{x},\epsilon/3)\cap[f\leq
r_{1}]\subset D.$

\smallskip

Let us now prove that $D$ is (relatively) closed in $[0<f\leq r_{1}]$. Let
$y_{k}\in D$ be a sequence converging to $y$ such that $f(y)\in(0,r_{1}]$.
Then there exist sequences $\{x_{n}\}_{n}\subset\bar{B}(\bar{x},\epsilon
/3)\cap\lbrack0<f\leq r_{1}]$ and $\{t_{n}\}_{n}\subset\mathbb{R}_{+}$ such
that $\chi_{x_{n}}(t_{n})=y_{n}$. Since $f$ is lower semicontinuous, there
exists $n_{0}\in\mathbb{N}$ and $\eta>0$ such that $f(y_{n})>\eta$ for all
$n\geq n_{0}$. By Theorem~\ref{Theorem_thibault}(ii),(iv), \eqref{crit} and
Proposition~\ref{slopes}(iii) (Annex), we obtain for all $n\geq n_{0}$
\[
0<\;t_{n\,}\,\inf_{z\in\lbrack\eta\leq f\leq r_{1}]\cap\bar{B}(\bar
{x},\epsilon)}||\partial f(z)||_{-}^{2}\;\leq\; \int_{0}^{t_{n}}||\dot{\chi
}_{x_{n}}(t)||^{2}dt\,\leq\,f(x_{n})\,\leq\,r_{1}.
\]
The above inequality shows that the sequence $\{t_{n}\}_{n}$ is bounded. Using
a standard compactness argument we therefore deduce that, up to an extraction,
$x_{n}\rightarrow\tilde{x}$ and $t_{n}\rightarrow\tilde{t}$ for some
$\tilde{x}\in\bar{B} (\bar{x},\epsilon/3)\cap\lbrack f\leq r_{1}]$ and
$\tilde{t}\in\mathbb{R}_{+} $. Theorem~\ref{continuous} (continuity of the
semiflow) implies that $y=\chi_{\tilde{x}}(\tilde{t})$ and consequently that
$f(\tilde{x})\geq f(y)>0$, yielding that $y\in D$. This shows that $D$ is
(relatively) closed in $[0<f\leq r_{0}]$.

\smallskip

Now we build a piecewise absolute continuous curve in the valley. According to
the notation of Proposition~\ref{slopes} (Annex) we set
\[
s_{D}(r):=\inf\{||\partial f(x)||_{-}:x\in D\cap\lbrack f=r]\},
\]
so that for any $R>1$ the $R$-valley around $\bar{x}$ (cf.
Definition~\ref{Definition_valley}) is given by
\[
\mathcal{V}_{R}(r):=\{x\in[f=r]\cap D:\;||\partial f(x)||_{-}\;\leq
\;R\;s_{D}(r)\}.
\]
If
$\bar{B}(\bar{x},\epsilon/3)\cap[f=r]=\emptyset$ for all $0<r\leq
r_{1},$ there is nothing to prove. Otherwise, there exists
$0<r_{2}\leq r_{1}$ and
$x_{2}\in\bar{B}(\bar{x},\epsilon/3)\cap[f=r_{2}]\subset D.$
From Theorem~\ref{Theorem_thibault} and Proposition
\ref{slopes}(iii) (Annex), we deduce that $\chi_{x_{2}}(t)\in[f=f(\chi_{x_2}(t))]\cap
D\cap\mbox{dom\,}\partial f$ for all $t\geq 0$ such that $[f\circ\chi_{x_{2}}](t)>0$ and that the {\em inverse} function $[f\circ \chi_{x_{2}} ]^{-1}(\cdot)$ is defined 
on an interval containing $(0,r_2)$. In other words the set $[f=r]\cap
D\cap\mbox{dom\,}\partial f$ is nonempty for each $r\in (0,r_2)$, which in turn implies that the valley is nonempty for small positive values of $r$, {\it i.e.} $\mathcal{V}_{R}(r)\neq \emptyset$  for all $r\in(0,r_2)$. With no loss of generality
 we assume that $\mathcal{V}_{R}(r_2)\neq \emptyset$.

\smallskip

Let further $R^{\prime}\in(1,R)$ and $x\in[f=r_{2}]\cap D$ be such that
$||\partial f(x)||_{-}\,\leq\,R^{\prime}\,s_{D}(r_{2})$ (therefore, in
particular, $x\in\mathcal{V}_{R}(r_{2})$). Take $\rho\in(R^{\prime},R)$. Since
the mapping $t\longmapsto||\partial f(\chi_{x}(t)||_{-}$ is right--continuous
(\textit{cf}. Theorem~\ref{Theorem_thibault}(iii)), there exists $t_{0}>0$
such that $||\partial f(\chi_{x}(t)||_{-}<\rho s_{D}(r_{2})$ for all
$t\in(0,t_{0})$. On the other hand $t\longmapsto s_{D}(f(\chi_{x}(t))$ is
lower semicontinuous (\textit{cf}. Proposition~\ref{slopes}--Annex), hence
there exists $t_{1}\in(0,t_{0})$ such that $R\,s_{D}(f(\chi_{x}(t))\,>\,\rho
\,s_{D}(r_{2}),$ for all $t\in(0,t_{2} )$. Using the continuity of the mapping
$\chi_{x}(\cdot)$ and the stability property (\ref{stable}), we obtain the
existence of $t_{2}>0$ such that
\begin{equation}
\chi_{x}(t)\in\mathcal{V}_{R}(f(x(t))\mbox{ for all
}t\in\lbrack0,t_{2}). \label{t's}
\end{equation}
By using arguments similar to those of [(iii)$\Rightarrow$(ii)] we define the
following absolutely continuous curve:
\[
(f\circ\chi_{x}(t_{2}),r_{2}]\ni r\longmapsto\theta(r)=\chi_{x} ([f\circ
\chi_{x}]^{-1}(r))\in D\cap\lbrack f=r].
\]
By Proposition~\ref{selection} based on Zorn's Lemma (see Annex), we obtain a piecewise subgradient
curve that we still denote by $\theta$, defined on $(0,r_{2}]$, satisfying
$\theta(r)\in\mathcal{V}_{R}(r)$ for all $r\in(0,r_{2}]$. Assumption (iii) now
yields
\[
\mbox{length}\;[\theta]<M<+\infty,
\]
completing the proof of the assertion.

\medskip

\noindent\textbf{[(iv)}$\Rightarrow$\textbf{(v)]} Fix $R>1$ and let
$\epsilon,r_{0}$ and $\theta:(0,r_{0}]\rightarrow H$ be as in (iv). Applying
Lemma~\ref{chain} (Annex), we get
\[
\frac{d}{dr}(f\circ\theta)(r)=1=\langle\dot{\theta}(r),p(r)\rangle\mbox{
a.e on }(0,r_{0}]\text{, \qquad for all }p(r)\in\partial f(\theta(r)).
\]
Using the Cauchy-Schwartz inequality together with the fact that $D\cap\lbrack
f=r]\supset\bar{B}(\bar{x},\epsilon)\cap\lbrack f=r]$, we obtain
\[
R\,||\dot{\theta}(r)||\,\geq\,u(r)\;=\,\frac{1}{\inf_{x\in\bar{B}(\bar{x}
,\epsilon)\cap\lbrack f=r]} ||\partial f(x)||_{-}},
\]
for almost all $r\in(0,r_{0}]$. Since $\theta$ has finite length we deduce
that $u\in L^{1}((0,r_{0})$.

\medskip

\noindent\textbf{[(v)}$\Rightarrow$\textbf{(i)]} Let $\epsilon$, $r_{0}$ and
$u$ be as in $(v)$. From Proposition~\ref{slopes} (Annex) we deduce that $u$
is finite-valued and upper semicontinuous. Applying Lemma~\ref{LemmaReg}
(Annex) we obtain a continuous function $\bar{u} :(0,r_{0}]\rightarrow
(0,+\infty)$ such that $\bar{u}(r)\geq u(r)$ for all $r\in(0,r_{0}]$. We set
\[
\varphi(r)=\int_{0}^{r}\bar{u}(s)ds.
\]
It is directly seen that $\varphi(0)=0$, $\varphi\in C([0,r])\cap
C^{1}(0,r_{0})$ and $\varphi^{\prime}(r)>0$ for all $r\in(0,r_{0})$. Let
$x\in\bar{B}(\bar{x},\epsilon)\cap[f=r]$ and $q\in\partial(\varphi\circ
f)(x)$. From Lemma~\ref{chain} (Annex) we deduce $p:=\frac{q}{\varphi^{\prime
}(r)}\in\partial f(x)$, and therefore
\[
||q||\,=\,\varphi^{\prime}(r)\,||\frac{q}{\varphi^{\prime}(r)}||\,\geq
\,u(r)\,||p||\,\geq1.
\]
The proof is complete.$\hfill\Box$

\bigskip

Under a stronger compactness assumption Theorem~\ref{Theorem_main} can be
reformulated as follows.

\begin{theorem}
[Subgradient inequality -- global characterization]\label{global} Let
$f:H\rightarrow\mathbb{R}\cup\{+\infty\}$ be a lower semicontinuous semiconvex
function. Assume that there exists $r_{0}>0$ such that
\[
[f\leq r_{0}]\mbox{ is compact and }0\notin\partial f(x),\; \forall
x\in[0<f<r_{0}].
\]

Then the following propositions are equivalent\smallskip

\noindent(i) \textbf{[Kurdyka-\L ojasiewicz inequality]} There exists a
$\varphi\in\mathcal{K}(0,r_{0})$ such that
\[
||\partial(\varphi\circ f)(x)||_{-}\geq1,\qquad\text{for all } x\in\lbrack0<f<
r_{0}].
\]

\noindent(ii) \textbf{[Length boundedness of subgradient curves]} There exists
an increasing continuous function $\sigma:[0,r_{0})\rightarrow\lbrack
0,+\infty)$ with $\sigma(0)=0$ such that for all subgradient curves $\chi
_{x}(\cdot)$ (where $x\in\lbrack0<f<r_{0}]$) we have
\[
\int_{0}^{T}||\dot{\chi}_{x}(t)||\,dt\,\leq\,\sigma(f(x))-\sigma(f(\chi
_{x}(T))),
\]
whenever $f(\chi_{x}(T))>0.$\smallskip

\noindent(iii) {\textbf{[Piecewise subgradient curves have bounded length]}}
There exists $M>0$ such that for all piecewise subgradient curves
$\gamma:[0,T)\rightarrow H$ such that $\gamma([0,T))\subset\lbrack0<f<r_{0}]$
we have
\[
\mbox{\rm length}[\gamma]<M.
\]

\noindent(iv) \textbf{[Talwegs of finite length] }For all $R>1$, there exists
a piecewise absolutely continuous curve (with countable pieces) $\theta
:(0,r_{0})\rightarrow\mathbb{R}^{n}$ with \emph{finite length} such that
\[
\theta(r)\in\left\{  x\in\lbrack f=r]:\,||\partial f(x)||_{-}\,\leq
\,R\,\underset{y\in\lbrack f=r]}{\inf}\,||\partial f(y)||_{-}\right\}
,\qquad\text{for all }r\in(0,r_{0}).
\]

\noindent(v) \textbf{[Integrability condition]} The function $u:(0,r_{0}
)\rightarrow\lbrack0,+\infty]$ defined by
\[
u(r)=\frac{1}{\underset{x\in[f=r]}{\inf}\, {||\partial f(x)||_{-}}},\qquad
r\in(0,r_{0})
\]
is finite-valued and belongs to $L^{1}(0,r_{0} )$.\smallskip

\noindent(vi) \textbf{[Lipschitz continuity of the sublevel mapping]} There
exists $\varphi\in\mathcal{K}(0,r_{0})$ such that
\[
\mbox{\rm Dist}([f\leq r],[f\leq s])\leq|\varphi(r)-\varphi(s)|\qquad\text{for
all }r,s\in(0,r_{0}).
\]
\end{theorem}

\noindent\textbf{Proof} The proof is similar to the proof of
Theorem~\ref{Theorem_main} and will be omitted. The equivalence
between (i) and (vi) is a consequence of
Corollary~\ref{Corollary_KurMetric}.$\hfill\Box$

\subsection{Application: convergence of the proximal algorithm}

In this subsection we assume that the function $f:H\rightarrow\mathbb{R}
\cup\{+\infty\}$ is \emph{semiconvex} (\textit{cf}.
Definition~\ref{Definition_semiconvex}). Let us recall the
definition of the proximal mapping (see
\cite[Definition~1.22]{Rock98}, for example).

\begin{definition}
[proximal mapping]\label{Definition_proximal}Let $\lambda\in(0,\alpha^{-1}).$
Then the proximal mapping $\mbox{\rm
prox}_{\lambda}:H\rightarrow H$ is defined by
\[
\mbox{\rm prox}_{\lambda}(x):=\mbox{\rm argmin\,}\left\{  f(y)+\frac
{1}{2\lambda}||y-x||^{2}\right\}  ,\;\forall x\in H.
\]
\end{definition}

\begin{remark}
The fact that $\mbox{\rm prox}_{\lambda}$ is well-defined and single-valued is
a consequence of the semiconvex assumption: indeed this assumption implies
that the auxiliary function appearing in the aforementioned definition is
strictly convex and coercive (see~\cite{Rock98}, \cite{combettes} for instance).
\end{remark}

\begin{lemma}
[Subgradient inequality and proximal mapping]\label{Lemma_prox}
Assume that $f:H\rightarrow\mathbb{R}\cup\{+\infty\}$ is a
semiconvex function that satisfies (i) of Theorem~\ref{global}.
Let $x\in\lbrack0<f<r_{0}]$ be such that $f(\mbox{\rm
prox}_{\lambda}x)>0.$ Then
\begin{equation}
||\mbox{\rm
prox}_{\lambda}x-x||\leq\varphi(f(x))-\varphi(f(\mbox{\rm
prox}_{\lambda}x)).\label{tame}
\end{equation}
\end{lemma}

\noindent\textbf{Proof.} Set $x^{+}=\mbox{\rm prox}_{\lambda}(x)$, $r=f(x)$,
and $r^{+}=f(x^{+})$. It follows from the definition of $x^{+}$ that
$0<r^{+}\leq r<r_{0}$. In particular, for every $u\in\lbrack f\leq r^{+}]$ we
have
\[
||x^{+}-x||^{2}\leq||u-x||^{2}+2\lambda[f(u)-r^{+}]\leq||u-x||^{2}.
\]
Therefore by Corollary~\ref{Corollary_KurMetric} (Lipschitz continuity of the
sublevel mapping) we obtain
\[
||x^{+}-x||=\mbox{dist}\;(x,[f\leq r^{+}])\leq\mbox{\rm Dist}\;([f\leq
r],[f\leq r^{+}])\leq\varphi(r)-\varphi(r^{+}).
\]
The proof is complete.$\qquad\hfill\Box$

\bigskip

The above result has an important impact in the asymptotic analysis of the
\emph{proximal algorithm} (see forthcoming Theorem~\ref{Theorem_Tprox}). Let
us first recall that, given a sequence of positive parameters $\{\lambda
_{k}\}\subset(0,\alpha^{-1})$ and $x\in H$ the proximal algorithm is defined
as follows:
\[
Y_{x}^{k+1}=\mbox{\rm prox}_{\lambda_{k}}Y_{x}^{k},\qquad Y_{x}^{0}=x,
\]
or in other words
\[
\{Y_{x}^{k+1}\}=\mbox{\rm argmin\,}\left\{  f(u)+\frac{1}{2\lambda_{k}}
||u-Y_{x}^{k}||^{2}\right\}  ,\qquad Y_{x}^{0}=x.
\]
If we assume in addition that $\inf f>-\infty$, then for any initial point $x$
the sequence $\{f(Y_{x}^{k})\}$ is decreasing and converges to a real number
$L_{x}$.

\begin{theorem}
[strong convergence of the proximal
algorithm]\label{Theorem_Tprox} Let
$f:H\rightarrow\mathbb{R}\cup\{+\infty\}$ be a semiconvex function
which is bounded from below. Let $x\in\mbox{dom\,}f,$
$\{\lambda_{k}\}\subset (0,\alpha^{-1})$ and
$L_{x}:=\,\underset{k\rightarrow\infty}{\lim}f(Y_{x} ^{k})$ and
assume that there exists $k_{0}\geq0$ and $\varphi\in
\mathcal{K}(0,f(Y_{x}^{k_{0}})-L_{x})$ such that
\begin{equation}
||\partial(\varphi\circ\lbrack
f(\cdot)-L_{x}])(x)||_{-}\geq1,\qquad\text{for all }x\in\lbrack
L_{x}<f\leq f(Y_{x}^{k_{0}})].\label{loj}
\end{equation}
Then the sequence $\{Y_{x}^{k}\}$ converges strongly to $Y_{x}^{\infty}$ and
\begin{equation}
||Y_{x}^{\infty}-Y_{x}^{k}||\,\leq\,\varphi(f(Y_{x}^{k})-L_x),\qquad\text{for
all }k\geq k_{0}.\label{ESTIM}
\end{equation}
\end{theorem}

\noindent\textbf{Proof} Since the sequence $\{Y_{x}^{k}\}_{k\geq k_{0}}$
evolves in $L_{x}\leq f<f(Y_{x}^{k_{0}})$, Lemma~\ref{Lemma_prox} applies.
This yields
\[
\sum_{k=p}^{q}||Y_{x}^{k+1}-Y_{x}^{k}||\leq\varphi(f(Y_{x}^{q+1}
)-L_x)-\varphi(f(Y_{x}^{p})-L_x),
\]
for all integers $k_{0}\leq p\leq q$. This implies that $Y_{x}^{k}$ converges
strongly to $Y_{x}^{\infty}$ and that inequality (\ref{ESTIM}) holds.$\hfill
\Box$

\bigskip

\begin{remark}
[Step-size]\textquotedblleft Surprisingly" enough the step-size sequence
$\{\lambda_{k}\}$ does not appear explicitly in the estimate (\ref{ESTIM}),
but it is instead hidden in the sequence of values $\{f(Y_{x}^{k})\}$. In
practice the choice of the step-size parameters $\lambda_{k}$ is however crucial to
obtain the convergence of $\{f(Y^{k} )\}$ to a critical value; standard
choices are for example sequences satisfying $\sum\lambda_{k}=+\infty$ or
$\lambda_{k}\in\lbrack\eta,\alpha^{-1})$ for all $k\geq0$ where $\eta
\in(0,\alpha^{-1})$, see \cite{combettes} for more details.
\end{remark}

\section{Convexity and K\L -inequality}

In this section, we assume that
$f:H\rightarrow\mathbb{R}\cup\{+\infty\}$ is a lower
semicontinuous proper convex function such that $\mathrm{inf}
\,f>-\infty$. Changing $f$ in $f-\inf f$, we may assume that
$\inf f=0$. Let us also denote the set of minimizers of
$f$ by
\[
C:=\mathrm{argmin}\,f=[f=0].
\]
When $C$ is nonempty, we may assume with no loss of
generality that $0\in C$.\smallskip

In this convex setting Theorem~\ref{Theorem_thibault} can be considerably
reinforced; related results are gathered in Section
\ref{sec:trajectoires-convexes}. We also recall well-known facts ensuring that
subgradient curves have finite length and provide a new result in that
direction (see Theorem~\ref{thm-hyperplan}). In
Section~\ref{sec:global-convexe}, we give some conditions which ensure that
$f$ satisfies the K\L -inequality and we show that the conclusions of
Theorem~\ref{global} can somehow be globalized. In
section~\ref{sec:contre-exemple} we build a counterexample of a $C^{2}$ convex
function in $\mathbb{R}^{2}$ which does not satisfy the K\L -inequality. This
counterexample also reveals that the uniform boundedness of the lengths of
subgradient curves is a strictly weaker condition than condition (iii) of
Theorem~\ref{Theorem_main}, which justifies further the introduction of
piecewise subgradient curves.

\subsection{Lengths of subgradient curves for convex functions}

\label{sec:trajectoires-convexes}

The following lemma gathers well known complements to Theorem
\ref{Theorem_thibault} when $f$ is convex.

\begin{lemma}
\label{prop-conv} Let $f: H\to\mathbb{R}\cup\{+\infty\}$ be a
lower semicontinuous proper convex function such that $0\in
C=[f=0].$ Let $x_{0} \in\mbox{\rm dom\,} f.$

\begin{itemize}
\item [(i)]If $a\in C,$ then
\begin{align*}
\frac{d}{dt}||\chi_{x_{0}}(t)-a ||^{2}\leq-2f(\chi_{x_{0}}(t))\leq
0\ \ \mathrm{a.e \ on } \ (0,+\infty).
\end{align*}
and therefore $t\mapsto||\chi_{x_{0}}(t)-a ||$ is nonincreasing.

\item[(ii)] The function $t\mapsto f(\chi_{x_{0}}(t))$ is nonincreasing and
converges to $0=\mathrm{min}\,f$ as $t\to+\infty.$

\item[(iii)] The function $t\in[0,+\infty)\longmapsto||\partial f(\chi_{x_{0}
}(t)||_{-}$ is nonincreasing.

\item[(iv)] The function $t\mapsto f(\chi_{x_{0}}(t))$ is convex and belongs
to $L^{1}([0,+\infty))$: for all $T>0,$
\begin{align}
\int_{0}^{T} f(\chi_{x_{0}}(t))dt = \frac{1}{2}||x_{0}||^{2} -
\frac{1} {2}||\chi_{x_{0}}(T)||^{2} \leq\frac{1}{2}||x_{0}||^{2}.
\end{align}

\item[(v)] For all $T>0,$
\begin{align}
\int_{0}^{T} ||\dot\chi_{x_{0}}(t)||dt \leq\left(  \int_{0}^{+\infty}
f(\chi_{x_{0}}(t))dt\right)  ^{1/2} (\mathrm{log}\,T)^{1/2}.
\end{align}
\end{itemize}
\end{lemma}

\noindent\textbf{Proof.} The proofs of these classical properties
can be found in \cite{Brezis,bruck75}. $\hfill\Box$

\bigskip

R. Bruck established in \cite{bruck75} that subgradient trajectories of convex
functions are always weakly converging to a minimizer in $C=\mbox{\rm
argmin\,}f$ whenever the latter is nonempty. However, as shown later on by
J.-B. Baillon \cite{baillon78}, strong convergence does not hold in general.

\smallskip

To the best of our knowledge, the problem of the characterization of length
boundedness of subgradient curves for convex functions is still open (see
\cite[Open problems, p.167]{Brezis}). In the present framework, the following
result of H. Br\'ezis \cite{brezis71, Brezis} is of particular interest.

\begin{theorem}
[Uniform boundedness of trajectory lengths \cite{brezis71}]\label{argmin-nonvide}
 Let $f:H\rightarrow\mathbb{R}\cup\{+\infty\}$ be a
lower semicontinuous proper convex function such that $0\in
C=\mathrm{argmin} \,f=[f=0].$ We assume that $C$ has nonempty
interior. Then, for all $x_{0} \in\mbox{\rm dom\,}f,$
$\chi_{x_{0}}(\cdot)$ has finite length. More precisely, if
$B(0,\rho)\subset C,$ we have, for all $T\geq0,$
\[
\int_{0}^{T}||\dot{\chi}_{x_{0}}(t))||dt\leq\frac{1}{2\rho}(||x_{0}
||^{2}-||\chi_{x_{0}}(T)||^{2}).
\]
\end{theorem}

\noindent\textbf{Proof.} We assume that $B(0,\rho)\subset C$ for
some $\rho>0$ and consider $x_{0}\in\mbox{\rm dom\,}f\backslash C$
(otherwise there is nothing to prove). Let $t\geq0$ such that
$\chi_{x_{0}}(t)\notin C$ and $\dot{\chi}_{x_{0}}(t)$ exists. By
convexity, we get
\[
\langle-(\chi_{x_{0}}(t)\!-\!\rho u),\dot{\chi}_{x_{0}}(t)\rangle\geq
f(\chi_{x_{0}}(t))-f(\rho u)>0
\]
for all $u$ in the unit sphere of $H$. As a consequence
$-\langle\chi_{x_{0}
}(t),\dot{\chi}_{x_{0}}(t)\rangle>\rho||\dot{\chi}_{x_{0}}(t)||$.
Therefore
$\int_{0}^{T}||\dot{\chi}_{x_{0}}(t)||dt\leq\frac{1}{2\rho}(||x_{0}
||^{2}-||\chi_{x_{0}}(T)||^{2}).$ $\hfill\Box$

\medskip

The following result is an extension of Theorem \ref{argmin-nonvide} under the
assumption that the vector subspace $\mathrm{span}(C)$ generated by $C$, has
codimension 1 in $H$. We denote by $\mathrm{ri}(C)$ the relative interior of
$C$ in $\mathrm{span}(C)$.

\begin{theorem}
\label{thm-hyperplan} Let
$f:H\rightarrow\mathbb{R}\cup\{+\infty\}$ be a lower
semicontinuous proper convex function such that $0\in
C=\mathrm{argmin} \,f=[f=0].$ Assume that $C\ $generates a
subspace of codimension 1 and that the relative interior
$\mathrm{ri}(C)$ of $C$ in $\mathrm{span}(C)$ is not empty. If
$x_{0}\in dom\,f$ is such that $\chi_{x_{0}}(t)$ converges
(strongly) to $a\in\mathrm{ri}(C)$ as $t\rightarrow+\infty,$ then
$\mathrm{length}\,[\chi_{x_{0}}]<+\infty.$
\end{theorem}

\noindent\textbf{Proof.} Let us denote by $a$ the limit point of
$\chi(t):=\chi_{x_{0}}(t)$ as $t$ goes to infinity. By assumption
$a$ belongs to $\mathrm{ri}(C)$, so that there exists $\delta>0$
such that $\bar {B}(a,\delta)\cap \mathrm{span}(C)\subset C$. Let $T>0$ be such
that $\chi(t)\in B(a,\delta)$ for all $t\geq T$. Write
$\mathrm{span}(C)=\{x\in H:\langle x,x^{\ast} \rangle=0\}$ with
$x^{\ast}\in H.$ We claim that the function $[T,+\infty)\ni
t\mapsto h(t)=\langle x^{\ast},\chi(t)\rangle$ has a constant
sign. Let us argue by contradiction and assume that there exist
$T<t_{1}<t_{2}$ such that $h(t_{1})<0<h(t_{2})$. Hence there
exists $t_{3}\in(t_{1},t_{2})$ such that $h(t_{3})=0$. Since
$\chi(t)\in B(a,\delta)$, this implies $\chi(t_{3})\in C$ and thus
by the uniqueness theorem for subgradient curves
(Theorem~\ref{Theorem_thibault}), we have $\chi(t)=\chi(t_{3})$
for all $t\geq t_{3}$ which is a contradiction. Note also that if
$h(t_{0})=0$ for some $t_{0}\geq T$, then $\chi$ has finite
length. Indeed applying once more Theorem~\ref{Theorem_thibault},
we deduce that $\chi(t)=\chi(t_{0})$ for all $t\geq t_{0}$, hence
\[
\int_{0}^{+\infty}||\dot{\chi}||=\int_{0}^{t_{0}}||\dot{\chi}||\leq\sqrt
{t_{0}}\;\sqrt{\int_{0}^{t_{0}}||\dot{\chi}||^{2}}<+\infty.
\]
Assume that $h$ is positive (the case $h$ negative can be treated
similarly) and define the following function
\[
\tilde{f}(x)=\left\{
\begin{array}
[c]{ll} 0 & \mbox{ if }\langle x,x^{\ast}\rangle<0\mbox{ and
}x\in\bar{B}(a,\delta)\\ f(x) & \mbox{ if }\langle
x,x^{\ast}\rangle\geq0\mbox{ and }x\in\bar{B}(a,\delta)\\ +\infty
& \mbox{ otherwise.}
\end{array}
\right.
\]
One can easily check that $\tilde{f}$ is proper, lower semicontinuous, convex
and that $\mbox{\rm argmin\,}\tilde{f}$ has non empty interior. Note also that
$\partial\tilde{f}(x)=\partial f(x)$ for all $x\in B(a,\delta)$ such that
$\langle x,x^{\ast}\rangle>0$. The conclusion follows from the previous result
and the fact that $\dot{\chi}(t)+\partial\tilde{f}(\chi(t))\ni0\mbox{ a.e. on
}(T,+\infty).$ $\hfill\Box$

\subsection{K\L -inequality for convex functions}

\label{sec:global-convexe}

The following result shows that if $f$ is convex, then the function $\varphi$
of Theorem~\ref{Theorem_main}(i) can be assumed to be concave and defined on
$[0,\infty)$.

\begin{theorem}
[Subgradient inequality -- convex case]\label{convex} Let $f: H\to
\mathbb{R}\cup\{+\infty\}$ be a lower semicontinuous proper convex function
which is bounded from below (recall that $\inf f=0$). The following statements are
equivalent:\smallskip\newline

(i) There exist $r_{0}>0$ and $\varphi\in\mathcal{K}(0,r_{0})$ such that
\[
||\partial(\varphi\circ f)(x)||_{-}\geq1,\qquad\text{for all }x\in
\lbrack0<f\leq r_{0}].
\]

\smallskip(ii) There exists a \textbf{concave} function $\psi\in
\mathcal{K}(0,\infty)$ such that
\begin{equation}
||\partial(\psi\circ f)(x)||_{-}\geq1,\qquad\text{for all
}x\notin\lbrack f=0].\label{Kurcon}
\end{equation}
\end{theorem}

\noindent\textbf{Proof.} The implication (ii)$\Longrightarrow$(i)
is obvious. To prove (i)$\Longrightarrow$(ii) let us first
establish that the function
\[
r\in(0,+\infty)\longmapsto u(r)=\frac{1}{\displaystyle\mathop{\rm
inf} _{x\in\lbrack f=r]}||\partial f(x)||_{-}}
\]
is finite-valued and nonincreasing. Let $0<r_{2}<r_{1}$ and let us
show that $u(r_{2})\geq u(r_{1}).$ To this end we may assume with
no loss of generality that $u(r_{1})>0$ (and therefore that
$[f=r_{1}]\cap\mbox{dom\,}\partial f$ is nonempty). Take
$\epsilon>0$ and let $x_{1}\in\lbrack f=r_{1}]$ and $p_{1}
\in\partial f(x_{1})$ such that
$u(r)\leq\frac{1}{||p_{1}||}+\epsilon$. Since the continuous
function $t\mapsto f(\chi_{x_{1}}(t))$ tends to $\mathrm{inf} f=0$
as $t$ goes to infinity (see \cite{Lemaire} for instance), there
exists $t_{2}>0$ such that $f(\chi_{x_{1}}(t_{2}))=r_{2}$. From
Lemma~\ref{prop-conv} (iii), we obtain
\[
\frac{1}{||\partial f(\chi_{x_{1}}(t_{2})||_{-}}\geq\frac{1}{||p_{1}||}\geq
u(r_{1})-\epsilon,
\]
which yields $u(r_{2})\geq u(r_{1})$. By (i) the function $u$ is finite-valued
on $(0,r_{0})$, thus, since $u$ is nonincreasing, it is also finite-valued on
$(0,+\infty)$.

\smallskip

It is easy to see that [(i)$\Rightarrow$(v)] of Theorem~\ref{Theorem_main} holds without
the compactness assumption \eqref{compacite-ss-niv} (see
Remark~\ref{rmk-thm111}). It follows that $u\in L^{1}(0,r_{0})$ and by
Lemma~\ref{LemmaReg} (Annex) that there exists a decreasing continuous
function $\tilde{u}\in L^{1}(0,r_{0})$ such that $\tilde{u}\geq u.$
Reproducing the proof of $(v)\Rightarrow(i)$ of Theorem~\ref{Theorem_main} we
obtain a strictly increasing, concave, $C^{1}$ function
\[
\psi(r):=\int_{0}^{r}\tilde{u}(s)ds
\]
for which (\ref{Kurcon}) holds for all $x\in\lbrack0<f<r_{0}]$. Fix $\bar
{r}\in(0,r_{0})$ and take $\psi$ as above. Applying (\ref{Kurcon}) and using
the fact that $u(r)$ is decreasing we obtain
\[
1\leq\psi^{\prime}(\bar{r})u(\bar{r})^{-1}\leq\psi^{\prime}(\bar{r}
)u(r)^{-1}\leq\psi^{\prime}(\bar{r})||p||,
\]
for all $p\in\partial f(x)$, $x\in\lbrack\bar{r}\leq f]$ and $r\in(\bar
{r},+\infty)$ such that $u(r)>0$. This shows that the function $\Psi
:\mathbb{R}_{+}\rightarrow\mathbb{R}_{+}$ defined by
\[
\Psi(r):=\left\{
\begin{array}
[c]{ll} \psi(r) & \mbox{ if }r\leq\bar{r},\\
\psi(\bar{r})+\psi^{\prime}(\bar{r})(r-\bar{r}) &
\mbox{otherwise}.
\end{array}
\right.
\]
satisfies the required properties.$\hfill\Box$

\bigskip

A natural question arises: when does a convex function $f$ satisfy
the K\L --inequality? In finite-dimensions a quick positive answer
can be given whenever $f$ belongs to an o-minimal structure
(convexity then becomes superflous). The following result gives an
alternative criterion when $f$ is not extremely
``flat'' around its set of minimizers. More precisely, we assume
the following growth condition:
\begin{equation}
\left\{
\begin{array}
[c]{l} \text{There exists
$m:[0,+\infty)\rightarrow\lbrack0,+\infty)$ and $S\subset H$ such
that}\\[6pt] \text{$m$ is continuous, increasing, $m(0)=0$, $f\geq
m(\mathrm{dist} (\cdot,C))$ on $S\cap\mathrm{dom}\,f$}\\[6pt]
\text{ and
$\displaystyle\int_{0}^{\rho}\,\frac{m^{-1}(r)}{r}dr<+\infty$
(for some $\rho>0$). }
\end{array}
\right.  \label{hyp-croissance}
\end{equation}

\begin{theorem}
[growth assumptions and Kurdyka-\L ojasiewicz
inequality]\label{cond-croiss} Let
$f:H\rightarrow\mathbb{R}\cup\{+\infty\}$ be a lower
semicontinuous proper convex function satisfying
(\ref{hyp-croissance}) and let us assume $0\in
C:=\mathrm{argmin}\,f$. Then the K\L --inequality holds,
\emph{i.e.}
\[
||\partial(\varphi\circ f)(x)||_{-}\geq1,\mbox{ for all }x\in
S\setminus \mbox{\rm argmin\,}f,
\]
with
\[
\varphi(r)=\int_{0}^{r}\frac{m^{-1}(s)}{s}ds.
\]
\end{theorem}

\noindent\textbf{Proof.} Let $x\in S\cap\mbox{\rm dom\,}\partial
f$ and $a$ be the projection of $x$ onto the convex subset
$C=\mbox{\rm argmin\,}f$. Using the
convex inequality we have
\[
f(x)-f(a)\leq\langle\partial^{0}f(x),x-a\rangle\leq\mbox{\rm
dist}\;(0,\partial f(x))\;\mbox{\rm dist}\;(x,C)\leq\mbox{\rm
dist}\;(0,\partial f(x))\;m^{-1}(f(x)-f(a)).
\]
Using the chain rule (see Lemma~\ref{chain}) an the fact that
$f(a)=0$, we obtain $\mbox{\rm dist}\;(0,\partial(\varphi\circ
f)(x))\geq1$ where $\varphi$ is as above (note that 
$\varphi\in \mathcal{K}(0,\rho)$).
$\hfill\Box$

\begin{remark}
Assume that $H$ is
finite-dimensional, and let $S$ be a compact convex subset of $H$
which satisfies $S\cap C\neq\emptyset$. Then there exists a convex
continuous increasing function
$m:\mathbb{R}_{+}\rightarrow\mathbb{R}_{+}$ with $m(0)=0$ such
that $f(x)\geq m(\mathrm{dist}(x,C))$ for all $x\in S$.

\noindent\lbrack\textit{Sketch of the proof}\textbf{.} With no
loss of generality we assume that $0\in S\cap C$. Using the
Moreau-Yosida regularization (see \cite{Brezis} for instance), we
obtain the existence of a finite-valued convex continuous function
$g:H\rightarrow\mathbb{R}$ such that $f\geq g$ and $\mbox{\rm
argmin\,}f=\mbox{\rm argmin\,}g$. Set $\alpha =\max\{\mbox{\rm
dist}\;(x,C):x\in S\}$ and $m_{0}(s)=\min\{g(x):x\in S,\;\mbox{\rm
dist}\;(x,C)\geq s\}\in\mathbb{R}_{+}$ for all $s\in
\lbrack0,\alpha]$. Let $0\leq s_{1}<s_{2}\leq\alpha$, and let
$x_{2}\in S$ be such that $\mbox{\rm dist}\;(x_{2},C)\geq s_{2}$
and $0<g(x_{2})=m(s_{2})$. Using the convexity of $g$ and the fact
that $0\in\mbox{\rm argmin\,}g\cap S$, we see that there exists
$\lambda\in(0,1)$ such that $g(\lambda x_{2})<g(x_{2})$, $\lambda
x_{2}\in S$ (recall that $S$ is convex and contains $0$), and
$\mbox{\rm dist}\;(\lambda x_{2},C)\geq s_{1}$. This shows that
the function $m_{0}$ is finite-valued increasing on $[0,\alpha]$
and satisfies $m_{0}(\mbox{\rm dist}\;(x,C))\leq g(x)\leq f(x)$
for any $x\in S$. Applying Lemma~\ref{minoration} (Annex) to
$m_{0}$, we obtain a smooth increasing finite-valued function $m$
such that $0<m(s)\leq m_{0}(s)$ for $s\in \lbrack0,\alpha]$ with
$m(0)=0$. The conclusion follows by extending $m$ to an increasing
continuous function on $\mathbb{R}_{+}$.]
\end{remark}

\begin{example}
\label{exple-croiss} Take $0< \alpha <1$. If $m(r)=\mathrm{exp}(-1/r^{\alpha})$ and $m(0)=0,$
then for $0\leq s\leq \rho <1$ we have
$m^{-1}
(s)=1/(-\,\mathrm{log}s)^{1/\alpha}$ and
$$\int_{0}^{\rho}\frac{m^{-1}(s)}{s}ds<+\infty.$$
 Therefore any convex
function which is minorized by the function $x\mapsto
\mathrm{exp}(-1/\mathrm{dist}(x,C)^{\alpha})$ in some neighborhood of
$C=\argmin f$ satisfies the K\L --inequality.
\end{example}

\subsection{A smooth convex counterexample to the K\L--inequality}

\label{sec:contre-exemple}

In this section we construct a $C^{2}$ convex function on
$\mathbb{R}^{2}$ with compact level sets that fails to satisfy the
K\L --inequality. This counterexample is constructed as follows:

\smallskip

- we first note that any sequence of sublevel sets of a convex
function that satisfies the K\L --inequality must comply with a
specific property ;

\smallskip

- we build a sequence $T_{k}$ of nested convex sets for which this
property fails ;

\smallskip

- we show that there exists a smooth convex function which admits
$T_{k}$ as sublevel sets.

\smallskip

The last part relies on the use of support functions and on a
result of Torralba \cite{Torralba}. For any closed convex subset
$T$ of $\mathbb{R}^{n} $, we define its support function by
$\sigma_{T}(x^{\ast})=\sup_{x\in T}\langle x,x^{\ast}\rangle$ for
all $x^{\ast}\in\mathbb{R}^{n}$. Let
$f:\mathbb{R}^{n}\rightarrow\mathbb{R}$ be a convex function and
$x^{\ast} \in\mathbb{R}^{n}$. Fenchel has observed, see
\cite{Fenchel}, that the function $\lambda\mapsto\sigma_{\lbrack
f\leq\lambda]}(x^{\ast})$ is concave and nondecreasing. The
following result asserts that this fact provides somehow a
sufficient condition to rebuild a convex function starting from a
countable family of nested convex sets.

\begin{theorem}
[Convex functions with prescribed level sets \cite{Torralba}]
\label{torralba} Let $\{T_{k}\}$ be a nonincreasing sequence of
convex compact subsets of $\mathbb{R}^{n}$ such that
$\mbox{int}\;T_{k}\supset T_{k+1}$ for all $k\geq 0$. For every $k>0$
we set:
\[
K_{k}=\max_{||x^{\ast}||=1}\frac{\sigma\newline
_{T_{k-1}}(x^{\ast}
)-\sigma_{T_{k}}(x^{\ast})}{\sigma_{T_{k}}(x^{\ast})-\sigma_{T_{k+1}}(x^{\ast
})}\in (0,+\infty).
\]
Then for every strictly decreasing sequence $\{\lambda_{k}\}$,
starting from $\lambda_{0}>0$ and satisfying
\[
0<K_{k}(\lambda_{k}-\lambda_{k+1})\leq\lambda_{k-1}-\lambda_{k},\mbox
{ for each }k>0,
\]
there exists a continuous convex function $f$ such that
\[
T_{k}=[f\leq\lambda_{k}],\quad\text{ for every }k\in\mathbb{N}
\]
and being maximal with this property.
\end{theorem}

\begin{remark}\ \\
\label{remTor}
(i) If $\{\lambda_{k}\}$ is as in the above
theorem and $x^{\ast}\in\mathbb{R}^{n}\backslash\{0\}$, we have
\[
\lambda_{k}-\lambda_{k+1}\leq\frac{\lambda_{0}-\lambda_{1}}{\sigma_{T_{0}
}(x^{\ast})-\sigma_{T_{1}}(x^{\ast})}(\sigma_{T_{k}}(x^{\ast})-\sigma
_{T_{k+1}}(x^{\ast})).
\]
Since the sum
$\sum(\sigma_{T_{k}}(x^{\ast})-\sigma_{T_{k+1}}(x^{\ast}))$
converges, so does the sum $\sum(\lambda_{k}-\lambda_{k+1})$,
yielding the existence of the limit $\lim\lambda_{k}$. Since $f$
is the greatest function admitting \{$T_{k}\}$ as prescribed
sublevel sets, we obtain $\min
f=\lim\lambda_{k}$.\smallskip\newline 
(ii) Let $k\geq0$
and $\lambda\in\lbrack\lambda_{k+1},\lambda_{k}]$. The function
$f$ satisfies further
\begin{equation}
\lbrack f\leq\lambda]=\left(
\frac{\lambda-\lambda_{k+1}}{\lambda_{k} -\lambda_{k+1}}\right)
T_{k}+\left(
\frac{\lambda_{k}-\lambda}{\lambda_{k}-\lambda_{k+1}}\right)
T_{k+1}, \label{deff}
\end{equation}
see \cite[Remark~5.9]{Torralba}.
\end{remark}

The following lemma provides a decreasing sequence of convex
compact subsets in $\mathbb{R}^{2}$ which can not be a
sequence of prescribed sublevel sets of a function satisfying the
K\L --inequality (see the {\it conclusion} part at the end of the proof of Theorem \ref{thm_counter}).

\begin{lemma}
\label{laurent}There exists a decreasing sequence of compact
subsets $\{T_{k}\}_{k}$ in $\mathbb{R}^{2}$ such that:

\begin{itemize}
\item [$(i)$]$T_{0}$ is the unit disk $D:=B(0,1)$ ;

\item[$(ii)$] $T_{k+1}\subset int\;T_{k}$ for every $k\in\mathbb{N}$ ;

\item[$(iii)$] $\displaystyle\bigcap_{k\in\mathbb{N}}T_{k}$ is the disk
$D_{r}:=B(0,r)$ for some $r>0$ ;

\item[$(iv)$] $\displaystyle\sum_{k=0}^{+\infty} \Dist(T_{k},T_{k+1})=+\infty$.
\end{itemize}
\end{lemma}

\noindent\textbf{Proof.} We proceed by constructing the boundaries
$\partial T_{k}$ of $T_{k}$ for each integer $k$. Let $C_{2,3}$
denote the circle of radius $1$ and let us define recursively a
sequence of closed convex curves $C_{n,m}$ for $n\geq3$ and $1\leq
m\leq n+1$; we assume that $C_{n-1,n}$ is the circle of radius
$R_{n}>0$. Let $\{\mu_{n}\}$ be a sequence in $(0,1)$ that will be
chosen later in order to satisfy $(iii)$. Then, for $1\leq m\leq
n$, let us define $C_{n,m}$ to be the union of the segments:

\begin{itemize}
\item [--]$\left[  \mu_{n}^{m}R_{n}\exp(2i\pi(\frac{j}{n})),\mu_{n}^{m}
R_{n}\exp(2i\pi(\frac{j+1}{n}))\right]  $ for $0\leq j\leq m-1$
(here $i$ stands for the imaginary unit) \newline and the
circle-arc:

\item[--] $\mu_{n}^{m}R_{n}\exp(i\theta)$ for $2\pi\frac{m}{n}\leq\theta
\leq2\pi$.
\end{itemize}

\begin{figure}[h]
\begin{center}
\resizebox{0.6\linewidth}{!}{\input{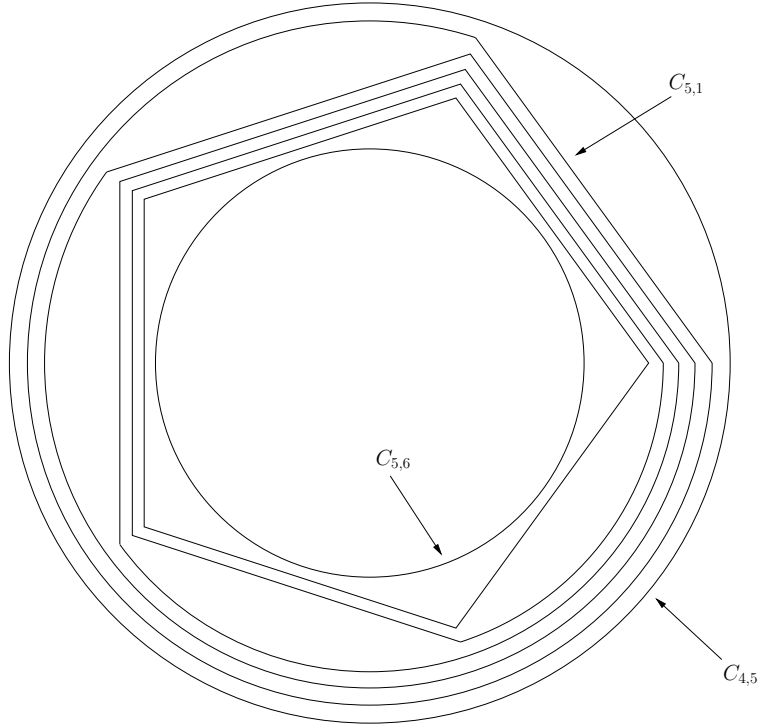}}
\end{center}
\caption{\textit{The curves $C_{4,5},C_{5,1}$ to $C_{5,6}$}}
\label{dessindomaine}
\end{figure}

In other words, $C_{n,m}$ consists of the first $m$ edges of a
regular convex $n$-gonon inscribed in a circle of radius
$\mu_{n}^{m}R_{n}$ and a circle-arc of the same radius to close
the curve. We then set
\[
R_{n+1}=\mu_{n}^{n+1}R_{n}\cos(\frac{\pi}{n})
\]
and define $C_{n,n+1}$ to be the circle of radius $R_{n+1}$.
Figure~\ref{dessindomaine} illustrates the curves $C_{4,5}$ and
$C_{5,m}$ for $m=1,\ldots,6$.

Ordering $\{(n,m):n\geq3,$ $1\leq m\leq n+1\}$ lexicographically
 we define succesively the convex subset
$T_{k}$ to be the convex envelope of the set $C_{n,m}$. 
By construction $(i)$ and $(ii)$ are
satisfied. Item $(iii)$ holds if $\lim R_{n}>0$ which is
equivalent to the fact that the infinite product $\displaystyle
\Pi_{n=3}^{+\infty}\mu_{n}^{n+1}\cos(\pi/n)$ does not converge to
$0$. This can be achieved by taking $\mu_{n}=1-1/n^{3}$. Let $r>0$
be the limit of $\{R_{n}\}$. The intersection of the convex sets
$T_{n}$ is the disk of radius $r$.

Take $n\geq3$. Considering the middle of the segment
$\displaystyle\left[
\mu_{n}R_{n},\mu_{n}R_{n}\exp(\frac{2i\pi}{n})\right]  $ in
$C_{n,1}$ and the point $R_{n}\exp(\frac{i\pi}{n})\in C_{n-1,n}$,
we obtain $\Dist (C_{n,1},C_{n-1,n})=R_{n}(1-\mu_{n}\cos(\pi/n))$.
If $2\leq m\leq n$, considering the middle of
$\displaystyle\left[  \mu_{n}^{m}R_{n}\exp
(\frac{2i\pi(m-1)}{n},\mu_{n}^{m}R_{n}\exp(\frac{2i\pi
m}{n})\right]  $ in $C_{n,m}$ and the point
$\mu_{n}^{m-1}R_{n}\exp(\frac{i\pi(2m-1)}{n})\in C_{n,m-1}$, we
get $\Dist(C_{n,m},C_{n,m-1})=\mu_{n}^{m-1}R_{n}(1-\mu_{n}
\cos(\pi/n))$. Finally considering the points $\mu_{n}^{n}R_{n}\in
C_{n,n}$ and $\mu_{n}^{n+1}\cos(\pi/n)R_{n}\in C_{n,n+1}$, we
obtain

\[
\Dist (C_{n,n},C_{n,n+1})=\mu_{n}^{n}R_{n}(1-\mu_{n}\cos(\pi/n)).
\]
Thus
\[
\Dist(C_{n,1},C_{n-1,n})+\sum_{m=2}^{n+1}\Dist(C_{n,m},C_{n,m-1})=\sum
_{m=1}^{n+1}\mu_{n}^{m-1}R_{n}(1-\mu_{n}\cos\frac{\pi}{n})\sim
nr\frac{\pi ^{2}}{2n^{2}}=\frac{\pi^{2}r}{2n}.
\]
Hence $(iv)$ holds. $\hfill\Box$

\bigskip
For $\theta\in \mathbb{R} /2\pi\mathbb{Z}$,
 set $n(\theta)=(\cos\theta,\sin \theta)$ and $\tau(\theta)=(-\sin\theta,\cos\theta)$. 
We say that a closed $C^{2}$ \emph{curve} $C$ in $\mathbb{R}^{2}$
is convex if its curvature has constant sign. 
If moreover the
curvature never vanishes, then there exists a $C^{1}$
\emph{parametrization} $c:\mathbb{R} /2\pi\mathbb{Z}\rightarrow C$
of $C$, called \emph{parametrization of $C$ by its normal}, such
that the unit tangent vector at $c(\theta)$ is $\tau (\theta)$. In
this case $n(\theta)$ is the outward normal to the convex envelope
of $C$ at $c(\theta)$. Moreover, $c$ is $C^{\infty}$, whenever
$C$ is so. In this case, we denote by $\rho_{c}(\theta)$ the
curvature radius of $c$ at $c(\theta)$ and we have
\[
\dot{c}(\theta)=\rho_{c}(\theta)\tau(\theta).
\]

Let us denote by $T$ the convex envelope of $C$. Using the fact
that $n$ defines the outward normals to $T$, we get
\[
\langle c(\theta),n(\theta)\rangle=\max_{x\in T}\langle x,n(\theta
)\rangle=\sigma_{T}(n(\theta)),\;\forall\theta\in\mathbb{R}/2\pi\mathbb{Z}.
\]

\begin{theorem}
[convex counterexample]\label{thm_counter}There exists a $C^{2}$
convex function $f:\mathbb{R}^{2}\rightarrow\mathbb{R}$ such that $\min f=0$ which does
not satisfy the K\L --inequality and whose set of minimizers is
compact with nonempty interior. 
More precisely, for each $r>0$ and for each desingularization function
$\varphi \in {\cal K}(0,r)$ we have
\[
\inf\left\{\Vert\nabla(\varphi\circ f)(x)\Vert : x\in[0<f<r]\right\}=0.
\]
\end{theorem}

\begin{remark}\ \\
\label{Remark38}
(i) It can be seen from the forthcoming
proof that $\mbox{\rm argmin\,}f$ is the closed disk centered at
$0$ of radius $r$, and that $f$ is actually $C^{\infty}$ on the
complement of the circle of radius $r$.\smallskip\newline
(ii) The fact that $f$ is $C^{2}$ shows that K\L
--inequality is not related to the smoothness of $f$. Besides, it
seems clear from the proof that a $C^{k}$ ($k$ arbitrary)
counterexample could be obtained.\smallskip\newline 
(iii)
Since $\mbox{\rm argmin\,}f$ has nonempty interior,
Theorem~\ref{argmin-nonvide} shows that the lengths of subgradient
curves are uniformly bounded. Using the notation and the results
of Theorem~\ref{global}, we see that the function $f$ shows that
the uniform boundedness of the lengths of the subgradient curves
(starting from a given level set $[f=r_{0}]$) does not
yield the uniform boundedness of the lengths of the piecewise
subgradient curves $\gamma$ lying in $[\min f<f<r_{0}]\}$.
\end{remark}

\bigskip

\noindent\textbf{Proof of Theorem \ref{thm_counter}.} Let $M,N$ be
 topological finite-dimensional manifolds. In this proof, a mapping
$F:M\rightarrow N$ is said to be \emph{proper} if for each compact
subset $K$ of $N$, $F^{-1}(K)$ is a compact subset of $M$.

\smallskip

\noindent\emph{Smoothing the sequence $T_{k}$.} Let us consider a
sequence of convex compact sets $\{T_{k}\}$ as in
Lemma~\ref{laurent}. Set $C_{k}=\partial T_{k}$ and consider a
positive sequence $\epsilon_{k}$ such that $\sum
\epsilon_{k}<+\infty$ with $\epsilon_{k}+\epsilon_{k+1}<\Dist
(T_{k},T_{k+1})=\Dist(C_{k},C_{k+1})$ for each integer $k$. The
$\epsilon_{k} $-neighborhood of $C_{k}$ can be seen to be disjoint
from the $\epsilon _{k^{\prime}}$-neighborhood of $C_{k^{\prime}}$
whenever $k\neq k^{\prime}$. We can deform $C_{k}$ into a
$C^{\infty}$ convex closed curve $\widetilde {C}_{k}$ whose
curvature never vanishes, lying in the $\epsilon_{k}
$-neighborhood of $C_{k}$. This smooth deformation can be achieved
by letting $C_{k}$ evolve under the mean-curvature flow during a
very short time, see \cite{Evansspruck3} for the smoothing aspects
and \cite{GageHamilton,Zhu} for the positive curvature results. We
set $\widetilde{T}_{k}$ to be the closed convex envelope of
$\widetilde{C}_{k}$. This process yields a decreasing sequence of
compact convex sets $\{\widetilde{T}_{k}\}$, that satisfies the
conditions of Lemma~\ref{laurent}. We note that the circle of
radius $1$ has non-zero curvature and we set
$C_{0}=\widetilde{C}_{0}$. Since $\Dist
(\widetilde{T}_{k},\widetilde{T}_{k+1})\geq\Dist
(T_{k},T_{k+1})-(\epsilon_{k}+\epsilon_{k+1})$ and
$\sum\epsilon_{k}<+\infty$, condition $(iv)$ holds. With no loss
of generality we may therefore assume that for each $k\geq0$ the
curve $\partial T_{k}$ is smooth and can be parametrized by its
normal.

\smallskip

\noindent Let $K_{k}$ be as in Theorem~\ref{torralba}, let $\lambda_{0}$ and
$\lambda_{1}$ be such that $\lambda_{0}>\lambda_{1}$. We define
$\lambda_{k}$ recursively by
\begin{equation}
K_{k}(\lambda_{k}-\lambda_{k+1})=\frac{1}{2}(\lambda_{k-1}-\lambda
_{k}).\label{deflambda}
\end{equation}
Because of \eqref{deflambda}, Theorem~\ref{torralba} yields a
continuous convex function $f:T_{0}\rightarrow\mathbb{R}$ such
that $T_{k}=[f\leq \lambda_{k}]$. Since $f$ is the greatest
function with this property, we deduce that $\min
f=\lim\lambda_{k}$ and $\mbox{\rm argmin\,}f=\cap
_{k\in\mathbb{N}}T_{k}$.

\smallskip

\noindent\emph{Smoothing the function $f$ on
$\mathbb{R}^{n}\setminus\mbox{\rm argmin\,}f$}. We can easily
extend $f$ outside $T_{0}$ into a smooth convex function. Let us
examine the restriction of $f$ to $T_{0}$. Since $\partial T_{k}$
can be parametrized by its normal, we denote by $c_{k}:\mathbb{R}
/2\pi\mathbb{Z}\rightarrow\mathbb{R}^{2}$ this parametrization.
Let us fix $k\in\mathbb{N}$. Let $\theta$ be in
$\mathbb{R}/2\pi\mathbb{Z}$. Using Remark~\ref{remTor} (b), we
obtain
\begin{align*}
\max_{x\in\lbrack f\leq\lambda]}\langle x,n(\theta)\rangle &
=\left(
\frac{\lambda-\lambda_{k+1}}{\lambda_{k}-\lambda_{k+1}}\right)
\max_{x\in T_{k}}\langle x,n(\theta)\rangle+\left(
\frac{\lambda_{k}-\lambda} {\lambda_{k}-\lambda_{k+1}}\right)
\max_{x\in T_{k+1}}\langle x,n(\theta )\rangle\\ &  =\left(
\frac{\lambda-\lambda_{k+1}}{\lambda_{k}-\lambda_{k+1}}\right)
\langle c_{k}(\theta),n(\theta)\rangle+\left(
\frac{\lambda_{k}-\lambda }{\lambda_{k}-\lambda_{k+1}}\right)
\langle c_{k+1}(\theta),n(\theta )\rangle\\ &  =\left\langle
\left(  \frac{\lambda-\lambda_{k+1}}{\lambda_{k}
-\lambda_{k+1}}\right)  c_{k}(\theta)+\left(
\frac{\lambda_{k}-\lambda }{\lambda_{k}-\lambda_{k+1}}\right)
c_{k+1}(\theta),n(\theta)\right\rangle .
\end{align*}
Using \eqref{deff} once more we obtain

\begin{equation}
\left(\frac{\lambda-\lambda_{k+1}}{\lambda_{k}-\lambda_{k+1}}\right)
c_{k}(\theta)+\left(
\frac{\lambda_{k}-\lambda}{\lambda_{k}-\lambda_{k+1} }\right)
c_{k+1}(\theta)\in\lbrack f\leq\lambda]. \label{frondG}
\end{equation}
\noindent
Since the above maximum is achieved in $[f=\lambda]$, it follows
that
\begin{equation}
f\left(  \left(
\frac{\lambda-\lambda_{k+1}}{\lambda_{k}-\lambda_{k+1} }\right)
c_{k}(\theta)+\left(  \frac{\lambda_{k}-\lambda}{\lambda_{k}
-\lambda_{k+1}}\right)  c_{k+1}(\theta)\right)  =\lambda.
\end{equation}
Let us define
$G:\mathbb{R}\times\mathbb{R}/2\pi\mathbb{Z}\rightarrow
\mathbb{R}^{2}$ by
\[
G(\lambda,\theta)=\left(  \frac{\lambda-\lambda_{k+1}}{\lambda_{k}
-\lambda_{k+1}}\right)  c_{k}(\theta)+\left(
\frac{\lambda_{k}-\lambda }{\lambda_{k}-\lambda_{k+1}}\right)
c_{k+1}(\theta).
\]
The map $G$ is clearly $C^{\infty}$. Since $\dfrac{\partial
G}{\partial
\lambda}=\frac{c_{k}(\theta)-c_{k+1}(\theta)}{\lambda_{k}-\lambda_{k+1}}$,
we have
\begin{align*}
\langle\dfrac{\partial
G}{\partial\lambda},n(\theta)\rangle=\langle\frac
{c_{k}(\theta)-c_{k+1}(\theta)}{\lambda_{k}-\lambda_{k+1}},n(\theta)\rangle
& =\frac{\langle c_{k}(\theta),n(\theta)\rangle-\langle
c_{k+1}(\theta ),n(\theta)\rangle}{\lambda_{k}-\lambda_{k+1}}\\ &
=\frac{\max_{x\in T_{k}}\langle x,n(\theta)\rangle-\max_{x\in
T_{k+1} }\langle x,n(\theta)\rangle}{\lambda_{k}-\lambda_{k+1}}\\
&  >0.
\end{align*}
On the other hand
\begin{equation}
\dfrac{\partial G}{\partial\theta}=\left(  \left(
\frac{\lambda-\lambda _{k+1}}{\lambda_{k}-\lambda_{k+1}}\right)
\rho_{c_{k}}(\theta)+\left(
\frac{\lambda_{k}-\lambda}{\lambda_{k}-\lambda_{k+1}}\right)
\rho_{c_{k+1} }(\theta)\right)  \tau(\theta). \label{dGdtheta}
\end{equation}
Since $\rho_{c_{k}}>0$ and $\rho_{c_{k+1}}>0$, $G$ is a local
diffeomorphism on
$(\lambda_{k+1}-\delta,\lambda_{k}+\delta)\times\mathbb{R}/2\pi\mathbb{Z}$
for any $\delta>0$ sufficiently small. In view of \eqref{frondG},
we have $G(\lambda,\theta)\in\lbrack\lambda_{k+1}\leq
f\leq\lambda_{k}]$ for $\lambda_{k+1}\leq\lambda\leq\lambda_{k}$
and $G(\lambda,\theta)\in \lbrack\lambda_{k+1}<f<\lambda_{k}]$ for
$\lambda_{k+1}<\lambda<\lambda_{k}$. Since the map
$\widetilde{G}:[\lambda_{k+1},\lambda_{k}]\times\mathbb{R}
/2\pi\mathbb{Z}\rightarrow\lbrack\lambda_{k+1}\leq
f\leq\lambda_{k}]$ defined by
$\widetilde{G}(\lambda,\theta)=G(\lambda,\theta)$ is proper,
$\widetilde {G}$ is a covering map from
$[\lambda_{k+1},\lambda_{k}]\times\mathbb{R} /2\pi\mathbb{Z}$ to
$[\lambda_{k+1}\leq f\leq\lambda_{k}]$. The set
$[\lambda_{k+1}\leq f\leq\lambda_{k}]$ is connected, thus
$\widetilde{G}$ is onto. Using \eqref{frondG} and
$G(\lambda_{k},\theta)=c_{k}(\theta)$, one sees that
$(\lambda_{k},\theta)$ is the only antecedent of $c_{k}(\theta)$
by $\widetilde{G}$ and, since
$[\lambda_{k+1},\lambda_{k}]\times\mathbb{R} /2\pi\mathbb{Z}$ is
connected, $\widetilde{G}$ is injective. Thus $\widetilde{G}$ is a
$C^{\infty}$ diffeomorphism (see \cite[Proposition 2.19]{Lee}). By
\eqref{frondG}, this implies that the restriction of $f$ to
$[\lambda_{k+1}\leq f\leq\lambda_{k}]$ is $C^{\infty}$. Using
\eqref{frondG}, we know that the level line $[f=\lambda]$ (for
$\lambda_{k+1}\leq\lambda \leq\lambda_{k}$) is parametrized by
$G(\lambda,\theta)$ for $\theta \in\mathbb{R}/2\pi\mathbb{Z}$; if
$c_{\lambda}$ denotes this parametrization, then
$c_{k}=c_{\lambda_{k}}$. Besides, by \eqref{dGdtheta},
$c_{\lambda}$ is a parametrization by the normal and
$\rho_{c_{\lambda}}$ is a convex combination of $\rho_{c_{k}}$ and
$\rho_{c_{k+1}}$, hence $\rho_{c_{\lambda}}>0$.

\smallskip

Let us compute $\nabla f$ at $c_{\lambda}(\theta)$. Equation
\eqref{frondG} yields $1=\langle\nabla
f(G(\lambda,\theta)),\dfrac{\partial G}{\partial
\lambda}(\lambda,\theta)\rangle$. Besides we also know that the
normal to $[f=\lambda]$ at $c_{\lambda}(\theta)$ is $n(\theta)$.
Since the gradient $\nabla f(G(\lambda,\theta))$ and the normal  
$n(\theta)$ are linearly dependent, we obtain
\begin{equation}
\nabla
f(c_{\lambda}(\theta))=\frac{\lambda_{k}-\lambda_{k+1}}{\langle
c_{\lambda_{k}}(\theta)-c_{\lambda_{k+1}}(\theta),n(\theta)\rangle}n(\theta).
\label{gradient}
\end{equation}
Note that this expression does not depend on
$\lambda\in\lbrack\lambda _{k+1}-\lambda_{k}]$.

\smallskip

Before going further let us observe/recall two facts.

\smallskip

-- First using the aforementioned result of Fenchel
\cite{Fenchel}, we deduce from the convexity of $f$ that the
function
\begin{equation}\label{clambdaconcave}
\lambda\mapsto\langle
c_{\lambda}(\theta),n(\theta)\rangle=\sigma_{\lbrack
f\leq\lambda]}(n(\theta))\text{ is concave and increasing}.
\end{equation}

-- Let $\lambda$ and $\lambda^{\prime}$ be such that
$\lambda_{k+1}\leq \lambda\leq\lambda^{\prime}\leq\lambda_{k}$. We
have :
\begin{align}
c_{\lambda}(\theta) &  =\left(
\frac{\lambda-\lambda_{k+1}}{\lambda^{\prime
}-\lambda_{k+1}}\right)  c_{\lambda^{\prime}}(\theta)+\left( \frac
{\lambda^{\prime}-\lambda}{\lambda^{\prime}-\lambda_{k+1}}\right)
c_{\lambda_{k+1}}(\theta),\label{moyenne}\\
c_{\lambda^{\prime}}(\theta) &  =\left(
\frac{\lambda^{\prime}-\lambda }{\lambda_{k}-\lambda}\right)
c_{\lambda_{k}}(\theta)+\left(  \frac
{\lambda_{k}-\lambda^{\prime}}{\lambda_{k}-\lambda}\right)
c_{\lambda} (\theta).
\end{align}

(\emph{Smoothing $f$ around $[f=\lambda_{k}]$}.) We have seen that
the function $f$ is $C^{\infty}$ on the complement of the union of
the level lines $[f=\lambda_{k}]$ for
$k\in\mathbb{N}$. In order to go further we need to
modify $f$ around each $[f=\lambda_{k}]$.

\smallskip

Consider a positive sequence $\{\epsilon_{k}\}$ such that
$\sum_{i} \epsilon_{i}<+\infty$ and
$\epsilon_{k}+\epsilon_{k+1}<\Dist(T_{k}
,T_{k+1})=\Dist([f=\lambda_{k}],[f=\lambda_{k+1}])$ for each
integer $k$. Let us assume that there exists a sequence
$f_{k}:\mathbb{R}^{2}\rightarrow \mathbb{R}$ of convex functions
such that:

\begin{enumerate}
\item [P1]$f_{0}=f$ ;

\item[P2] $f_{k}=f_{k-1}$ outside an $\epsilon_{k}$-neighborhood of
$[f=\lambda_{k}]$ ;

\item[P3] $f_{k}$ is $C^{\infty}$ in $[f>\lambda_{k+1}$] ;

\item[P4] $\Vert\nabla f_{k}\Vert$ is bounded in $[f\leq\lambda_{k}]$ by the
maximum of $\Vert\nabla f\Vert$ in $[\lambda_{k}\leq
f\leq\lambda_{k-1}].$
\end{enumerate}

Let us choose $k\geq1$ and $\lambda,\lambda^{\prime}$ such that
$\lambda_{k+1}\leq\lambda\leq\lambda_{k}\leq\lambda^{\prime}\leq\lambda_{k-1}$.
Then by \eqref{deflambda} and \eqref{gradient} we have:
\[
\Vert\nabla
f(c_{\lambda}(\theta))\Vert=\frac{\lambda_{k}-\lambda_{k+1}
}{\langle
c_{\lambda_{k}}(\theta)-c_{\lambda_{k+1}}(\theta),n(\theta)\rangle
}\leq\frac{1}{2}\frac{\lambda_{k-1}-\lambda_{k}}{\langle
c_{\lambda_{k-1}
}(\theta)-c_{\lambda_{k}}(\theta),n(\theta)\rangle}=\frac{1}{2}\Vert\nabla
f(c_{\lambda^{\prime}}(\theta)\Vert.
\]
Hence
\begin{equation}
\max_{\lbrack\lambda_{k+1}\leq f\leq\lambda_{k}]}\Vert\nabla
f\Vert\leq \frac{1}{2}\max_{[\lambda_{k}\leq
f\leq\lambda_{k-1}]}\Vert\nabla f\Vert. \label{majgradient}
\end{equation}
Combining with (P4), the above implies that the sequence
$(f_{k})_{k\in\mathbb{N}}$ is uniformly Lipschitz continuous.
Applying Ascoli compactness theorem we obtain that $f_{k}$
converge to a continuous function $\tilde{f}$  which is convex. From (P2) and (P3), we obtain
successively that $\tilde{f}$ has the same set of minimizers as
$f$, $f$ is $C^{\infty}$ outside $\mbox{\rm argmin\,}\tilde{f}$,
$[\tilde{f} =\lambda_{k}]$ is in the $\epsilon_{k}$-neighborhood
of $[f=\lambda_{k}]$. Moreover \eqref{majgradient} and (P4) imply
that $\Vert\nabla\tilde{f} (x)\Vert$ goes to zero as $x$
approaches $\mbox{\rm argmin\,}\tilde{f}$, hence $\tilde{f}$ is
globally $C^{1}$. Note also, that the sequence of level sets
$[\tilde{f}\leq\lambda_{k}]$ satisfies the hypothesis $(iv)$ of
Lemma~\ref{laurent}. As shown in the conclusion, $\tilde{f}$ provides a  $C^{1}$ counterexample
to the K\L --inequality.

\smallskip

Let us define such a sequence $\{f_{k}\}$ by induction. Assume
that $f_{k-1}$ is defined. In order to construct $f_{k}$, it
suffices to proceed in the $\epsilon_{k}$-neighborhood of
$[f=\lambda_{k}]$. Let $\epsilon>0$ such that
$[\lambda_{k}-2\epsilon\leq f\leq\lambda_{k}+2\epsilon]$ is in the
$\epsilon_{k}$-neighborhood of $[f=\lambda_{k}]$. Let us consider
a $C^{\infty}$ function
$\mu_{-}:[-2\epsilon,2\epsilon]\rightarrow\mathbb{R}$ which
satisfies the following properties:
\begin{align*}
1. &  \ \mu_{-}\text{ is nonincreasing,} & 2. &  \
\mu_{-}^{\prime\prime} \geq0,\\ 3. &  \
\mu_{-}(\lambda)=-\lambda/\epsilon\text{ on
}[-2\epsilon,-\epsilon/2], & \ 4. &  \ \mu_{-}(\lambda)=0\text{ on
}[\epsilon/2,2\epsilon].
\end{align*}
Let us then define
$\mu_{+}(\lambda):=\lambda/\epsilon+\mu_{-}(\lambda)$ and
$\mu_{0}=1-(\mu_{-}+\mu_{+})$. The function $\mu_{+}$ satisfies
\begin{align*}
1^{\prime}. &  \ \mu_{+}\text{ is nondecreasing,} & 2^{\prime}. &
\ \mu _{+}^{\prime\prime}=\mu_{-}^{\prime\prime}\geq0,\\
3^{\prime}. &  \ \mu_{+}(\lambda)=0\text{ on
}[-2\epsilon,-\epsilon/2], & 4^{\prime}. &  \
\mu_{+}(\lambda)=\lambda/\epsilon\text{ on }[\epsilon
/2,2\epsilon].
\end{align*}
Set $c_{-}=c_{\lambda_{k}-\epsilon}$, $c_{0}=c_{\lambda_{k}}$,
$c_{+} =c_{\lambda_{k}+\epsilon}$ and
\begin{align*}
M_{-}(\theta) &  =\langle
c_{-}(\theta),n(\theta)\rangle=\max_{x\in\lbrack
f\leq\lambda_{k}-\epsilon]}\langle x,n(\theta)\rangle,\\
M_{0}(\theta) &  =\langle
c_{0}(\theta),n(\theta)\rangle=\max_{x\in\lbrack
f\leq\lambda_{k}]}\langle x,n(\theta)\rangle,\\ M_{+}(\theta) &
=\langle c_{+}(\theta),n(\theta)\rangle=\max_{x\in\lbrack
f\leq\lambda_{k}+\epsilon]}\langle x,n(\theta)\rangle.
\end{align*}
For
$(\lambda,\theta)\in\lbrack-2\epsilon,2\epsilon]\times\mathbb{R}
/2\pi\mathbb{Z}$, we define:
\[
H(\lambda,\theta)=\mu_{-}(\lambda)c_{-}(\theta)+\mu_{0}(\lambda)c_{0}
(\theta)+\mu_{+}(\lambda)c_{+}(\theta).
\]
Then $H$ is a $C^{\infty}$ map and for any
$\lambda\in\lbrack-\epsilon ,\epsilon]$, we have
$\mu_{-}(\lambda),\mu_{0}(\lambda)$ and $\mu_{+} (\lambda)$ in
$[0,1]$. Since $H(\lambda,\theta)$ is a convex combination of
points in $[f\leq\lambda_{k}+\epsilon]$, we deduce
$H(\lambda,\theta )\in\lbrack f\leq\lambda_{k}+\epsilon]$ and
$H(\lambda,\theta)\in\lbrack f<\lambda_{k}+\epsilon]$ whenever
$\lambda<\epsilon$ and $\mu_{+}(\lambda)<1$. Since
\[
\langle
H(\lambda,\theta),n(\theta)\rangle=\mu_{-}(\lambda)M_{-}(\theta
)+\mu_{0}(\lambda)M_{0}(\theta)+\mu_{+}(\lambda)M_{+}(\theta)\geq
M_{-} (\theta),
\]
we get $H(\lambda,\theta)\in\lbrack f\geq\lambda_{k}-\epsilon]$,
and $H(\lambda,\theta)\in\lbrack f>\lambda_{k}-\epsilon]$ whenever
$\lambda
>\epsilon$, $\mu_{-}(\lambda)<1$. It follows that
\[
\dfrac{\partial
H}{\partial\lambda}=\mu_{-}^{\prime}(\lambda)c_{-}(\theta
)+\mu_{0}^{\prime}(\lambda)c_{0}(\theta)+\mu_{+}^{\prime}(\lambda)c_{+}
(\theta).
\]
Since $\mu_{0}^{\prime}=-\mu_{-}^{\prime}-\mu_{+}^{\prime}$, items
$1$ and $1^{\prime}$ entail
\begin{align*}
\langle\dfrac{\partial H}{\partial\lambda},n(\theta)\rangle &
=\mu _{+}^{\prime}(\lambda)\langle
c_{+}(\theta)-c_{0}(\theta),n(\theta)\rangle
-\mu_{-}^{\prime}(\lambda)\langle
c_{0}(\theta)-c_{-}(\theta),n(\theta )\rangle\\ &
=\mu_{+}^{\prime}(\lambda)(M_{+}(\theta)-M_{0}(\theta))-\mu_{-}^{\prime
}(\lambda)(M_{0}(\theta)-M_{-}(\theta))\\ &  >0.
\end{align*}
On the other hand
\begin{equation}\label{dHdtheta}
\dfrac{\partial H}{\partial\theta}=\left(
\mu_{-}(\lambda)\rho_{c_{-}}
(\theta)+\mu_{0}(\lambda)\rho_{c_{0}}(\theta)+\mu_{+}(\lambda)\rho_{c_{+}
}(\theta)\right)  \tau(\theta),
\end{equation}
so that $\langle\dfrac{\partial
H}{\partial\theta},n(\theta)\rangle=0$ and $\langle\dfrac{\partial
H}{\partial\theta},\tau(\theta)\rangle>0$ for
$\lambda\in]-\epsilon^{\prime},\epsilon^{\prime}[$ with
$\epsilon^{\prime }>\epsilon$. Thus $H$ is a local diffeomorphism
on $]-\epsilon^{\prime
},\epsilon^{\prime}[\times\mathbb{R}/2\pi\mathbb{Z}$. The map
$\widetilde
{H}:[-\epsilon,\epsilon]\times\mathbb{R}/2\pi\mathbb{Z}\rightarrow
\lbrack\lambda_{k}-\epsilon\leq f\leq\lambda_{k}+\epsilon]$
defined by $\widetilde{H}(\lambda,\theta)=H(\lambda,\theta)$ is
proper, therefore $\widetilde{H}$ is a covering map from
$[-\epsilon,\epsilon]\times \mathbb{R}/2\pi\mathbb{Z}$ to
$[\lambda_{k}-\epsilon\leq f\leq\lambda _{k}+\epsilon]$. Since
$[\lambda_{k}-\epsilon\leq f\leq\lambda_{k}+\epsilon]$ is
connected, $\widetilde{H}$ is onto. Besides, since
$c_{+}(\theta)\in\lbrack f=\lambda_{+}\epsilon]$,
$(\epsilon,\theta)$ is the only antecedent of $c_{+}(\theta)$ by
$H$, $\widetilde{H}$ is injective by connectedness of
$[-\epsilon,\epsilon]\times\mathbb{R}/2\pi\mathbb{Z}$.
$\widetilde{H}$ is therefore a $C^{\infty}$ diffeomorphism from
$[-\epsilon,\epsilon ]\times\mathbb{R}/2\pi\mathbb{Z}$ into
$[\lambda_{k}-\epsilon\leq f\leq \lambda_{k}+\epsilon]$.

We then define $f_{k}$ to be $f_{k-1}$ outside of
$[\lambda_{k}-\epsilon\leq f\leq\lambda_{k}+\epsilon]$ and by
$f_{k}(H(\lambda,\theta))=\lambda _{k}+\lambda$ in
$[\lambda_{k}-\epsilon\leq f\leq\lambda_{k}+\epsilon]$. When
$\lambda\in\lbrack\lambda_{k}-\epsilon,\lambda_{k}-\epsilon/2]$,
Properties $3$, $3^{\prime}$ and equation \eqref{moyenne} yield
\begin{align*}
H(\lambda-\lambda_{k},\theta) &
=-\frac{\lambda-\lambda_{k}}{\epsilon}
c_{-}(\theta)+(1+\frac{\lambda-\lambda_{k}}{\epsilon})c_{0}(\theta)\\
&
=\frac{\lambda_{k}-\lambda}{\lambda_{k}-(\lambda_{k}-\epsilon)}c_{-}
(\theta)+\frac{\lambda-(\lambda-\epsilon)}{\lambda_{k}-(\lambda_{k}-\epsilon
)}c_{0}(\theta)\\ &  =c_{\lambda}(\theta).
\end{align*}
Thus $f_{k}=f=f_{k-1}$ in $[\lambda_{k}-\epsilon\leq
f\leq\lambda_{k} -\epsilon/2]$ and for similar reasons
$f_{k}=f_{k-1}$ in $[\lambda _{k}+\epsilon/2\leq
f\leq\lambda_{k}+\epsilon]$. The ``gluing'' of $f_{k-1}$ and
$f_{k}$ is therefore $C^{\infty}$ along $[f=\lambda_{k}-\epsilon]$
and $[f=\lambda_{k}+\epsilon]$. Hence, $f_{k}$ satisfies (P3).

\smallskip

Let us compute $\nabla f_{k}$ in $[\lambda_{k}-\epsilon\leq
f\leq\lambda _{k}+\epsilon]$. By definition of $f_{k}$,
$1=\langle\nabla f_{k} (H(\lambda,\theta)),\dfrac{\partial
H}{\partial\lambda}\rangle$. Besides
$H(\lambda-\lambda_{k},\theta)$ is a parametrization of the level
line $[f_{k}=\lambda]$ by its normal (see \eqref{dHdtheta}), hence
$\nabla f_{k}(H(\lambda,\theta))=\alpha n(\theta)$ with
$\alpha>0$. Using both formulae, we finally get
\[
\nabla
f_{k}(H(\lambda,\theta))=\frac{1}{\mu_{+}^{\prime}(\lambda)\langle
c_{+}(\theta)-c_{0}(\theta),n(\theta)\rangle-\mu_{-}^{\prime}(\lambda)\langle
c_{0}(\theta)-c_{-}(\theta),n(\theta)\rangle}\;n(\theta).
\]
From the definition of $\mu_{+}$,
$\mu_{+}^{\prime}(\lambda)-\mu_{-}^{\prime }(\lambda)=1/\epsilon$.
Besides, for $\lambda\in\lbrack-\epsilon,-\epsilon/2]$ we have
\[
\frac{\epsilon}{\langle
c_{0}(\theta)-c_{-}(\theta),n(\theta)\rangle} =\Vert\nabla
f(c_{\lambda+\lambda_{k}}(\theta))\Vert,
\]
while for $\lambda\in\lbrack\epsilon/2,\epsilon]$ we get
\[
\frac{\epsilon}{\langle
c_{+}(\theta)-c_{0}(\theta),n(\theta)\rangle} =\Vert\nabla
f(c_{\lambda+\lambda_{k}}(\theta))\Vert.
\]
Hence by \eqref{clambdaconcave}:
\[
\Vert\nabla f_{k}(H(\lambda,\theta))\Vert\leq\Vert\nabla
f(c_{\lambda _{k}+\epsilon}(\theta))\Vert.
\]
(P4) is therefore satisfied.

\smallskip

The last assertion we need to establish is the convexity of $f_{k}$.
By construction, it suffices to prove that the Hessian $Q_{f_{k}}$
of $f$ is nonnegative in $[\lambda_{k}-\epsilon\leq
f\leq\lambda_{k}+\epsilon]$. Let us denote by $Q_{H}$ the Hessian
of $H$ (observe that $Q_{H}$ takes its values in
$\mathbb{R}^{2}$). For $-\epsilon\leq\lambda\leq\epsilon$, we have
$\lambda+\lambda_{k}=f_{k}(H(\lambda,\theta))$, thus
\[
0=Q_{f_{k}}(H(\lambda,\theta))(DH(\lambda,\theta)(\cdot),DH(\lambda
,\theta)(\cdot))+\langle\nabla
f_{k}(H(\lambda,\theta)),Q_{H}(\lambda
,\theta)(\cdot,\cdot)\rangle
\]
where $DH$ denotes the differential map of $H$. To prove that
$Q_{f_{k}}$ is nonnegative, it suffices to prove that
$\langle\nabla f_{k}(H(\lambda
,\theta)),Q_{H}(\lambda,\theta)(\cdot,\cdot)\rangle\leq0$. We
have
\begin{align*}
\dfrac{\partial^{2}H}{\partial\lambda^{2}} &
=\mu_{-}^{\prime\prime}
(\lambda)c_{-}(\theta)+\mu_{0}^{\prime\prime}(\lambda)c_{0}(\theta)+\mu
_{+}^{\prime\prime}(\lambda)c_{+}(\theta)\\ &
=\mu_{-}^{\prime\prime}(\lambda)(c_{-}(\theta)-c_{0}(\theta))+\mu
_{+}^{\prime\prime}(\lambda)(c_{+}(\theta)-c_{0}(\theta))\\ &
=\mu_{+}^{\prime\prime}(\lambda)\big((c_{+}(\theta)-c_{0}(\theta
))-(c_{0}(\theta)-c_{-}(\theta))\big),
\end{align*}
where the last equality is due to item $2^{\prime}$. On the other
hand
\[
\langle\nabla
f_{k}(H(\lambda,\theta)),\dfrac{\partial^{2}H}{\partial
\lambda^{2}}\rangle=\mu_{+}^{\prime\prime}(\lambda)\Vert\nabla
f_{k} (H(\lambda,\theta)\Vert\big(\langle
c_{+}(\theta)-c_{0}(\theta),n(\theta )\rangle-\langle
c_{0}(\theta)-c_{-}(\theta),n(\theta)\rangle\big)
\]
which is nonpositive because of \eqref{clambdaconcave}. Besides we
have
\[
\dfrac{\partial^{2}H}{\partial\lambda\partial\theta}=\big(\mu_{-}^{\prime
}(\lambda)\rho_{c_{-}}(\lambda)+\mu_{0}^{\prime}(\lambda)\rho_{c_{0}}
(\lambda)+\mu_{+}^{\prime}(\lambda)\rho_{c_{+}}(\lambda)\big)\tau(\theta),
\]
thus $\langle\nabla f_{k}(H(\lambda,\theta)),\dfrac{\partial^{2}H}
{\partial\lambda\partial\theta}\rangle=0$. Finally
\[
\dfrac{\partial^{2}H}{\partial\theta^{2}}=\big(\mu_{-}(\lambda)\rho_{c_{-}
}(\theta)+\mu_{0}(\lambda)\rho_{c_{0}}(\theta)+\mu_{+}(\lambda)\rho_{c_{+}
}(\theta)\big)(-n(\theta))+\big(\cdots \big)\tau(\theta),
\]
hence the quantity 
$$
\displaystyle\langle\nabla
f_{k}(H(\lambda,\theta
)),\dfrac{\partial^{2}H}{\partial\theta^{2}}\rangle=-\big(\mu_{-}(\lambda
)\rho_{c_{-}}(\theta)+\mu_{0}(\lambda)\rho_{c_{0}}(\theta)+\mu_{+}
(\lambda)\rho_{c_{+}}(\theta)\big)\Vert\nabla
f_{k}(H(\lambda,\theta))\Vert
$$ 
is negative since all the $\mu$ and
$\rho$ are nonnegative. Hence $Q_{f_{k}}$ is nonnegative and the
function $f_{k}$ is convex.

\medskip

\noindent\emph{$C^{2}$ smoothing}. For
$\lambda\in(\min\tilde{f},\lambda_{0} ]$, define
\[
h(\lambda)=(\lambda-\min\tilde{f})(1+\max_{[\lambda\leq\tilde{f}\leq
\lambda_{0}]}\Vert Q_{\tilde{f}}\Vert)^{-1}.
\]
Since $\tilde{f}$ is $C^{\infty}$ in $[\min\tilde{f}<\tilde{f}]$,
$h$ is a continuous, positive, increasing function. Then there
exists $\psi\in C^{\infty}(\mathbb{R},\mathbb{R}_{+})$ which
vanishes on $(-\infty,\min \tilde{f}]$, increases on $(0,+\infty)$ and for
$\lambda\in(\min\tilde{f},\lambda_{0}]$, $0<\psi(\lambda)\leq
h(\lambda)$ (see Lemma~\ref{minoration}). Let $g$ be the primitive
of $\psi$ with $g(\min\tilde{f})=0$. The function $g$ is a strictly increasing
convex $C^{\infty}$-function on $[\min\tilde{f},+\infty)$. The
function $\bar{f}=g\circ\tilde{f}$ is therefore a $C^{1}$ convex
function. Moreover $\bar{f}$ is $C^{\infty}$ at each point outside
the boundary of $\mbox{\rm argmin\,}f$. For $x\in\mbox{\rm
argmin\,}f$, we have
\[
\frac{\nabla\bar{f}(x+h)-\nabla\bar{f}(x)}{\Vert
h\Vert}=\frac{g^{\prime
}(\tilde{f}(x+h))\nabla\tilde{f}(x+h)}{\Vert
h\Vert}=\frac{g^{\prime} (\tilde{f}(x)+o(\Vert h\Vert))o(1)}{\Vert
h\Vert}=\frac{o(\Vert h\Vert)}{\Vert h\Vert}=o(1).
\]
Thus $Q_{\bar{f}}(x)=0$. On the other hand
\begin{align*}
\Vert Q_{\bar{f}}(x+h)\Vert &  \leq
g^{\prime}(\tilde{f}(x+h))\Vert
Q_{\tilde{f}}(x+h)\Vert+g^{\prime\prime}(\tilde{f}(x+h))\Vert\nabla\tilde
{f}(x+h)\Vert^{2}\\ &  \leq h(\tilde{f}(x+h))\Vert
Q_{\tilde{f}}(x+h)\Vert+o(1)\\ &  \leq(f(x+h)-f(x))+o(1)=o(1).
\end{align*}
Thus $Q_{\bar{f}}$ is continuous at $x$ and thus $\bar{f}$ is
$C^{2}$.

\medskip

\noindent\textit{Conclusion}\label{conc}. Let us prove finally that $\bar{f}$
does not satisfy the K\L--inequality. Towards a contradiction,
let us assume that there exist $R>\inf\bar{f}=\min\bar{f}$, a
continuous function
$\varphi:[\min\bar{f},R)\rightarrow\mathbb{R}_{+}$ which satisfies
$\varphi(\min\bar{f})=0$, $\varphi$ is $C^{1}$ on
$(\min\bar{f},R)$ with $\varphi^{\prime}>0$, such that we have
\[
\Vert\nabla(\varphi\circ\bar{f})(x)\Vert\geq1,\;\forall
x\in\lbrack\min f<f<R].
\]
Applying Theorem~\ref{global} [(i)$\Leftrightarrow$(vi)], we
obtain
\[
\Dist([\bar{f}\leq g(\lambda_{k})],[\bar{f}\leq
g(\lambda_{k+1})])\leq
\varphi(g(\lambda_{k}))-\varphi(g(\lambda_{k+1})).
\]
and, as a consequence, \newline $\displaystyle\sum_{k=0}^{+\infty}
\Dist([\tilde{f}\leq\lambda_{k}],[\tilde{f}\leq\lambda_{k+1}])=\sum
_{k=0}^{+\infty}\Dist([\bar{f}\leq g(\lambda_{k})],[\bar{f}\leq
g(\lambda _{k+1})])\leq\varphi(g(\lambda_{0}))$. This contradicts
the fact that $\sum\Dist(T_{k},T_{k+1})=+\infty$.$\hfill\Box$

\subsection{Asymptotic equivalence for discrete and continuous dynamics}

\label{explicit}

In this part we assume that $f:H\rightarrow\mathbb{R}$ is a
$C^{1,1}$ \emph{convex} function, that is, continuously differentiable with
gradient $\nabla f$ Lipschitz continuous. Let $L$ be a Lipschitz constant of
$\nabla f$.\smallskip

Fix $\beta>0$ and $x\in\mathbb{R}^{n}$ and consider any sequence $\{Y_{x}
^{k}\}$ satisfying
\begin{equation}
\left\{
\begin{array}
[c]{l} \beta\,||\nabla
f(Y_{x}^{k})||\;||Y_{x}^{k+1}-Y_{x}^{k}||\,\leq\,f(Y_{x}
^{k})-f(Y_{x}^{k+1}),\quad k=1,2,\ldots\\
\\
Y_{x}^{0}=x
\end{array}
\right.  \label{expgrad}
\end{equation}
This condition has been considered in \cite{Absil} for nonconvex functions
defined in finite-dimensional spaces. It is easily seen that (\ref{expgrad})
is a descent sequence, that is, $f(Y_{x}^{k})\geq f(Y_{x}^{k+1})$, which
implies in particular that $\{f(Y_{x}^{k})\}$ converges as $k$ goes to infinity.

\smallskip

Condition (\ref{expgrad}) is fulfilled by several explicit gradient--like
methods, including trust region methods, line--search gradient methods and
some Riemannian variants; see \cite{Absil} for examples and references.

\smallskip

The following theorem establishes connections between length boundedness
properties of continuous gradient methods and length boundedness of discrete
gradient iterations.

\begin{theorem}
[discrete vs continuous]\label{T:exp} Let $f$ be a $C^{1,1}$ convex function
with compact sublevel sets such that $\min\,f=0$. Let us denote by $L$ a
Lipschitz constant of $\nabla f$. Then the following statements are equivalent:\smallskip

\noindent(i) \textbf{[Kurdyka-\L ojasiewicz inequality]} There exist $r_{0}>0$
and $\varphi\in\mathcal{K}(0,r_{0})$ such that
\begin{equation}
||\nabla(\varphi\circ f)(x)||\geq1,\qquad\text{for all }
x\in\lbrack0<f\leq r_{0}]. \label{KL}
\end{equation}

\noindent(ii) \textbf{[Length boundedness of piecewise gradient iterates]} For
all $\beta>0$ and all $R>0$, there exists $\mathcal{L}(\beta)>0$ such that for
any sequence of gradient iterates of the form
\[
Y_{x_{0}}^{0},\,Y_{x_{0}}^{1},\ldots,\,Y_{x_{0}}^{k_{0}},\,Y_{x_{1}}
^{0},\,\ldots Y_{x_{1}}^{k_{1}},\,\ldots
\]
with $f(x_{0})<R$, $f(Y_{x_{i+1}}^0)=f(x_{i+1})\leq f(Y_{x_{i}}^{k_{i}})\ $ and $\{Y_{x_{i}}
^{j}:j=0,\ldots,k_{i}\}$ satisfying (\ref{expgrad}) for all $i\in\mathbb{N}$
we have
\[
\sum_{i=0}^{+\infty}\sum_{l=0}^{k_{i}}\,||Y_{x_{i}}^{l+1}-Y_{x_{i}}
^{l}||\,\leq\,\mathcal{L}(\beta).
\]

\noindent(iii) {\textbf{[Length boundedness of piecewise gradient curves]}}
For every $R>0$ there exists $\mathcal{L}>0$ such that
\[
\mbox{\rm length}\;(\gamma)\leq\mathcal{L},
\]
for all piecewise subgradient curves $\gamma:[0,+\infty)\rightarrow H$ with
$f(\gamma(0))<R$.
\end{theorem}

\noindent\textbf{Proof}. Let us first prove that (i)$\Rightarrow$(ii). By
Theorem~\ref{convex}[(i)$\Rightarrow$(ii)] (subgradient inequality -- convex
case) we may assume that $\varphi$ is concave, defined on $(0,+\infty)$ and
(\ref{KL}) holds for all $x\in\lbrack0<f]$. We now proceed in the spirit of
\cite{Absil}. Let $\beta>0$, $x\in\lbrack0<f]$ and let $Y_{x}^{0}
,\ldots,Y_{x}^{k}$ be a (finite) sequence of gradient--type iterations that
satisfies (\ref{expgrad}). For simplicity we set $Y_{x}^{j}=Y^{j}$ for all
$j\in\{0,\ldots,k\}$, so that
\[
f(Y^{j})-f(Y^{j+1})\,\geq\,\beta\,||\nabla f(Y^{j})||\;||Y^{j+1}-Y^{j}||.
\]
Multiplying both parts with $\varphi^{\prime}(f(Y^{j}))$ and applying (i) we
get
\[
\varphi^{\prime}(f(Y^{j}))[f(Y^{j})-f(Y^{j+1})]\,\geq\,\beta\,||Y^{j+1}
-Y^{j}||.
\]
Since $\varphi$ is concave we have
\[
\varphi(f(Y^{j+1}))\leq\varphi(f(Y^{j}))+\varphi^{\prime}(f(Y^{j}
))\,[f(Y^{j+1})-f(Y^{j})],
\]
and therefore
\[
\varphi(f(Y^{j}))-\varphi(f(Y^{j+1}))\,\geq\,\beta\,||Y^{j+1}-Y^{j}||.
\]
Adding the above inequalities for $j=0,\ldots,k$ we obtain
\begin{equation}
\varphi(f(Y^{0}))-\varphi(f(Y^{k}))\,\geq\,\beta\,\sum_{j=0}^{k}
||Y^{j+1}-Y^{j}||. \label{sum}
\end{equation}
Let us now consider a sequence of the form $\{Y_{x_{0}}^{0},\,Y_{x_{0}}
^{1},\ldots,\,Y_{x_{0}}^{k_{0}},\,Y_{x_{1}}^{0},\,\ldots Y_{x_{1}}^{k_{1}
},\,\ldots\}$ as in (ii). Then applying (\ref{sum}) to each subsequence
$\{Y_{x_{i}}^{j},\;j=0,\ldots,k_{i}\}$ we deduce
\[
\sum_{i=0}^{+\infty}\sum_{l=0}^{k_{i}}\,||Y_{x_{i}}^{l+1}-Y_{x_{i}}
^{l}||\,<\,\frac{1}{\beta}\varphi(f(Y_{x_{0}}^{0}))\,\leq\,\frac{1}{\beta
}\varphi(R),
\]
which proves the assertion.

\smallskip

The equivalence (i)$\Longleftrightarrow$(iii) follows from
Theorem~\ref{Theorem_main} and Theorem~\ref{convex}. To complete the proof it
suffices to establish that (ii) implies the assertion (iv) of
Theorem~\ref{Theorem_main} (valley selection of finite length) (in fact we
prove (iv') with $R=2$). So let us assume that (ii) holds and let $r_{0} >m$.
We aim to construct a piecewise absolutely continuous curve $\theta
:(0,r_{0}]\rightarrow\mathbb{R}^{n}$ of finite length that satisfies
\[
\theta(r)\in\mathcal{V}_{2}(r):=\left\{  x\in[f=r]:\quad||\nabla
f(x)||\,\leq\,2\,\underset{y\in[f=r]}{\inf}\,||\nabla f(y)||\right\}
,\;\forall r\in(0,r_{0}].
\]
We shall use the explicit gradient method described in
Subsection~\ref{Explicit}. Let $x_{0}\in\mathcal{V}_{2}(r_{0})$ be such that
\[
||\nabla f(x_{0})||\,\leq\frac{3}{2}\underset{y\in f^{-1}(r_{0})}{\inf
}\,||\nabla f(y)||,
\]
and consider the $C^{1}$ curve
\[
\lbrack0,\frac{1}{3L})\ni t\longmapsto x_{0}(t):=x_{0}-t\nabla f(x_{0}).
\]
Set $t_{0}=\sup A_{0}$ where
\[
A_{0}:=\left\{  \;t\in(0,\frac{1}{3L}):\;\left[
\begin{array}
[c]{c} f\circ x_{0}\text{ }\quad\text{strictly decreasing on
}[0,t],\\
x_{0}(\tau)\in\mathcal{V}_{2}(f(x_{0}(\tau))\quad\text{for
}\tau\in \lbrack0,t].
\end{array}
\right.  \quad\right\}  .
\]
Clearly $A_{0}$ is nonempty and $0<t_{0}\leq(3L)^{-1}$. Set $r_{1}
=f(x_{0}(t_{0}))<r_{0}$ and take $x_{1} \in\mathcal{V}_{2}(r_{1})$
such that
\[
||\nabla f(x_{1})||\,\leq\,\frac{3}{2}\underset{y\in[f=r_{1}]}{\inf}\,||\nabla
f(y)||.
\]
Proceeding by induction we obtain a sequence
$\{(t_{k},r_{k},x_{k})\}$ where $\{r_{k}\}\subset\lbrack0,r_{0}]$
is strictly decreasing, $x_{n} (t):=x_{n}-t\nabla f(x_{n})$ with
$f(x_{n})=r_{n}$ and
\[
||\nabla f(x_{n})||\,\leq\,\frac{3}{2}\underset{y\in[f=r_{n}]}{\inf}\,||\nabla
f(y)||.
\]
Let us denote by $r_{\infty}$ the limit of $\{r_{k}\}$ and let us assume,
towards a contradiction, that $r_{\infty}>0$. Set
\[
s(r):=\underset{x\in f^{-1}(r)}{\inf}||\partial f(x)||_{-}\quad\text{and}\quad
s_{\infty}=\,\underset{n\rightarrow\infty}{\liminf}\,s(r_{n})\,=\,\underset
{n\rightarrow\infty}{\lim}s(r_{n})
\]
(note that convexity of $f$ guarantees that $s(r_{1})\leq s(r_{2})$ whenever
$r_{1}\leq r_{2}$) and observe that $r_{\infty}>0$ implies that $s_{\infty}>0$
(use the compactness of the sublevel set $[f\leq r_{0}]$). Let $n_{0}
\in\mathbb{N}$ be such that $s(r_{n})\leq\frac{5}{4}s_{\infty}$ for all $n\geq
n_{0}$. For $n\geq n_{0}$ and $t\in\lbrack0,t_{n})$,
Proposition~\ref{estimation} (Annex) yields
\[
||\nabla f(x_{n}(t))||\,\leq\,(Lt+1)\,||\nabla f(x_{n})||,
\]
which implies
\[
||\nabla f(x_{n}(t))||\,\leq\,(Lt+1)\,||\nabla f(x_{n})||\,\leq\;\frac{3}
{2}(Lt+1)s(r_{n})\;\leq\;\frac{15}{8}(Lt+1)s_{\infty}.
\]
A sufficient condition to have $x_{n}(t)\in\mathcal{V}_{2}(f(x_{n}(t)))$ is
therefore
\begin{equation}
\frac{15}{8}(Lt+1)s_{\infty}\,\leq\,2s_{\infty}\;\Longleftrightarrow\;0\leq
t\leq(15L)^{-1}. \label{vall}
\end{equation}
Similarly we can estimate the rate of decrease of $f(x_{n}(t))$. Since
\[
\frac{d}{dt}f(x_{n}(t))=-\langle\nabla f(x_{n}),\nabla f(x_{n}(t))\rangle,
\]
the condition $\frac{d}{dt}f(x_{n}(t))<0$ is satisfied whenever
\[
||\nabla f(x_{n})||^{2}>||\nabla f(x_{n})||\;||\nabla f(x_{n}(t))-\nabla
f(x_{n})||
\]
But since $\nabla f$ is Lipschitz continuous, $\|\nabla f(x_{n}(t))-\nabla
f(x_{n})\|\le Lt\|\nabla f(x_{n})\|$. Thus the condition is satisfied if
\[
\|\nabla f(x_{n})\|^{2}>Lt\|\nabla f(x_{n})\|^{2}
\]
This last inequality is equivalent to $t<L^{-1},$ which implies in particular
that for all $n\in\mathbb{N}$ such that $s(r_{n})\leq\frac{5}{4}s_{\infty}$,
we have
\[
t_{n}\geq(15L)^{-1}.
\]
In this case Proposition~\ref{estimation}\,(Annex) yields
\[
f(x_{n}(t_{n}))\leq f(x_{n})+(\frac{Lt_{n}^{2}}{2}-t_{n})||\nabla
f(x_{n})||^{2}\leq r_{n}+\frac{9}{4}(\frac{Lt_{n}^{2}}{2}-t_{n})s(r_{n}
)^{2}\leq r_{n}+\frac{9}{4}(\frac{Lt_{n}^{2}}{2}-t_{n})(\frac{5}{4}s_{\infty
})^{2}.
\]
Thus in order to have $f(x_{n}(t))<r_{\infty}$, it suffices to require
\[
t_{n}-\frac{Lt_{n}^{2}}{2}>\frac{64}{225}\left(  \frac{r_{n}-r_{\infty}
}{s_{\infty}}^{2}\right)  .
\]
Using the fact that $(3L)^{-1}\geq t_{n}\geq(15L)^{-1}$, we see that
\[
t_{n}-\frac{Lt_{n}^{2}}{2}\geq(15L)^{-1}-(18L)^{-1}=(90L)^{-1}.
\]
Since $({r_{\infty}-r_{n}})/{s_{\infty}}$ tends to zero, we have that
$f(x_{n}(t_{n}))<r_{\infty}$ for $n$ sufficiently large, which is a contradiction.

\smallskip

We thus conclude that $\{r_{k}\}\rightarrow r_{\infty}=0$ and $(0,r_{0}
]=\cup_{n}(r_{n+1},r_{n}].$ We define $\theta:(0,r_{0}]\rightarrow H$ as
follows: $\theta(r):=x_{n}([f\circ x_{n}]^{-1}(r))$ whenever $r\in
(r_{n+1},r_{n}]$. Clearly $\theta$ defines a piecewise absolutely continuous
curve. To see that $\theta$ has finite length it suffices to observe that the
sequence $\{x_{n}\}_{n}$ is a sequence of gradient iterates that satisfies
(\ref{expgrad}). Using Remark~\ref{descentcond} and the fact that the
step--sizes in the construction of the $x_{n}$'s do not exceed $(3L)^{-1}$ we
infer that
\[
\frac{5}{6}\;||x_{n+1}-x_{n}||\;||\nabla f(x_{n})||\leq f(x_{n})-f(x_{n+1}).
\]
Hence the curve $\theta$ has a finite length. This completes the
proof.$\hfill\Box$

\bigskip

\begin{remark}
The assumption that $f$ is convex has been used to apply Theorem~\ref{convex}
(cf. concavity of $\varphi$ which seems to be crucial for the proof of
implication (i)$\Rightarrow$(ii)) and to assert that $f(Y_{0}^{k}
)\rightarrow\inf f$. These are the reasons for which Theorem~\ref{T:exp} is
not stated for general semiconvex functions (in a local version). It would
therefore be interesting to figure out under which type of conditions (other
than convexity or o-minimality of $f$) the function $\varphi$ of (\ref{KL})
can be taken concave.
\end{remark}

\section{Annex}

\label{Annex}In this Annex section we give several technical results which are
needed in the text.

\subsection{Technical results}

\begin{proposition}
[closed graph of the subdifferential]\label{graph-closed} Let
$f:H\rightarrow \mathbb{R}\cup\{+\infty\}$ be a lower
semicontinuous semiconvex function. Let $\{x_{k}\}$ and
$\{p_{k}\}$ be two sequences in $H$ such that $p_{k} \in\partial
f(x_{k})$, $x_{k}$ converges strongly to $x$ and $p_{k}$ converges
weakly to $p$. Then as $k\rightarrow+\infty$ we obtain
\[
\left\{
\begin{array}
[c]{l} f(x_{k})\rightarrow f(x)\\ p\in\partial f(x)
\end{array}
\right.
\]
\end{proposition}

\noindent\textbf{Proof.} This is a standard property. For a proof (in the more
general setting of primer--lower--nice functions) we refer the reader to
\cite{Thibault}.$\hfill\Box$

\bigskip

\begin{proposition}
[slope functions and semicontinuity]\label{slopes} Let
$f:H\to\mathbb{R} \cup\{+\infty\}$ be a lower semicontinuous
semiconvex function.

\smallskip

\noindent(i) The extended-real-valued function
\begin{equation}
H\ni x\longmapsto||\partial f(x)||_{-}:=\,\inf_{p\in\partial
f(x)}\,||p|| \tag{slope at $x$}
\end{equation}
is lower semicontinuous.\smallskip

\noindent(ii) Take $r_{0}\in\mathbb{R}$ and let $D$ be a nonempty compact
subset of $[f\leq r_{0}]$. Then the function
\begin{equation}
(-\infty,r_{0}]\ni r\longmapsto
s_{D}(r):=\,\underset{x\in[f=r]\cap D}{\inf }\,||\partial
f(x)||_{-} \tag{minimal slope of the {\it r} level-line}
\end{equation}
is lower semicontinuous. \smallskip

\noindent(iii) Assume that \eqref{crit} and \eqref{compacite-ss-niv} hold for
some $\bar{r}, \bar{\epsilon} >0.$ If $0<r_{1}\leq r_{2}\leq\bar{r},$ then
there exists $\eta_{r_{1},r_{2}} >0$ such that
\begin{align*}
\mathop{\rm inf}_{x\in[r_{1}\leq f\leq r_{2}]\cap\bar{B}(\bar{r},
\bar{\epsilon})} ||\partial f(x)||_{-} \geq\eta_{r_{1},r_{2}} >0.
\end{align*}
\end{proposition}

\noindent\textbf{Proof.} (ii) Take $r\in(-\infty,r_{0}]$ and let
$\{r_{k}\}\subset(-\infty,r_{0}]$ be a sequence such that
$r_{k}\rightarrow r$
and $\liminf_{k}s_{D}(r_{k})<+\infty$. Fix $\eta>0$ and let $(x_{k},p_{k}%
)\in\mbox{graph\,}\partial f$ be such that $f(x_{k})=r_{k}$, $p_{k}\in\partial
f(x_{k})$ and $||p_{k}||<s_{D}(r_{k})-\eta$. Using a standard compactness
argument together with the fact that $\liminf_{k}s_{D}(r_{k})<+\infty$ we can
assume, with no loss of generality, that $x_{k}$ converges (strongly) to $x\in
D$ and that $p_{k}$ converges weakly to $p$. Using
Proposition~\ref{graph-closed}, we obtain that $(x,p)\in\mbox{graph\,}\partial
f$ and $f(x)=r.$ The conclusion follows from the (weak) lower semicontinuity
of the norm. Indeed
\[
\liminf_{k\rightarrow+\infty}s_{D}(r_{k})-\eta\geq\liminf_{k\rightarrow
+\infty}||p_{k}||\geq||p||\geq s_{D}(r).
\]
The proof of (i) and (iii) involve similar arguments. $\hfill\Box$

\begin{lemma}
[strong slope]\label{Lemma_slope} Let $f$ be a proper lower semicontinuous
semiconvex function. Then for all $x$ in $H$
\[
||\partial f(x)||_{-}=|\nabla f|(x).
\]
\end{lemma}

\noindent\textbf{Proof}. Let $x\in H$ and $p=\partial^{0} f(x)$ the projection
of 0 on $\partial f(x).$ By \eqref{croiss-grad-semiconv}, for any $y\in H,$ we
have
\begin{align*}
\frac{(f(x)-f(y))^{+}}{||y-x||} \leq( -\langle p,
\frac{y-x}{||y-x||}
\rangle+\alpha||y-x||^{2})^{+}\leq(||p||+\alpha||y-x||^{2})^{+}.
\end{align*}
By taking the limsup as $y\to x,$ we get $|\nabla f|(x)\leq||p||=||\partial
f(x)||_{-}.$ To prove the opposite inequality, we consider the subgradient
trajectory $\chi_{x}.$ If $x$ is a critical point of $f,$ then $0=||\partial
f(x)||_{-}\geq|\nabla f|(x).$ Otherwise, $\chi_{x}(t)\not = x$ for all $t>0.$
By Theorem~\ref{Theorem_thibault}(iv), we have
\begin{align*}
\frac{(f(x)-f(\chi_{x}(t)))^{+}}{||x-\chi_{x}(t)||}\geq\frac{1}{||x-\chi
_{x}(t)||} \int_{0}^{t} ||\partial f(\chi_{x}(\tau))||_{-}^{2} d\tau.
\end{align*}
Taking the limsup as $t\downarrow0$ and using the continuity of the semiflow
and Theorem \ref{Theorem_thibault}(ii),(iii) we obtain the desired result.
\hfill$\Box$

\begin{lemma}
[chain rules]\label{chain} Let $f:H\rightarrow\mathbb{R}\cup\{+\infty\}$ be a
extended-real-valued function.\newline (i) Let $\varphi:(0,1)\rightarrow
\mathbb{R}$ be a $C^{1}$ function. Then
\[
\partial(\varphi\,of)(x)=\varphi^{\prime}(f(x))\partial f(x),\mbox{ for
all }x\in\lbrack0<f<1].
\]
(ii) Let $\gamma:(0,1)\rightarrow H$ be a $C^{1}$ curve. For all $t\in(0,1),$
we have
\[
\partial(f\circ\gamma)(t)\supset\{\langle\dot{\gamma}(t),p(t)\rangle
:p(t)\in\partial f(\gamma(t))\}.
\]
\end{lemma}

\noindent\textbf{Proof} For the proof see \cite{Rock98} for example.
$\hfill\Box$

\begin{lemma}
[continuous integrable majorant]\label{LemmaReg} Let $u:(0,r_{0}
]\rightarrow\mathbb{R}_{+}$ be an upper semicontinuous function such that
$u\in L^{1}(0,r_{0})$. Then there exists a \emph{\ continuous} function
$w:(0,r_{0}]\rightarrow\mathbb{R}_{+}$ such that $w\geq u$ and $w\in
L^{1}(0,r_{0})$. If moreover $u$ is assumed to be nonincreasing, $w$ can be
chosen to be decreasing.
\end{lemma}

\noindent\textbf{Proof} With no loss of generality we assume $r_{0}=1$.
Replacing if necessary $u(\cdot)$ by the function $u(\cdot)+1$ we may also
assume that $u\geq1$. Let $a_{k}>0$ be a strictly decreasing sequence such
that $a_{0}=1$ and $(0,1]=\cup_{k\in\mathbb{N}}[a_{k+1},a_{k}]$. Let us assume
that there exists a sequence of continuous functions $w_{k}:[a_{k+1}
,a_{k}]\rightarrow\mathbb{R}$ such that $w_{k}\geq u$ on $[a_{k+1},a_{k}]$ and
$\int_{a_{k+1}}^{a_{k}}w_{k}\leq\int_{a_{k+1}}^{a_{k}}u+\frac{1}{(k+1)^{2}}$.
To establish the existence of $w$, we proceed by induction on $k$. Fix
$k\geq1$ and assume that $w$ is defined on $[a_{k},1]$ with $w\geq u$, $w$
continuous and
\[
\int_{a_{k}}^{1}w\leq\int_{a_{k}}^{1}u+\sum_{i=1}^{k}\frac{2}{i^{2}}.
\]
There is no loss of generality to assume $w_{k}(a_{k})\leq w(a_{k})$ (the case
$w_{k}(a_{k})>w(a_{k})$ can be treated analogously). Let us define
\[
0< \epsilon_{k}=\frac{w_{k}(a_{k})(a_{k}-a_{k+1})}{(k+1)^{2}\;w(a_{k}
)\;\max_{[a_{k+1},a_{k}]}w_{k}}< a_{k}-a_{k+1},
\]
and let us consider the functions
\[
\lambda_{k}:[a_{k}-\epsilon_{k},a_{k}]\rightarrow\lbrack1,\frac{w(a_{k}
)}{w_{k}(a_{k})}]
\]
defined by
\[
\lambda_{k}(r)=\frac{1}{\epsilon_{k}}\left(  (a_{k}-r)+(r-(a_{k}-\epsilon
_{k}))\frac{w(a_{k})}{w_{k}(a_{k})}\right)  .
\]
The function $w$ can be now extended to $[a_{k+1},1]$ by setting
\[
w(r)=\left\{
\begin{array}
[c]{ll} w_{k}(r), & \text{if }r\in\lbrack
a_{k+1},a_{k}-\epsilon_{k}),\\ \lambda_{k}(r)w_{k}(r), & \text{if
}r\in\lbrack a_{k}-\epsilon_{k},a_{k}]\\ w(r), & \text{if
}r\in(a_{k},1].
\end{array}
\right.
\]
It is easily seen that the function $w$ is continuous (by definition of
$\lambda_{k}$), it satisfies $w\geq u$ on $[a_{k+1},a_{k}]$ (thus on
$(a_{k+1},1]$) and moreover
\begin{align*}
\int_{a_{k+1}}^{1}w  &  =\int_{a_{k+1}}^{a_{k}-\epsilon_{k}}w_{k}+\int
_{a_{k}-\epsilon_{k}}^{a_{k}}\lambda_{k}w_{k}+\int_{a_{k}}^{1}w\\
&  \leq\int_{a_{k+1}}^{a_{k}}u+\frac{1}{(k+1)^{2}}+\epsilon_{k}\frac{w(a_{k}
)}{w_{k}(a_{k})}\max_{[a_{k+1},a_{k}]}w_{k}+\int_{a_{k}}^{1}u+\sum_{i=1}
^{k}\frac{2}{i^{2}}.,\\
&  \leq\int_{a_{k+1}}^{1}u+\frac{2}{(k+1)^{2}}+\sum_{i=1}^{k}\frac{2}{i^{2}}.
\end{align*}
This proves the existence of a continuous function $w$ that satisfies the
required properties.

\smallskip

To complete the proof it suffices to prove the existence of such a sequence
$\{w_{k}\}$. To this end, fix $k\in\mathbb{N}^{\ast}$ and set
\[
u^{\epsilon}(r)=\sup_{\rho\in\lbrack a_{k+1},a_{k}]}\{ u(\rho)-\frac
{||r-\rho||^{2}}{2\epsilon}\}.
\]
It is easily seen that $u^{\epsilon}$ is continuous, $u(r)\leq
u^{\epsilon }(r)\leq\max_{\rho\in\lbrack
a_{k+1},a_{k}]}u:=M_{k}<+\infty$ and
$\lim_{\epsilon\rightarrow0}u^{\epsilon}(r)=u(r)$ for all
$r\in\lbrack a_{k+1},a_{k}]$ (see \cite{Rock98}, for example).
Note that the upper semicontinuity of $u$ on the compact set
$[a_{k+1},a_{k}]$ guarantees that $M_{k}$ is finite. Applying the
Lebesgue domination convergence theorem we conclude that
$u^{\epsilon}$ converges to $u$ in the norm topology of
$L^{1}(a_{k+1},a_{k})$. Thus there exists $\epsilon_{0}>0$ such
that
\[
\int_{\lbrack a_{k+1},a_{k}]}u^{\epsilon_{0}}\leq\int_{\lbrack a_{k+1},a_{k}
]}u+\frac{1}{(k+1)^{2}}.
\]
Thus the function $w_{k}:=u^{\epsilon_{0}}$ satisfies the requirements stated
above. This completes the proof of the first part of the statement. The case
where $u$ is assumed decreasing, can be treated with similar (and occasionally
simpler) arguments. $\hfill\Box$

\begin{lemma}
\label{minoration} Let $h\in C^{0}((0,r_{0}],\mathbb{R}_{+}^{*})$
be an increasing function, then there exists a function $\psi\in
C^{\infty }(\mathbb{R},\mathbb{R}_{+})$ such that $\psi=0$ on
$\mathbb{R}_{-}$, $0<\psi(s)\le h(s)$ for all $s\in(0,r_{0})$, and
$\psi$ is increasing on $(0,r_{0})$.
\end{lemma}

\noindent\textbf{Proof.} Let us extend the definition of $h$ by
$0$ on $\mathbb{R}_{-}$ and $h(r_{0})$ for $s>r_{0}$. Consider
$\phi\in C^{\infty }(\mathbb{R},\mathbb{R}_{+})$ with $[0,1]$ as
support and $\int_{\mathbb{R} }\phi=1$. Then we define $\psi$ by
$\psi=\phi*h$; \textit{i.e.} $\psi
(s)=\int_{\mathbb{R}}\phi(t)h(s-t)dt$. It is then straightforward
to verify that $\psi$ satisfies the expected properties.
$\hfill\Box$

\begin{proposition}
[Piecewise absolutely continuous selections]\label{selection} Let $r_{0}>0$
and $\mathcal{V}:(0,r_{0}]\rightrightarrows H$ be a set-valued mapping with
nonempty values. Assume that for each $r\in(0,r_{0}]$ there exists
$\epsilon_{r}\in(0,r)$ and an absolutely continuous curve $\theta
_{r}:(r-\epsilon_{r},r]\rightarrow H$ such that
\[
\theta_{r}(s)\in\mathcal{V}(s)\mbox{ for all }s\mbox{ in }(r-\epsilon_{r},r].
\]
Then there exist a countable partition $\{I_{n}\}_{n\in\mathbb{N}}$ of
$(0,r_{0}]$ into intervals $I_{n}$ of nonempty interior and a selection
$\theta:(0,r_{0}]\rightarrow\mathbb{R}^{n}$ of $\mathcal{V}$ such that
$\theta$ is absolutely continuous on each $I_{n}$.
\end{proposition}

\noindent\textbf{Proof.} Let $\Omega$ be the set of couples
$(\alpha
:I_{\alpha}\subset(0,r_{0}]\rightarrow\mathbb{R}^{n},\;\{I_{\alpha,j}\}_{j\in
J_{\alpha}})$ where $\{I_{\alpha,j}\}_{j\in J_{\alpha}}$ is a
countable partition of $I_{\alpha}$ into (disjoint) intervals
$I_{\alpha,j},\;{j\in J_{\alpha}}$ with nonempty interior such
that:\smallskip

(a) for each $j\in J_{\alpha}$, $\alpha$ is absolutely continuous on
$I_{\alpha,j}$,\smallskip

(b) for each $r\in I_{\alpha}$, $\alpha(r)\in\mathcal{V}(r)$.

\smallskip

\noindent We define a partial order $\preccurlyeq$ on $\Omega$ by
\[
\begin{array}
[c]{c}
\alpha_{1}\preccurlyeq\alpha_{2}\quad\Leftrightarrow\quad\forall
j\in J_{\alpha_{1}},\;\exists k\in
J_{\alpha_{2}},\;I_{\alpha_{1},j}\subset I_{\alpha_{2},k}\mbox{
and }\alpha_{1}(r)=\alpha_{2}(r)\mbox{ for all }r\in
I_{\alpha_{1}}.
\end{array}
\]
Note that $(\Omega,\preccurlyeq)$ is nonempty partially ordered. Let us check
that each totally ordered subset of $\Omega$ has an upper bound in $\Omega$.
To this end, let
\[
\omega=\{(\alpha_{l},\{I_{\alpha_{l},j}\}_{j\in J_{\alpha_{l}}})\}_{l\in
\mathcal{L}}
\]
be a totally ordered subset of $\Omega$. For each $r\in\cup_{l\in\mathcal{L}
}I_{\alpha_{l}}$ define $\alpha(r)$ by
\[
\alpha(r):=\alpha_{l}(r),
\]
whenever $r\in I_{l},$ and set $I_{\alpha}=\cup_{l\in L}I_{\alpha_{l}}$. Since
$\omega$ is totally ordered, the mapping $\alpha:I_{\alpha}\rightarrow
\mathbb{R}^{n}$ is well defined and (b) is clearly satisfied. For $l\in L$ and
$j\in J_{l}$, set $J_{l}:=J_{\alpha_{l}}$, $I_{\alpha_{l},j}=I_{l,j}$ and
$D:=\{(m,k):m\in L,\;k\in J_{m}\}$. For each $(l,j)\in D$, let us define
\begin{equation}
M_{l,j}:=\bigcup_{(m,k)\in D,\;I_{l,j}\subset I_{m,k}}I_{m,k}.
\label{ens}
\end{equation}
Observe that $I_{\alpha}=\cup_{(l,j)\in D}M_{l,j}$ and that each $M_{l,j}$ is
an interval with nonempty interior.

\medskip

Let us prove that for all $(l,j),(l^{\prime},j^{\prime})\in D$, we have either
$M_{l^{\prime},j^{\prime}}=M_{l,j}$ or $M_{l^{\prime},j^{\prime}}\cap
M_{l,j}=\emptyset$. In order to establish this result, let us beforehand show
that for all $(l,j)$, $(l^{\prime},j^{\prime})$ in $D$ such that $I_{l,j}\cap
I_{l^{\prime},j^{\prime}}\neq\emptyset$, we have $M_{l,j}=M_{l^{\prime
},j^{\prime}}$. Indeed, since $\omega$ is totally ordered, we have for
instance $I_{l^{\prime},j^{\prime}}\subset I_{l,j}$ and so $M_{l,j}
\subset\ M_{l^{\prime},j^{\prime}}$. Conversely, take $(m,k)\in D$ such that
$I_{m,k}\supset I_{l^{\prime},j^{\prime}}$. Since $I_{m,k}\cap I_{l,j}
\neq\emptyset$, we have either $I_{m,k}\subset I_{l,j}$ or $I_{m,k}\supset
I_{l,j}$, in any case we see (\textit{cf.}~definition~(\ref{ens})) that
$I_{m,k}\subset M_{l,j}$ and thus $M_{l^{\prime},j^{\prime}}\subset M_{l,j}$.

\smallskip

If $M_{l,j}\cap M_{l^{\prime},j^{\prime}}\neq\emptyset$, take $r$ in the
intersection, and observe that by definition there exist $(m,k)$ and
$(m^{\prime},k^{\prime})$ in $D$ such that $I_{m,k}\supset I_{l,j}$ with $r\in
I_{m,k}$ and $I_{m^{\prime},k^{\prime}}\supset I_{l^{\prime},j^{\prime}}$ with
$r\in I_{m^{\prime},k^{\prime}}$. Using the previous remark, we obtain that
$M_{m,k}=M_{l,j}$ and $M_{m^{\prime},k^{\prime}}=M_{l^{\prime},j^{\prime}}$.
But since $I_{m,k}\cap I_{m^{\prime},k^{\prime}}\neq\emptyset$, we also have
$M_{m,k}=M_{m^{\prime},k^{\prime}}$ and thus $M_{l,j}=M_{l^{\prime},j^{\prime
}}$.

\smallskip

Let us define an equivalence relation $\simeq$ on $D$ by
\[
(l,j)\simeq(l^{\prime},j^{\prime})\Leftrightarrow M_{l,j}=M_{l^{\prime
},j^{\prime}}.
\]
This equivalence relation defines a partition of $D$ into equivalence classes.
By the axiom of choice we can pick one and only one element in each
equivalence class and this defines a nonempty subset $D^{\prime}$ of $D$. By
construction we have $I_{\alpha}=\cup_{(l,j)\in D^{\prime}}M_{l,j}$ and
$M_{l,j}\cap M_{l^{\prime},j^{\prime}}=\emptyset$ for each $(l,j)\neq
(l^{\prime},j^{\prime})$ in $D^{\prime}$. Besides since each $M_{l,j}$ (for
$\;(l,j)\in D^{\prime}$) has a nonempty interior, we see that $D^{\prime}$ is
a countable set. This shows that $(\alpha,\{M_{l,j},\;(l,j)\in D^{\prime}\})$
is in $\Omega$ with in addition $\alpha\geq\alpha_{l}$ for all $l\in
\mathcal{L}$.

\smallskip

Applying Zorn's lemma to $\Omega$, we obtain the existence of a maximal
element $(\theta:I_{\theta}\rightarrow\mathbb{R}^{n},\{I_{\theta,j},\;j\in
J_{\theta}\})$. Arguing by contradiction, we see immediately that $I_{\theta
}=(0,r_{0}]$.$\hfill\Box$

\subsection{Explicit gradient method}

\label{Explicit} We recall the following useful result

\begin{lemma}
[Descent lemma]\label{descent} Let $f$ be a $C^{1,1}$ function (that is,
$\nabla f$ is $L$-Lipschitz continuous). Then
\[
f(y)\leq f(x)+\langle\nabla f(x),y-x\rangle+\frac{L}{2}||y-x||^{2}.
\]
\end{lemma}

\noindent\textbf{Proof} Set $x(t)=x+t(y-x)$ and notice that
\[
f(y)-f(x)=\int_{0}^{1}\frac{d}{dt}f(x(t))dt=\langle\nabla f(x),y-x\rangle
+\int_{0}^{1}\langle\nabla f(x(t))-\nabla f(x),y-x\rangle dt.
\]
The assertion follows easily.$\hfill\Box$

\bigskip

Given $x\in H$, let us consider the following recursion rule
\begin{equation}
x^{+}:=X(t,x)=x-t\nabla f(x),\;t>0. \label{x+}
\end{equation}
Choosing a starting point $x^{0}$ in $H$, and $\lambda_{k}>0$ a sequence of
step size, the explicit gradient method writes
\[
x^{k+1}=X(\lambda_{k},x^{k}).
\]
A part of the convergence analysis of this method (and some of its variants)
is based on the following elementary results.

\begin{proposition}
\label{estimation} Let $f$ be a $C^{1,1}$ function, $x\in H$, $t\in
\lbrack0,2L^{-1})$ and $x^{+}$ be given by \eqref{x+}. Then\smallskip

(i) $(1-\frac{Lt}{2})\;||x^{+}-x||\;||\nabla f(x)||\,\leq\,f(x)-f(x^{+})$ ;\smallskip

(ii) $||\nabla f(x^{+})||\,\leq\,(Lt+1)\,||\nabla f(x)||.$
\end{proposition}

\noindent\textbf{Proof} Assertion (i) follows directly from
Lemma~\ref{descent} while assertion (ii) is a consequence of the fact that
$\nabla f$ is Lipschitz continuous on $[x,x(t)]$ of constant $L.\hfill\Box$

\begin{remark}
\label{descentcond} Condition (\ref{expgrad}) of Section~\ref{explicit}
corresponds of course to the inequality (i) above.
\end{remark}

\begin{center}
---------------------------------------
\end{center}

\bigskip

\begin{center}
----------------------------------------------------
\end{center}

\newpage

\noindent J\'{e}r\^{o}me BOLTE\smallskip

\noindent UPMC Univ Paris 06 - Equipe Combinatoire et Optimisation (UMR 7090), Case 189\newline
Universit\'{e} Pierre et Marie Curie\newline 4 Place Jussieu, F--75252 Paris
Cedex 05.

\medskip

\noindent INRIA Saclay, CMAP, Ecole Polytechnique, 91128 Palaiseau, France.

\smallskip

\noindent E-mail: \texttt{bolte@math.jussieu.fr}\newline \texttt{http://www.ecp6.jussieu.fr/pageperso/bolte}

\bigskip

\noindent Aris DANIILIDIS\smallskip

\noindent Departament de Matem\`{a}tiques, C1/308\newline Universitat
Aut\`{o}noma de Barcelona\newline E--08193 Bellaterra (Cerdanyola del
Vall\`{e}s), Spain.

\medskip

\noindent Laboratoire de Math\'{e}matiques et Physique Th\'{e}orique \newline 
Universit\'{e} Fran\c{c}ois Rabelais, Tours, France.

\smallskip

\noindent E-mail: \texttt{arisd@mat.uab.es}\newline
\texttt{http://mat.uab.es/\symbol{126}arisd}

\bigskip

\noindent Olivier LEY\smallskip

\noindent Laboratoire de Math\'{e}matiques et Physique Th\'{e}orique (CNRS UMR
6083)\newline 
F\'ed\'eration Denis Poisson\newline
Facult\'{e} des Sciences et Techniques, Universit\'{e} Fran\c
{c}ois Rabelais \newline Parc de Grandmont, F--37200 Tours, France.

\smallskip

\noindent E-mail: \texttt{ley@lmpt.univ-tours.fr}\newline
\texttt{http://www.phys.univ-tours.fr/\symbol{126}ley}

\bigskip

\noindent Laurent MAZET\smallskip

\noindent Universit\'{e} Paris-Est,\newline Laboratoire d'Analyse et
Math\'{e}matiques Appliqu\'{e}es, UMR 8050\newline UFR des Sciences et
Technologie, D\'{e}partement de Math\'{e}matiques\newline 61 avenue du
G\'en\'eral de Gaulle 94010 Cr\'{e}teil cedex, France.

\smallskip

\noindent E-mail: \texttt{laurent.mazet@univ-paris12.fr}\newline \texttt{http://perso-math.univ-mlv.fr/users/mazet.laurent/}
\end{document}